\documentclass[11pt,twoside]{article}
\newcommand{\ifims}[2]{#1} 
\newcommand{\ifAMS}[2]{#1}   
\newcommand{\ifau}[4]{#1}  
\newcommand{\ifbook}[2]{#1}   
\newcommand{\ifapp}[2]{#2}  
\newcommand{\ifshort}[2]{#2}  
\newcommand{\iffourG}[2]{#2}  
\newcommand{\ifsqnorm}[2]{#2}  

\usepackage{tocloft}

\advance\cftsecnumwidth 0.3em\relax
\advance\cftsubsecnumwidth 0.3em\relax

\usepackage{ifthen}
\usepackage[ampersand]{easylist}
\usepackage{amsmath,amssymb,amsthm}
\usepackage[bb=stixtwo]{mathalpha}
\usepackage{mathtools}
\usepackage{natbib}
\usepackage{epsfig,graphicx}
\usepackage{comment}
\usepackage{color}
\usepackage{srcltx}
\usepackage[mathscr]{eucal}
\usepackage[math]{easyeqn}
\usepackage{etoolbox}
\usepackage{hyperref}
\hypersetup{
            colorlinks,
            linkcolor=hookersgreen,
            linktoc=blue,
            citecolor=ultramarine,
            urlcolor=black,
            filecolor=black
            }
\usepackage{nameref}
\usepackage{multirow}

\usepackage{accents}
\usepackage{mathrsfs}
\usepackage{bbm}
\usepackage{algorithm}
\usepackage{algpseudocode}

\ifims{
\textheight=22cm
\textwidth=14.8cm
\topmargin=0pt
\oddsidemargin=1.0cm
\evensidemargin=1.0cm
\linespread{1.3}

}{ 
}

\numberwithin{equation}{section}
\numberwithin{figure}{section}
\newcounter{example}[section]
\numberwithin{example}{section}
\newcounter{remark}[section]
\numberwithin{remark}{section}
\newtheorem{theorem}{Theorem}[section]
\newtheorem{proposition}[theorem]{Proposition}
\newtheorem{lemma}[theorem]{Lemma}

\newtheorem{exmp}[example]{Example}
\newtheorem{rmrk}[remark]{Remark}
\newenvironment{example}{\begin{exmp}\rm}{\end{exmp}}
\newenvironment{remark}{\begin{rmrk}\rm}{\end{rmrk}}

\ifbook{
    \newcommand{\Chapter}[1]{\section{#1}}
    \newcommand{\Section}[1]{\subsection{#1}}
    \newcommand{\Subsection}[1]{\subsubsection{#1}}
    \def\Chname{Section }

  }
  {
    \newcommand{\Chapter}[1]{\chapter{#1}}
    \newcommand{\Section}[1]{\section{#1}}
    \newcommand{\Subsection}[1]{\subsection{#1}}
    \def\Chname{Chapter}
    
  }

\renewcommand{\(}{$\,}
\renewcommand{\)}{\,$}

\def\nquad{\hspace{-1cm}}
\def\eqdef{\stackrel{\operatorname{def}}{=}}

\DeclareMathAlphabet{\mathbbmsl}{U}{bbm}{bx}{sl}

\DeclareMathSymbol{\Alpha}{\mathalpha}{operators}{"41}
\DeclareMathSymbol{\Beta}{\mathalpha}{operators}{"42}
\DeclareMathSymbol{\Epsilon}{\mathalpha}{operators}{"45}
\DeclareMathSymbol{\Zeta}{\mathalpha}{operators}{"5A}
\DeclareMathSymbol{\Eta}{\mathalpha}{operators}{"48}
\DeclareMathSymbol{\Iota}{\mathalpha}{operators}{"49}
\DeclareMathSymbol{\Kappa}{\mathalpha}{operators}{"4B}
\DeclareMathSymbol{\Mu}{\mathalpha}{operators}{"4D}
\DeclareMathSymbol{\Nu}{\mathalpha}{operators}{"4E}
\DeclareMathSymbol{\Omicron}{\mathalpha}{operators}{"4F}
\DeclareMathSymbol{\Rho}{\mathalpha}{operators}{"50}
\DeclareMathSymbol{\Tau}{\mathalpha}{operators}{"54}
\DeclareMathSymbol{\Chi}{\mathalpha}{operators}{"58}
\DeclareMathSymbol{\omicron}{\mathord}{letters}{"6F}

\newcommand{\cc}[1]{\mathscr{#1}}
\newcommand{\bb}[1]{\boldsymbol{#1}}

\DeclareFontFamily{U}{mathx}{\hyphenchar\font45}
\DeclareFontShape{U}{mathx}{m}{n}{
<5><6><7><8><9><10>
<10.95><12><14.4><17.28><20.74><24.88>
mathx10
}{}
\DeclareSymbolFont{mathx}{U}{mathx}{m}{n}
\DeclareFontSubstitution{U}{mathx}{m}{n}
\DeclareMathAccent{\widebar}{0}{mathx}{"73}

\renewcommand{\bar}[1]{\widebar{#1}}
\renewcommand{\hat}[1]{\widehat{#1}}
\renewcommand{\tilde}[1]{\widetilde{#1}}

\makeatletter
\def\mathcenterto#1#2{\mathclap{\phantom{#1}\mathclap{#2}}\phantom{#1}}
\let\old@widetilde\widetilde
\def\widetildeto#1#2{\mathcenterto{#2}{\old@widetilde{\mathcenterto{#1}{#2\,}}}}
\let\old@widehat\widehat
\def\widehatto#1#2{\mathcenterto{#2}{\old@widehat{\mathcenterto{#1}{#2\,}}}}
\makeatother

\newcommand{\thankstitle}[1]{\ifthenelse{\equal{#1}{}}{}{\thanks{#1}}}
\newcommand{\thanksau}[1]{\ifthenelse{\equal{#1}{}}{}{\thanks{#1}}}

\newcommand{\aua}[6]
{\def\authora{#1}
\def\runauthora{#2}
\def\addressa{#3}
\def\emaila{#4}
\def\affiliationa{#5}
\def\thanksa{#6}}

\def\theauthors{
\ifau{ 
  \author{
    \authora
    \thanksau{\thanksa}
    \\[5.pt]
    \addressa \\
    \texttt{ \emaila}
  }
}
{  
  \author{
    \authora
    \thanksau{\thanksa}
    \\[5.pt]
    \addressa \\
    \texttt{ \emaila}
    \and
    \authorb
    \thanksau{\thanksb}
    \\[5.pt]
    \addressb \\
    \texttt{ \emailb}
  }
}
{   
  \author{
    \authora
    \thanksau{\thanksa}
    \\[5.pt]
    \addressa \\
    \texttt{ \emaila}
    \and
    \authorb
    \thanksau{\thanksb}
    \\[5.pt]
    \addressb \\
    \texttt{ \emailb}
    \and
    \authorc
    \thanksau{\thanksc}
    \\[5.pt]
    \addressc \\
    \texttt{ \emailc}
  }
} {   
  \author{
    \authora
    \thanksau{\thanksa}
    \\[5.pt]
    \addressa \\
    \texttt{ \emaila}
    \and
    \authorb
    \thanksau{\thanksb}
    \\[5.pt]
    \addressb \\
    \texttt{ \emailb}
    \and
    \authorc
    \thanksau{\thanksc}
    \\[5.pt]
    \addressc \\
    \texttt{ \emailc}
    \and
    \authord
    \thanksau{\thanksd}
    \\[5.pt]
    \addressd \\
    \texttt{ \emaild}
  }
}
}

\renewcommand{\Gamma}{\varGamma}
\renewcommand{\Pi}{\varPi}
\renewcommand{\Sigma}{\varSigma}
\renewcommand{\Delta}{\varDelta}
\renewcommand{\Lambda}{\varLambda}
\renewcommand{\Psi}{\varPsi}
\renewcommand{\Phi}{\varPhi}
\renewcommand{\Theta}{\varTheta}
\renewcommand{\Omega}{\varOmega}
\renewcommand{\Xi}{\varXi}
\renewcommand{\Upsilon}{\varUpsilon}

\def\argmax{\operatornamewithlimits{argmax}}
\def\argmin{\operatornamewithlimits{argmin}}

\def\ev{\bb{e}}
\def\fv{\bb{f}}

\def\hv{\bb{h}}

\def\uv{\bb{u}}

\def\wv{\bb{w}}
\def\xv{\bb{x}}

\def\zv{\bb{z}}

\def\Uv{\bb{U}}

\def\Wv{\bb{W}}
\def\Xv{\bb{X}}
\def\Yv{\bb{Y}}
\def\Zv{\bb{Z}}

\def\alphav{\bb{\alpha}}

\def\epsv{\bb{\varepsilon}}
\def\etav{\bb{\eta}}

\def\thetav{\bb{\theta}}

\def\xiv{\bb{\xi}}

\def\Psiv{\bb{\Psi}}

\def\sumi{\sum_{i=1}^{n}}

\usepackage{color}

\definecolor{blue(pigment)}{rgb}{0.2, 0.2, 0.6}
\definecolor{ultramarine}{rgb}{0.07, 0.04, 0.56}
\definecolor{darkspringgreen}{rgb}{0.09, 0.45, 0.27}
\definecolor{hookersgreen}{rgb}{0.0, 0.44, 0.0}
\definecolor{hgreen}{rgb}{0.0, 0.44, 0.0}
\definecolor{plum(traditional)}{rgb}{0.56, 0.27, 0.52}
\definecolor{purple(html/css)}{rgb}{0.5, 0.0, 0.5}
\definecolor{magenta(dye)}{rgb}{0.79, 0.08, 0.48}

\newcommand{\scorem}[1]{\score_{\!#1}}

\def\jJ{J}
\def\mM{M}
\def\gpks{\mathbbmsl{k}}
\def\eEO{\delta_{0}}
\def\nEOa{\mathbbmsl{n}}
\def\CGPz{\CONSTi_{\cgp}}

\def\cgp{w}

\def\DFM{\mathcal{D}}

\def\KSn{\KS^{\circ}}

\def\nEO{N}

\def\F{F}
\def\Fb{\breve{\F}}
\def\Fb{\Phi}
\def\aIF{\alpha}
\def\bIF{\beta}

\def\dmax{\kappa}

\def\bvn{\bb{\mu}}
\def\svn{\bb{\phi}}

\def\dprmtv{\uv} 
\def\dtarpv{\alphav}
\def\dimav{\hv}
\def\dKS{\Delta}

\def\REO{R}

\def\Tens{\mathcal{T}}

\def\TensU{\Tens}

\def\prmt{\ups}

\def\prmtv{\bb{\prmt}}

\def\prmtvs{\prmtv^{*}}

\def\targ{x}
\def\targv{\bb{\targ}}

\def\tarp{\theta}
\def\tarpv{\bb{\tarp}}
\def\tarps{\tarp^{*}}
\def\tarpvs{\tarpv^{*}}

\def\hmax{\mathsf{c}}

\def\hL{h}

\def\elll{\ell}

\def\lgdL{\elll}

\def\dagg{\prime}

\def\proj{P}

\def\Matr{\mathfrak{M}}

\def\Eta{\mathcal{H}}

\def\HL{\mathbb{m}}

\def\dltw{\delta}

\def\dltwb{\omega}

\def\dltwu{\tau}
\def\dltwd{\dltw^{\dagg}}
\def\dltwbd{\dltwb^{\dagg}}

\def\fgblk{J}

\def\DFN{\DVL}

\def\AFblk{\fgblk}

\def\R{\mathbbmsl{R}}
\def\E{\mathbbmsl{E}}

\def\P{\mathbbmsl{P}}

\def\kappa{\varkappa}

\def\blk{\operatorname{block}}

\def\diag{\operatorname{diag}}

\def\Fr{\operatorname{Fr}}

\def\Var{\operatorname{Var}}

\def\T{\top}
\def\tr{\operatorname{tr}}

\def\cond{\, \big| \,}

\def\nsize{{n}}

\def\sumi{\sum_{i=1}^{\nsize}}

\def\ex{\mathrm{e}}

\def\Id{\mathbbmsl{I}}
\def\Ind{\operatorname{1}\hspace{-4.3pt}\operatorname{I}}

\def\alp{\alpha}

\def\bias{\mathsf{b}}

\def\BB{\mathbbmsl{B}}

\def\BBB{\cc{B}}

\def\BBH{B}

\def\CA{\mathcal{A}}

\def\CAs{\CA^{*}}

\def\CGP{w}

\def\CONST{\mathtt{C} \hspace{0.1em}}
\def\CONSTi{\mathtt{C}}

\def\DPN{D}

\def\DVL{\mathbb{D}}

\def\dimA{\mathbb{p}}

\def\dimAfull{\dimA}

\def\dimfull{\bar{\dimp}}
\def\dimAfull{\bar{\dimA}}

\def\dimttl{\bar{\dimp}}

\def\dimp{p}

\def\dimG{\dimA_{\GP}}

\def\dimq{q}

\def\dimd{d}
\def\dimm{m}

\def\dimm{M}

\def\dens{f}

\def\Eta{\cc{H}}

\def\errKS{\nuv}
\def\errY{U}
\def\errYv{\bb{\errY}}

\def\errKS{\Uv}
\def\errZv{\bb{\omega}}
\def\errKSe{\nu}
\def\errZ{\omega}

\def\eps{\epsilon}			
\def\eps{\varepsilon}

\def\epsA{\nu}

\def\fvs{\fv}  

\def\fs{f}

\def\gp{g}

\def\gpks{\kappa}

\def\gpks{k}

\def\GP{G}

\def\GPY{\Gamma}
\def\GPKS{\mathcal{K}}

\def\GPT{\mathcal{G}}

\def\hKS{\hat{\KS}}

\def\IF{\mathbbmsl{F}}

\def\IFN{\IFL}

\def\IFL{{\mathbbmss{F}}}

\def\IFtotal{F}
\def\IFT{\mathscr{\IFtotal}}

\def\IFTb{\Phi}

\def\IFL{\mathbb{F}}

\def\ima{z}
\def\imav{\bb{\ima}}

\def\imavs{\imav^{*}}
\def\Ima{Z}
\def\Imav{\bb{\Ima}}

\def\KS{A}
\def\KSs{\KS^{*}}

\def\KSm{\KS}

\def\Kappa{\cc{K}}

\def\LL{\cc{L}}

\def\LLp{\mathbbmsl{L}}

\def\mm{m}

\def\muA{\mu}

\def\nuov{\bb{a}}

\def\nup{\eta}
\def\nupv{\bb{\nup}}
\def\nupvs{\nupv^{*}}

\def\nui{s}

\def\nuiv{\bb{\nui}}

\def\nuo{\tau}
\def\nuov{\bb{\nuo}}

\def\Phiv{\bb{\Phi}}

\def\Proj{\Pi}

\def\QP{Q}

\def\rhoIF{\rho}

\def\riskt{\cc{R}}

\def\rr{\mathtt{r}}

\def\score{\nabla}

\def\thetav{\bb{\theta}}
\def\thetavs{\thetav^{*}}

\def\Tau{T}

\def\ups{\upsilon}
\def\upsv{\bb{\ups}}

\def\ups{\upsilon}
\def\upsv{\bb{\ups}}

\def\UV{\mathcal{U}}

\def\UVz{\UV}

\def\Ups{\varUpsilon}
\def\Upsd{\Ups^{\circ}}

\def\VP{V}

\def\wv{\bb{w}}

\def\Wv{\bb{W}}

\def\xivr{\breve{\xiv}}

\def\xx{\mathtt{x}}

\def\zq{z}

\def\thetitle{Estimation and inference in error-in-operator model}
\def\theruntitle{Estimation and inference in error-in-operator model}

\def\theabstract{
Many statistical problems can be reduced to a linear inverse problem in which 
only a noisy version of the operator is available. 
Particular examples include random design regression, 
deconvolution problem, 
instrumental variable regression, 
functional data analysis,
error-in-variable regression, drift estimation in stochastic diffusion, and many others.
The pragmatic plug-in approach can be well justified in the classical asymptotic setup 
with a growing sample size.
However, recent developments in high dimensional inference reveal some new features of this problem.
In high dimensional linear regression with a random design, the plug-in approach is questionable
but the use of a simple ridge penalization yields a benign overfitting phenomenon; see \cite{baLoLu2020},
\cite{ChMo2022}, \cite{NoPuSp2024}.
This paper revisits the general Error-in-Operator problem for finite samples and high dimension 
of the source and image spaces. 
A particular focus is on the choice of a proper regularization.
We show that a simple ridge penalty (Tikhonov regularization) works properly in the case
when the operator is more regular than the signal. 
In the opposite case, some model reduction technique like spectral truncation should be applied. 
}

\def\kwdp{62F10,62E17}
\def\kwds{62J12}

\def\thekeywords{semiparametric MLE, risk bounds, ridge penalization, spectral cut-off}

\def\thankstitle{}

\aua
{Vladimir Spokoiny}
{Spokoiny, V. }
{
    Weierstrass Institute and HU Berlin,  
    \\
    Mohrenstr. 39, 10117 Berlin, Germany
    }
{spokoiny@wias-berlin.de}
{Weierstrass-Institute and Humboldt University Berlin}
{
  Financial support by the German Research Foundation (DFG) through the Collaborative Research Center 1294 ``Data assimilation'' is gratefully acknowledged.
}

\begin{document}
\thispagestyle{empty}
{
\title{\thetitle}
\theauthors

\maketitle
\begin{abstract}
{\footnotesize \theabstract}
\end{abstract}

\ifAMS
    {\par\noindent\emph{AMS 2010 Subject Classification:} Primary \kwdp. Secondary \kwds}
    {\par\noindent\emph{JEL codes}: \kwdp}

} 

\tableofcontents


\Chapter{Introduction}
\label{SEiOintr}

%
Let \( \KS \) be a linear mapping (operator) of the source signal \( \tarpv \in \R^{\dimp} \) 
to the image space \( \R^{\dimq} \).
Consider the problem of recovering \( \tarpv \in \R^{\dimp} \) from the noisy image 
\( \Zv \) satisfying \( \E \Zv = \KS \tarpv \).
The standard least square approach suggests to recover \( \tarpv \) by maximizing 
the negative fidelity function (log-likelihood) \( \lgdL(\cdot) \) of the form 
\begin{EQA}
	\lgdL(\tarpv) 
	&=& 
	- \frac{1}{2 \sigma^{2}} \| \Zv - \KS \tarpv \|^{2} ,
\label{o3kv6ejfgiuhbdfjkse}
\end{EQA}
where \( \sigma^{2} \) describes the image noise.
However, this approach assumes the operator \( \KS \) to be precisely known.
In many applications, this assumption is not fulfilled and the operator \( \KS \) is known up to some error, or, equivalently,
only an estimate \( \hKS \) of \( \KS \) is available.
Particular examples include regression with random design and
instrumental variable regression among  many others.
A pragmatic plug-in approach just replaces \( \KS \) in \eqref{o3kv6ejfgiuhbdfjkse} by its estimate \( \hKS \) leading 
to the solution
\begin{EQA}
	\hat{\tarpv}
	&=&
	\argmin_{\tarpv} \| \Zv - \hKS \tarpv \|^{2}
	=
	\bigl( \hKS^{\T} \hKS \bigr)^{-1} \hKS^{\T} \zv .
\label{wkvoigb7y634hjgtkuiif}
\end{EQA}
Unfortunately, inverting a large matrix \( \hKS^{\T} \hKS \) can be tricky, 
especially if the operator \( \KS \) is smooth and the dimensions \( \dimp, \dimq \)
are large or infinite.
\cite{HoRe2008} considered simultaneous wavelet estimation of the signal \( \tarpv \) and the operator \( \KS \).
\cite{Trabs_2018} extended these results to Bayesian inference using a parametric assumption
about the unknown operator \( \KS = \KS_{\thetav} \) and provided some examples
from imaging and linear PDE/SDEs. 
In those papers, one can also find a literature overview and more references on this topic.

In this note, we make a further step by including the operator and its image \( \imav = \KS \thetav \) 
in the set of parameters and by considering the original problem in a semiparametric setup. 
This leads to the new fidelity function
\begin{EQA}
	\LL(\tarpv,\imav,\KS)
	&=&
	- \frac{1}{2 \sigma^{2}} \| \Zv - \imav \|^{2} 
	- \frac{\muA^{2}}{2 \sigma^{2}} \| \hKS - \KS \|_{\Fr}^{2} 
	- \frac{1}{2} \| \imav - \KS \tarpv \|^{2} \, , 
\label{lgbhu3hdf6vjghhksdemde}
\end{EQA}
where \( \muA \) describes the level of the operator noise and \( \| \KS \|_{\Fr}^{2} = \tr(\KS^{\T} \KS) \) stands
for the squared Frobenius norm.
An advantage of such an extension is that the random data \( \Zv \) and \( \hKS \) only enter
in the nice quadratic terms \( \| \Zv - \imav \|^{2} \) and \( \| \hKS - \KS \|_{\Fr}^{2} \).
This gradually simplifies the statistical analysis and allows to avoid the heavy burden
of empirical process theory by reducing the problem to the SLS setup from \cite{Sp2024}.
Another benefit of this extension is that the auxiliary parameters \( \imav \) and \( \KS \)
can be of independent interest and one can use some regularization techniques to control the complexity
of the related variables. 
This especially concerns the operator \( \KS \) which is typically smooth and its complexity
can be reduced by proper penalization.
Later we consider the penalized least squares method based on maximization of the penalized 
objective function \( \LL_{\GPT}(\tarpv,\imav,\KS) \) defined as
\begin{EQA}[c]
	\LL_{\GPT}(\tarpv,\imav,\KS)
	\eqdef
	\LL(\tarpv,\imav,\KS) 
	- \frac{1}{2} \| \GP \tarpv \|^{2} - \frac{1}{2} \| \KS \|_{\GPKS}^{2}
	\, ,
\end{EQA}
where \( \| \GP \tarpv \|^{2} \) mimics complexity of the target parameter \( \tarpv \)
while \( \| \KS \|_{\GPKS}^{2} \) is the operator penalty.
The choice of the penalization operator \( \GP^{2} \) is of special 
importance and will be discussed later in detail for two extreme cases:
a ridge penalty \( \GP^{2} = \gp^{2} \Id_{\dimp} \) and a truncation penalty 
\( \GP^{2} = \diag\{ \gp_{1}^{2}, \ldots, \gp_{\dimp}^{2} \} \)
with \( \gp_{j}^{2} = 0 \) for \( j \leq \jJ \) and \( \gp_{j}^{2} = \infty \) 
for \( j > \jJ \).
Ridge or Tikhonov regularization is widely used in practical applications, especially,
for ill-posed inverse problems. 
In the recent time, it gained new impulses in connection to 
``benign overfitting'' effect in high dimensional regression; \cite{baLoLu2020},
\cite{ChMo2022}.
An important benefit of ridge penalty is that it is coordinate free and does not rely 
on the eigenvalue decomposition of the operator \( \KSs \).
This is particularly crucial for the error-in-operator problem when the operator 
\( \KSs \) itself is not known. 
From the other side, it is well known that ridge regression is inefficient 
in the classical ``direct problem'' setup when the conditional number 
of the matrix \( {\KSs}^{\T} \KSs \) is bounded by a fixed constant.
In such cases, model reduction is usually achieved using roughness penalty 
or spectral cut-off approach.
One of the message of this paper is that the ridge regression procedure with 
a carefully selected penalty constant \( \gp^{2} \)
provides nearly optimal accuracy of estimation in situations 
when the operator \( \KSs \) is ``smoother'' as the signal \( \tarpvs \).
In other cases, one can use a spectral cut-off method which relies either 
on the singular value decomposition of the unknown operator \( \KSs \)
or on the approximation space representation. 
The following two examples illustrate our setup.

\paragraph{Example: Random design regression}
\label{Srandomdesign}
In many applications, the design \( X_{1},\ldots,X_{n} \) in regression models 
\( Y_{i} = \fs(X_{i}) + \eps_{i} \) may be naturally assumed to be random, often i.i.d.
This particularly concerns econometric, biological, chemical, sociological studies etc. 
Under linear parametric assumption \( \fs(\xv) = \sum_{j} \theta_{j} \psi_{j}(\xv) = \Psi(\xv)^{\T} \thetav \) and noise homogeneity, 
the log-likelihood reads exactly as in the case of a deterministic design:
\begin{EQA}
	L(\thetav)
	&=&
	- \frac{1}{2 \sigma^{2}} \| \Yv - \Psiv(\Xv)^{\T} \thetav \|^{2} + \CONST ,
\label{4556799kbvdfee4566y} 
\end{EQA}
where \( \Psiv(\Xv) \) is the \( \dimp \times n \) matrix with entries \( \psi_{j}(X_{i}) \).
The corresponding MLE \( \tilde{\thetav} \) is 
\begin{EQA}
	\tilde{\thetav}
	&=&
	\bigl\{ \Psiv(\Xv) \, \Psiv(\Xv)^{\T} \bigr\}^{-1} \Psiv(\Xv) \Yv .
\label{kvfduyerftywefope0rg}
\end{EQA}
However, a study of this estimator is quite involved due to the random nature of the matrix 
\( \Psiv(\Xv) \), especially in the case of a high design dimension \( \dimd \);
see recent papers 
\cite{baLoLu2020},
\cite{ChMo2022},
and references therein 
for the case of linear regression with \( \dimd = n \).
%
%
Suppose to be given another collection of features \( (\phi_{m}) \) for \( m=1,\ldots,\dimq \).
The use of \( (\phi_{m}) = (\psi_{j}) \) is one of possible options.
However, one can use many others including biorthogonal bases, wavelets, etc.
Define \( \Phiv(X_{i}) = (\phi_{m}(X_{i})) \in \R^{\dimq} \) and \( \Phiv(\Xv) = (\Phiv(X_{1}), \ldots,\Phiv(X_{n})) \).
Introduce a linear operators \( \hKS \): \( \R^{\dimp} \to \R^{\dimq} \) by 
\begin{EQA}
\label{iiu87654e3wsdrtyuiop}
	\hKS
	& \eqdef &
	\Phiv(\Xv) \, \Psiv(\Xv)^{\T}
	=
	\sumi \Phiv(X_{i}) \, \Psiv(X_{i})^{\T} ,
\end{EQA}
and define \( \E \hKS = \KSs \), so \( \hKS \) might be viewed as an empirical version of \( \KSs \).
If the \( X_{i} \)'s are i.i.d. and \( \dens_{X} \) is the design density then
\begin{EQA}
	\KSs
	&=&
	n \E \, \Phiv(X_{1}) \, \Psiv(X_{1})^{\T}
	=
	n \int \Phiv(\xv) \, \Psiv(\xv)^{\T} \, \dens_{X}(\xv) \, d\xv \, ,
\label{plkmnbvcxcvbnmnbvcvbnm}
\end{EQA}  
In general, with \( \hKS_{i} = n \Phiv(X_{i}) \, \Psiv(X_{i})^{\T} \), it holds 
\( \hKS = n^{-1} \sumi \hKS_{i} \) and also
\begin{EQA}
	\sumi \| \hKS_{i} - \KS \|_{\Fr}^{2}
	&=&
	\sumi \| \hKS_{i} - \hKS \|_{\Fr}^{2} + n \| \hKS - \KS \|_{\Fr}^{2} \, .
\label{iujwefy7wefyikjt5r433e}
\end{EQA}
As \( \sumi \| \hKS_{i} - \hKS \|_{\Fr}^{2} \) only depends on the design \( \Xv \),
one may equally use \( \sumi \| \hKS_{i} - \KS \|_{\Fr}^{2} \) and 
\( n \| \hKS - \KS \|_{\Fr}^{2} \)
to measure the data misfit of a possible guess \( \KS \).
Further, denote
\begin{EQA}
	Z_{m}
	& \eqdef &
	\sumi Y_{i} \, \phi_{m}(X_{i}) 
\label{dcirf7rft64rt6rfgh}
\end{EQA}
and let \( \Zv \) be the vector in \( \R^{\dimp} \) with entries \( Z_{m} \).
The LPA \( \fvs(\xv) = \Psiv(\xv)^{\T} \, \thetav \) together with \( \Yv = \fvs + \epsv \) and  \( \Uv = \Phiv(\Xv) \epsv \) implies
\begin{EQA}
	\Zv
	&=&
	\Phiv(\Xv) \fvs + 
	\Phiv(\Xv) \epsv
	=
	\Phiv(\Xv) \, \Psiv(\Xv)^{\T} \, \thetav + 
	\Phiv(\Xv) \epsv
	=
	\hKS \, \thetav + \Uv .
\label{poiuytrr43ewsdere4dr}
\end{EQA}
The vector \( \Zv \) might be viewed as new ``observations'' of the ``response'' \( \imav = \KS \thetav \)
corrupted by two sources of errors: the additive error \( \Uv \) and the error in operator \( \hKS - \KS \).
These errors are not independent, however, \( \E (\Uv \cond \Xv) = 0 \).
If the errors \( \eps_{i} = Y_{i} - \fs(X_{i}) \) are Gaussian independent of the design variables \( X_{i} \), then
the new ``noise'' \( \Uv \) is Gaussian conditionally on \( \Xv \).
This enables us to decouple these sources in the likelihood by considering in place of \( \| \Zv - \hKS \thetav \|^{2} \)
two fidelity terms \( \| \Zv - \imav \|^{2} \) and 
\( n \| \hKS - \KS \|_{\Fr}^{2} \) subject to the structural 
equation \( \imav = \KS \thetav \).
Using the Lagrange multiplier idea, we put all together in one expression
\begin{EQA}
	\LL(\thetav,\imav,\KS)
	& \eqdef &
	- \frac{1}{2 \sigma^{2}} \| \Zv - \imav \|^{2}
	- \frac{n}{2 \sigma_{X}^{2}} \| \hKS - \KS \|_{\Fr}^{2} 
	- \frac{\lambda}{2 \sigma^{2}} \| \imav - \KS \thetav \|^{2} \, .
\label{3ewsazsplokfr5}
\end{EQA}
The first term specifies the quality of fitting the ``data'' \( \Zv \) by the ``response'' \( \imav \).
The second term is a similar fidelity term for the operator \( \KS \) observed with the operator noise of variance 
\( \sigma_{\KS}^{2} = n^{-1} \sigma_{X}^{2} \).
The last term transfers the structural equation \( \imav = \KS \thetav \) into a penalty term.
Note that this equation is a special case of \eqref{lgbhu3hdf6vjghhksdemde}.
We refer to the forthcoming paper \cite{NoPuSp2024} for a rigorous study of a high-dimensional
linear regression problem with random design.


\paragraph{Example: Instrumental variable regression}
Nonparametric IV model assumes that we are interested in the functional dependence \( Y \) on \( X \) 
in the form \( Y = \fs(X) \), however, this model is probably incomplete and hence, 
\( \E \bigl\{ Y - \fs(X) \cond X \bigr\} \neq 0 \).
Instead, one assumes that \( W \) is a proper instrument and \( \E\bigl\{ Y - \fs(X) \cond W \bigr\} = 0 \).
The goal is to recover \( \fs \).
Let triples \( (Y_{i},X_{i},W_{i}) \) be observed, \( i=1,\ldots,n \). 
We apply a linear expansion of the target \( \fs \) using a dictionary \( (\psi_{j}) \):
\begin{EQA}
	\fs(\xv)
	&=&
	\sum_{j=1}^{\dimp} \tarp_{j} \, \psi_{j}(\xv)
	=
	\Psiv(\xv)^{\T} \tarpv .
\label{jdtrderewertytrerertgh}
\end{EQA}
Let \( \bigl\{ \phi_{m}(w), m=1,\ldots,\dimq \bigr\} \) be a function basis for the instrument. 
Define 
\begin{EQA}
	\Ima_{m}
	& \eqdef &
	\sumi Y_{i} \, \phi_{m}(W_{i}) ,
	\qquad
	\errY_{m}
	\eqdef 
	\sumi \eps_{i} \, \phi_{m}(W_{i}) ,
	\qquad
	m=1,\ldots,\dimq.
\label{KjmEspjXimpW}
\end{EQA}
The model assumption yields \( \E \errY_{m} = 0 \) for each \( m \).	
Introduce the linear operator \( \hKS \colon \R^{\dimp} \to \R^{\dimq} \) with
\begin{EQA}
	\hKS 
	& \eqdef &
	\biggl( \sumi \psi_{j}(X_{i}) \phi_{m}(W_{i}) \biggr)_{j \leq \dimp, \, m \leq \dimq} \, 
\label{EXiWim1qfX}
\end{EQA}
and define \( \KSs = \E \hKS \).
The operator \( \KSs \) expresses the joint distribution of the regressor \( X \) and the instrument \( W \).
Further we consider the vector \( \Imav \) and the matrix \( \hKS \) as our new observations.
Then the IV model can be rewritten as
\begin{EQA}
	\Imav
	&=&
	\hKS \tarpv + \errYv ,
	\qquad
	\hKS
	=
	\KSs + \epsA .
\label{KdKtUZtKkK}
\end{EQA}
%
This repeats the random design case \eqref{poiuytrr43ewsdere4dr}.
The calming device with the auxiliary ``image'' vector \( \imav \) 
linking the observation \( \Imav \) and the product \( \KS \tarpv \) by the structural relation 
\( \imav = \KS \tarpv \) leads to the following expression
\begin{EQA}
	\LL(\tarpv,\imav,\KS)
	&=&
	- \frac{1}{2} \| \Imav - \imav \|^{2}
	- \frac{\muA^{2}}{2} \| \hKS - \KS \|_{\Fr}^{2} 
	- \frac{1}{2} \| \imav - \KS \tarpv \|^{2}
	\, ,
\label{LfA1212sZttA}
\end{EQA}
with \( \muA^{2} \) reflecting the error of operator estimation \( \hat{\KS} - \KSs \).
Due to \eqref{EXiWim1qfX}, it is natural to apply \( \muA^{2} \asymp n \).
It is important to keep in mind that \( \Imav, \hKS \) are not our original data, 
they are computed by \eqref{KjmEspjXimpW} and \eqref{EXiWim1qfX} out of the triple \( (\Yv,\Xv,\Wv) \).
We refer to the recent paper \cite{XiWa2024} for an introduction and a detailed literature overview 
for the IV problem. 
The standard plug-in approach applied in that paper allows to obtain rather sharp finite sample accuracy results
about accuracy of estimating the parameter \( \tarpv \) in \eqref{jdtrderewertytrerertgh}
in a linear parametric setup.
However, it requires some conditions limiting the scope of applicability of this approach.
The most critical condition is that the dimensions \( \dimp \) and \( \dimq \) should coincide.
Smoothness of the operator is not addressed as well as the estimation bias due to linear parametric expansion
of the regression function. 
Also, the results do not provide any dependence of the estimation error on \( \dimp \),
and the case of a high dimensional parameter \( \tarpv \) is not included.
The results obtained later in this paper refine and improve 
the accuracy guarantees from \cite{XiWa2024} and apply to the nonparametric IV problem.

\medskip
Many other models can be reduced to the EiO-setup. 
Particular examples include deconvolution problem, 
functional data analysis,
error-in-variable regression, drift estimation in stochastic diffusion.

\paragraph{This paper contribution}
This paper offers a new approach to inference in the Error-in-Operator model
based on extending the parameter space by including the operator \( \KS \)
and the image \( \imav \) in the list of parameters.
The extended log-likelihood is stochastically linear and a polynomial of order four 
w.r.t. the set of parameters and hence very smooth. 
Moreover, after a minor regularization, it appears to be strongly concave in 
a rather large vicinity of the zero point in the full-dimensional parameter space.
Therefore, it can be studied by general methods developed in \cite{Sp2024} 
for stochastically linear smooth (SLS) models. 
The main results provide finite sample expansions for the profiled MLE \( \tilde{\tarpv} \)
as well as sharp risk bounds with the explicit leading term. 
In particular, we present sufficient conditions ensuring the oracle type properties 
of this estimator: it behaves essentially as well as the usual LSE with a precisely known operator. 
The theoretical guarantees are specified for two important classes of penalties:
a ridge (Tikhonov) regularization and a spectral cut-off method. 
We explain a phase transition effect: 
a simple ridge penalty with a proper choice of the ridge parameter ensures
accurate and rate optimal estimation in the situation when the operator is 
more regular than the source signal.
In the opposite case, one should use a model reduction technique like spectral cut-off
to achieve the optimal rate of estimation. 
A data-driven choice of the penalty parameters is not discussed in this paper.
Note however, that a standard cross-validation approach is well applicable in this setup, and 
the established Wilks expansion can be very useful.

\paragraph{Outlook}
Some questions important for practical applications, are not addressed by this study.
A data-driven choice of a ridge parameter \( \gp^{2} \) or a truncation level \( \jJ \) for 
a spectral cut-off methods can be done by usual cross-validation technique, and
the obtained Wilks expansion are very useful for justifying their validity.
However, we refer to a forthcoming paper for treating this problem 
by an extension of our semiparametric framework when the tuning parameter is considered
as a part of the full dimensional parameter.

Another issue is a check of the condition \nameref{EU2ref}, especially for the operator noise.
This is application dependent and
a careful treatment for each particular case requires a substantial amount of work.
We refer to \cite{NoPuSp2024} and \cite{PuNoSp2024} for the case of random design,
and to \cite{XiWa2024} for the IV model.
However, as soon as \nameref{EU2ref} has been checked, 
the general results of Theorem~\ref{TsemieffEO} can be applied. 

\paragraph{Organization of the paper}
Section~\ref{SEiOMLE} considers a general 
Error-in-Operator framework and presents some general 
results on accuracy of estimation in terms of \emph{efficient dimension}; see Theorem~\ref{TsemieffEO}.
We also introduce the \emph{critical dimension} condition
meaning the relation between the effective sample size 
and efficient dimension which enables consistent estimation; see Section~\ref{ScritdimEO}. 
Section~\ref{SpenaltyEO} addresses a proper choice of penalization for direct and mildly inverse problems.
%
The proofs of the results on estimation for Error-in-Operator models from Section~\ref{SEiOMLE}
are presented in Appendix~\ref{SEiOtools}.


\Chapter{Theoretical study for the error-in-operator model}
\label{SEiOMLE}

This section establishes some theoretical guarantees for the proposed semiparametric approach 
in Error-in-Operator models.
We adopt the setup of \ifshort{Section~\ref{SEiOintr}}{\cite{Sp2024}}.
For a noisy image \( \Imav \) and a pilot estimate \( \hKS \) of the operator \( \KS \),
the task is to recover the source signal \( \tarpv \) following the equation
\( \KS \tarpv = \E \Imav \).
We treat the image \( \imav \in \R^{\dimq} \) and the operator \( \KS \in \R^{\dimq\times\dimp} \) as nuisance parameters.
This leads to the extended model 
\begin{EQA}
	\LL(\prmtv)
	=
	\LL(\tarpv,\imav,\KS)
	&=&
	- \frac{1}{2} \| \Imav - \imav \|^{2}
	- \frac{\muA^{2}}{2} \| \KS - \hKS \|_{\Fr}^{2}
	- \frac{1}{2} \| \imav - \KS \tarpv \|^{2} \, .
\label{dhfdfededrsdsrdsresswddd}
\end{EQA}
The full-dimensional vector \( \prmtv \) is a point in \( \R^{\dimttl} \), where \( \dimttl = \dimp + \dimq + \dimp \dimq \).
The factor \( \muA^{2} \) describes the ``operator noise''.
Later we assume the error in operator to be smaller than the observation noise.
This leads to a large value \( \muA^{2} \), which is essential for parameter identifiability.
Another important issue is smoothness of the source signal \( \tarpv \) and the operator \( \KS \).
A standard way to address this point is given by regularization in terms of two quadratic penalties
\( \| \GP \tarpv \|^{2}/2 \) and \( \| \KS \|_{\GPKS}^{2}/2 \).
We denote 
\begin{EQA}
	\LL_{\GPT}(\prmtv)
	&=&
	\LL_{\GPT}(\tarpv,\imav,\KS)
	=
	\LL(\tarpv,\imav,\KS) - \frac{1}{2} \| \GP \tarpv \|^{2} - \frac{1}{2} \| \KS \|_{\GPKS}^{2} \, . 
\label{iduejef634yhdjdekee}
\end{EQA}
Typical examples of choosing \( \GP^{2} \) correspond to 
the ridge regression case with \( \GP^{2} = \gp^{2} \Id_{\dimp} \) and
the roughness penalty \( \GP^{2} = \CGP^{2} \diag(j^{2\beta}) \), while a truncation penalty assumes 
\( \GP_{j}^{2} = 0 \) for \( j \leq \jJ \) and \( \GP_{j}^{2} = \infty \) for \( j > \jJ \)
with a fixed spectral cut-off parameter \( \jJ \).
The elementwise operator penalty \( \GPKS = \diag\{ \gpks_{mj} \} \) leads to
\begin{EQA}[c]
	\| \KS \|_{\GPKS}^{2}
	=
	\sum_{\mm=1}^{\dimq} \sum_{j=1}^{\dimp} \gpks_{mj}^{} \KS_{mj}^{2}
	\, .
\label{4f8g8trthwgitiehrEO}
\end{EQA}
The approximation space approach as in \cite{HoRe2008,Trabs_2018} with 
\( \gpks_{mj} = \gpks_{\mm} \) for \( \gpks_{1}^{2} \leq \ldots \leq \gpks_{\dimq}^{2} \)
results in
\begin{EQA}
	\| \KS \|_{\GPKS}^{2}
	&=&
	\sum_{\mm=1}^{\dimq} \gpks_{\mm}^{2} \| \KSm_{\mm} \|^{2} 
	\, ,
\label{4f8g8trthtiehrEOr}
\end{EQA}
where \( \KS_{\mm}^{\T} \) is the \( \mm \)th row of \( \KS \), \( \KS_{\mm} \in \R^{\dimp} \).
%

Define the penalized MLE \( \tilde{\prmtv}_{\GPT} \), its population counterpart \( \prmtvs_{\GPT} \), 
and the truth \( \prmtvs \) as
\begin{EQA}
	\tilde{\prmtv}_{\GPT}
	&=&
	(\tilde{\tarpv}_{\GPT},\tilde{\imav}_{\GPT},\tilde{\KS}_{\GPT})
	=
	\argmax_{\prmtv} \LL_{\GPT}(\prmtv) ,
\label{7cxjue8e3hfft3y3hwfire}
	\\
	\prmtvs_{\GPT}
	&=&
	(\tarpvs_{\GPT},\imavs_{\GPT},\KSs_{\GPT})
	=
	\argmax_{\prmtv} \E \LL_{\GPT}(\prmtv) .
	\\
	\prmtvs
	&=&
	(\tarpvs,\imavs,\KSs)
	=
	\argmax_{\prmtv} \E \LL(\prmtv) .
\label{7cxjue8e3hfft3y3hwfireE}
\end{EQA}
The main focus of the study is the \( \tarpv \)-component \( \tilde{\tarpv}_{\GPT} \) of
\( \tilde{\prmtv}_{\GPT} \) and its loss \( \tilde{\tarpv}_{\GPT} - \tarpvs \).
%
As a benchmark, consider the simple linear problem \( \Imav = \KSs \tarpv + \errZv \) with the known operator \( \KSs \).
The use of penalization \( \| \GP \tarpv \|^{2}/2 \) leads to the penalized MLE 
\begin{EQA}
	\hat{\tarpv}_{\GP} 
	&=&
	\argmax_{\tarpv} \bigl\{ - \frac{1}{2} \| \Imav - \KSs \tarpv \|^{2} - \frac{1}{2} \| \GP \tarpv \|^{2} \bigr\}
	= 
	({\KSs}^{\T} \KSs + \GP^{2})^{-1} {\KSs}^{\T} \Imav = \IF_{\GP}^{-1} {\KSs}^{\T} \Imav 
\label{oc8e7r46gyhr4ufwwwfc}
\end{EQA}
with \( \IF_{\GP}^{2} = {\KSs}^{\T} \KSs + \GP^{2} \).
Then \( \hat{\tarpv}_{\GP} \) concentrates on the set
\begin{EQA}[rcl]
	\CAs_{\GP}
	& \eqdef &
	\bigl\{ \tarpv \colon \| \IF_{\GP}^{1/2} (\tarpv - \tarpvs_{\GP}) \| \leq 3 \rr_{\GP}/2 \bigr\} ,
	\qquad
	\rr_{\GP} = \zq(\BB_{\GP},\xx) ,
\label{EGeHGeeG3rG2EOk}
\end{EQA}
where \( \BB_{\GP} = \IF_{\GP}^{-1/2} {\KSs}^{\T} \Var(\Imav) \KSs \, \IF_{\GP}^{-1/2} \).
Our main results state that the profile MLE \( \tilde{\tarpv}_{\GPT} \) performs essentially as the ideal
estimator \( \hat{\tarpv}_{\GP} \).

\Section{Main definitions and model conditions}
This section defines our main objects, including the Fisher information matrix, and presents the conditions on the model.

\Subsection{Full dimensional information matrix}

For the extended model \eqref{dhfdfededrsdsrdsresswddd}, 
consider the function \( \fs(\prmtv) = \E \LL(\prmtv) \). 
Its most important characteristic is the corresponding information matrix 
\( \IFT(\prmtv) = - \nabla^{2} \E \LL(\prmtv) \).
We write \( \IFT = \IFT(\prmtvs_{\GPT}) \).
Alternatively, one can use \( \IFT = \IFT(\prmtvs) \).
The diagonal blocks of \( \IFT \) will be denoted as 
\( \IFT_{\tarpv\tarpv} \), \( \IFT_{\imav\imav} \), and \( \IFT_{\KS\KS} \).
Similarly, the off-diagonal blocks are \( \IFT_{\tarpv\imav} \), \( \IFT_{\tarpv\KS} \), 
and \( \IFT_{\imav\KS} \) leading to the block structure
\begin{EQA}
	\IFT
	&=&
	- \nabla^{2} \E \LL(\prmtvs_{\GPT})
	=
	\begin{pmatrix}
		\IFT_{\tarpv\tarpv} & 	\IFT_{\tarpv\imav} 	& \IFT_{\tarpv\KS} \\ 
		\IFT_{\imav\tarpv} 	&	\IFT_{\imav\imav} 	& \IFT_{\imav\KS}  \\
		\IFT_{\KS\tarpv} 	& 	\IFT_{\KS\imav} 	& \IFT_{\KS\KS}
	\end{pmatrix}
	\, .
\label{uy8csadjdscfa8jkmwmvus}
\end{EQA}
The roughness penalty on the target parameter \( \tarpv \) and the operator \( \KS \) leads to the penalized information matrix
\begin{EQA}
	\IFT_{\GPT}
	& \eqdef &
	\IFT + \GPT^{2}
	=
	\IFT + \blk\{ \GP^{2}, 0, \GPKS^{2} \} \, .
\label{t5dhxcji9e6fchfcdheyt3}
\end{EQA}
The main results involve the inverse of this matrix and, in particular, its \( \tarpv\tarpv \)-block.
By Schur complement; see Lemma~\ref{LSchur}; it can be represented as \( \IFTb_{\GPT,\tarpv\tarpv}^{-1} \),
where 
\begin{EQA}
	\IFTb_{\GPT,\tarpv\tarpv} 
	&=& 
	\IFT_{\tarpv\tarpv} + \GP^{2} - \IFT_{\tarpv\nupv} \, \IFT_{\GPT,\nupv\nupv}^{-1} \IFT_{\nupv\tarpv} \, ,
\label{kf8r484r5iu8t7t7eurggEO}
\end{EQA}
\( \nupv \) codes the bunch of nuisance parameters \( \nupv = (\imav,\KS) \),
and \( \IFT_{\GPT,\nupv\nupv} \) is the corresponding sub-matrix of \( \IFT_{\GPT} \) composed of the blocks 
\( \IFT_{\imav\imav} \), \( \IFT_{\imav\KS} \), \( \IFT_{\KS\imav} \), and \( \IFT_{\GPT,\KS\KS} = \IFT_{\KS\KS} + \GPKS^{2} \).
This matrix \( \IFTb_{\GPT,\tarpv\tarpv} \) is referred to as \emph{semiparametric information matrix}.
The basic fact about the extended model \eqref{dhfdfededrsdsrdsresswddd} is that 
the objective function \( \LL(\prmtv) \) is strongly concave after restricting to a rather large convex region of the parameter space;
see \eqref{rughur47478ytfg84g6gh} later.
This key funding will be systematically used for all the results.
A sufficient condition for concavity is that the full dimensional matrix \( \IFT_{\GPT}(\prmtv) \)
is well posed and bounded from below by a multiple of the identity matrix; see
Lemma~\ref{LIFTDFMEO}. 
In particular, \( \IFTb_{\GPT,\tarpv\tarpv} \geq \dmax^{-2} {\KSs}^{\T} \KSs + \GP^{2} \) 
for \( \dmax = 2 \).

\Subsection{Local concavity and local smoothness}
The generic result of \ifshort{Proposition~\ref{PconcMLEgenc}}{Theorem~2.7 of \cite{Sp2024}}
requires one more important condition \nameref{LLref} about concavity of \( \fs(\prmtv) = \E \LL(\prmtv) \).
This condition is not fulfilled on the whole parameter space \( \Ups \) for the full dimensional parameter 
\( \prmtv = (\tarpv,\imav,\KS) \) due to the nonlinear product term
\( \KS \tarpv \) in the structural penalty \( \| \imav - \KS \tarpv \|^{2} \).
However, restricting the parameter space to a local subset \( \Upsd \) enables us to state local concavity of 
\( \LL_{\GPT}(\prmtv) \) and parameter identifiability.

Let a matrix \( \DPN^{2} \in \Matr_{\dimp} \) satisfy 
\begin{EQA}[c]
	\DPN^{2} 
	\geq 
	{\KSs}^{\T} \KSs 
	\, .
\label{tfeuhfuwr365vdyge34def}
\end{EQA}
A natural candidate is \( \DPN^{2} = {\KSs}^{\T} \KSs \) provided that the latter product
is well posed.
Another proposal could be 
\( \DPN^{2} = {\KSs}^{\T} \KSs + \gp_{0}^{2} \Id_{\dimp} \) for some small \( \gp_{0} \).
Also, define
\begin{EQA}
	\nEO
	& \eqdef &
	\lambda_{\min}(\DPN^{2}) 
	\, .
\label{tt7t7yrs3esfhiyfsftgh}
\end{EQA}
The quantity \( \nEO \) is important and will be referred to as \emph{effective sample size}.
Our results state root-\( \nEO \) the accuracy of estimation.
With some \( \eEO < 1 \), 
define the radius \( \REO = \eEO \muA \sqrt{\nEO} \) and the local set
\begin{EQA}
	\Upsd
	& \eqdef &
	\bigl\{ (\tarpv,\imav,\KS) \colon 
		\| \DPN \tarpv \| \leq \REO, \,
		\| \imav \| \leq \REO, \,
		\| (\KS - \KSs) \DPN^{-1} \|_{\Fr} \leq \eEO
	\bigr\} \, .
	\qquad
\label{rughur47478ytfg84g6gh}
\end{EQA}
Later we use \( \eEO \leq 1/10 \) and assume that the full dimensional estimator \( \tilde{\prmtv}_{\GPT} \) is limited to \( \Upsd \) 
and this set contains \( \prmtvs \) and hence \( \prmtvs_{\GPT} \).
%
The function \( \fs(\prmtv) = \LL(\tarpv,\imav,\KS) \) is very smooth in the scope of parameters, in fact,
it is a polynomial of degree 4. 
In Section~\ref{SsmoothEO}, we show that 
the full dimensional smoothness conditions \nameref{LLsT3ref} and \nameref{LLsT4ref} are fulfilled for the metric tensor
\begin{EQA}
	\DFM^{2}
	&=&
	\blk\{ \DFM_{\tarpv\tarpv}^{2}, \DFM_{\imav\imav}^{2},\DFM_{\KS\KS}^{2} \} 
	=
	\blk\{ \DPN^{2}, \Id_{\imav}, \muA^{2} \Id_{\KS} \} \, 
\label{iwdfhw787e3urne6fyh3o}
\end{EQA}
and constants
\begin{EQA}
	\dltwu_{3}
	&=&
	4.5 \, \nEO^{-1/2} \muA^{-1} ,
	\qquad
	\dltwu_{4}
	=
	3 \nEO^{-1} \muA^{-2} \, .
\label{gdjdgredtgebegxgEO}
\end{EQA}
Lemma~\ref{LIFTDFMEO} later states strong concavity of \( \fs(\prmtv) = \E \LL(\prmtv) \) on \( \Upsd \) in the form
\( \IFT(\prmtv) = - \nabla^{2} \fs(\prmtv) \geq \dmax^{-2} \DFM^{2} \) 
for \( \dmax = 2 \) and all \( \prmtv \in \Upsd \).

\Subsection{Conditions on the errors}
\label{SEOOnoise}
The stochastic term \( \zeta(\prmtv) = \LL(\prmtv) - \E \LL(\prmtv) \)
is linear and its gradient \( \nabla \zeta(\prmtv) \) is constant in \( \prmtv \) with 
\begin{EQA}
	\scorem{\tarpv} \zeta(\prmtv)
	& \equiv &
	0,
	\qquad
	\scorem{\imav} \zeta(\prmtv)
	\equiv
	\Imav - \E \Imav
	\eqdef
	\errZv,
	\qquad
	\scorem{\KS} \zeta(\prmtv)
	\equiv
	\muA^{2} (\hKS - \E \hKS) 
	=
	\muA \, \errKS 
	\, ,
	\qquad
\label{yte3gf6rf64y4577bu7gtuir}
\end{EQA}
for \( \errKS = \muA (\hKS - \E \hKS) \).
This yields \nameref{Eref}.
We do not further specify the structure of the errors \( \errZv \) and \( \errKS \).
Instead, we only assume that 
condition \nameref{EU2ref} holds with \( \VP^{2} = \Var(\nabla \zeta) \).
General results from \ifsqnorm{Section~\ref{Sdevboundgen}}{\cite{Sp2023c,Sp2023d}} reduce \nameref{EU2ref} 
to locally uniform bounds on the cumulant generating function of \( \nabla \zeta \)
which can be studied independently for \( \errZv \) and \( \errKS \).
The case of Gaussian and sub-Gaussian errors is studied in \ifsqnorm{Section~\ref{SdevboundGauss} and Section~\ref{SdevboundnonGauss}}{\cite{Sp2023c}}
while the sub-exponential case is considered in \ifsqnorm{Section~\ref{Sdevboundexp}}{\cite{Sp2023d}}.
A careful check of such bounds is application-dependent, and later we just assume that \nameref{EU2ref} is fulfilled. 
Examples of checking this condition can be found in \cite{NoPuSp2024} and it is based 
on recent results on concentration of quadratic forms in Frobenius norm from \cite{PuNoSp2024}.

The general bounds from \ifshort{Section~\ref{SgenBounds}}{\cite{Sp2024}} involve the norm of the full dimensional stochastic term
\( \DFM \, \IFT_{\GPT}^{-1} \nabla \zeta \).
We will use the bound \( \IFT_{\GPT} \geq \dmax^{-2} \DFM^{2} + \GPT^{2} \)
with \( \DFM^{2} = \blk\{ \DPN^{2}, \Id_{\imav}, \muA^{2} \Id_{\KS} \} \),
\( \GPT^{2} = \blk\{ \GP^{2},0,\GPKS^{2} \} \),
and the structure \( \nabla \zeta = (0,\errZv,\muA \, \errKS) \) of the score vector.
Condition \nameref{EU2ref} implies with \( \BBH_{\DFM} \eqdef \Var(\DFM \, \IFT_{\GPT}^{-1} \nabla \zeta) \)
on a random set \( \Omega(\xx) \) with \( \P(\Omega(\xx)) \geq 1 - 3\ex^{-\xx} \)
\begin{EQA}[c]
	\| \DFM \, \IFT_{\GPT}^{-1} \nabla \zeta \|
	\leq 
	\rr_{\DFM}
	\eqdef
	\zq(\BBH_{\DFM},\xx)
	\leq 
	\sqrt{\tr \BBH_{\DFM}} + \sqrt{2 \xx \, \| \BBH_{\DFM} \|} \,\, .
\label{vhyvw323sdwe4w3rdesf23EO}
\end{EQA}
Some upper bounds on the value \( \dimAfull_{\DFM} \eqdef \tr \Var(\DFM \, \IFT_{\GPT}^{-1} \nabla \zeta) = \tr \BBH_{\DFM} \)
will be given in Section~\ref{SscoreEO} later.
The inequality \( \IFT_{\GPT} \geq \dmax^{-2} \DFM^{2} + \GPT^{2} \) yields
\begin{EQA}[rcl]
	\| \DFM \, \IFT_{\GPT}^{-1} \nabla \zeta \|^{2}
	& \leq &
	\dmax^{2} \| \errZv \|^{2} 
	+ \| (\dmax^{-2} \Id_{\KS} + \muA^{-2}\GPKS^{2})^{-1} \errKS \|^{2}  
	\\
	&=&
	\dmax^{2} \| \errZv \|^{2} 
	+ \sum_{\mm=1}^{\dimq}
	\| (\dmax^{-2} \Id_{\dimp} + \muA^{-2}\GPKS_{\mm}^{2})^{-1} \errKS_{\mm} \|^{2}
	\, .
\label{uvobkjyy643enhfgvtygyjyef}
\end{EQA}
Similarly 
\( \| \DFM_{\nupv\nupv} \, \IFT_{\GPT,\nupv\nupv}^{-1} \nabla_{\nupv} \zeta \|^{2} 
\leq \| (\dmax^{-2} \Id_{\KS} + \muA^{-2}\GPKS^{2})^{-1} \errKS \|^{2} \).

 
\Section{Loss and risk of estimation\iffourG{. 3S bounds}{}}
Here we collect our main results about the behavior of \( \tilde{\tarpv}_{\GPT} \), its loss and risk.
The results will be derived from the general statements in \ifshort{\Chname \ref{SgenBounds}}{\cite{Sp2024}}.
We use \( \nabla \zeta = (0,\errZv,\muA \, \errKS) \) and \( \GPT^{2} = \blk\{ \GP^{2},0,\GPKS^{2} \} \).

\begin{theorem}
\label{TsemieffEO}
Assume \nameref{EU2ref} for \( \nabla  \zeta = (0,\errZv,\muA \, \errKS) \).
Let \( \DPN^{2} \) satisfy \eqref{tfeuhfuwr365vdyge34def} and let,
with \( \nEO = \lambda_{\min}(\DPN^{2}) \) and \( \REO = \eEO \muA \sqrt{\nEO} \) for \( \eEO \leq 1/10 \), 
the set \( \Upsd \) from \eqref{rughur47478ytfg84g6gh} contain \( \prmtvs \).
For \( \DFM^{2} \) from \eqref{iwdfhw787e3urne6fyh3o}, \( \rr_{\DFM} \) from \eqref{vhyvw323sdwe4w3rdesf23EO}, 
\( \dltwu_{3} \) from \eqref{gdjdgredtgebegxgEO}, 
and \( \bias_{\DFM} \eqdef \| \DFM \IFT_{\GPT}^{-1} \, \GPT^{2} \prmtvs \| \),
assume 
\begin{EQA}[c]
	\REO \geq \frac{3}{2} (\rr_{\DFM} \vee \bias_{\DFM}) \, ,
	\qquad
	\dmax^{2} \dltwu_{3} \, (\rr_{\DFM} \vee \bias_{\DFM}) < \frac{4}{9} 
	\, 
\label{y7jekc72hwjmgihy3h2ftt}
\end{EQA}
with \( \dmax = 2 \).
Then the following statements hold for any linear mapping \( \QP \) on \( \R^{\dimp} \):
\begin{enumerate}
\item Fisher expansion: on \(\Omega(\xx) \)
\begin{EQA}
\label{vdud6ehry3hrf67EO}
	\| \QP \{ \tilde{\tarpv}_{\GPT} - \tarpvs_{\GPT} - (\IFT_{\GPT}^{-1} \nabla \zeta)_{\tarpv} \} \|
	& \leq &
	\| \QP \, \DPN^{-1} \| \, \frac{3 \dmax^{2} \dltwu_{3}}{4} \, 
	\| \DFM \IFT_{\GPT}^{-1} \nabla \zeta \|^{2} \, ,
\label{cdugec4e435dtfce525rs}
\end{EQA}
where \( \| \DFM \IFT_{\GPT}^{-1} \nabla \zeta \|^{2} \) satisfies 
\eqref{vhyvw323sdwe4w3rdesf23EO} and \eqref{uvobkjyy643enhfgvtygyjyef}.

\item Wilks expansion: on \(\Omega(\xx) \),
\( \LLp_{\GP}(\thetav) = \sup_{\etav} \LL_{\GP}(\thetav,\etav) \) satisfies
\begin{EQA}
	\Bigl| 
		2 \LLp_{\GPT}(\tilde{\tarpv}_{\GPT}) - 2 \LLp_{\GPT}(\tarpvs_{\GPT}) - \bigl\| \xivr_{\GPT} \bigr\|^{2} 
	\Bigr|
	& \leq &
	\frac{\dltwu_{3}}{2} 
	\Bigl( \| \DFM \IFT_{\GPT}^{-1} \nabla \zeta \|^{3} + \| \DFM_{\nupv\nupv} \IFT_{\GPT,\nupv\nupv}^{-1} \scorem{\nupv} \zeta \|^{3} \Bigr) \, ,
	\qquad
\label{vdud6ehry3hrf67ruy233EO}
\end{EQA}
where with \( \IFTb_{\GPT,\thetav\thetav} \) from \eqref{kf8r484r5iu8t7t7eurggEO}
\begin{EQA}
	\xivr_{\GPT}
	& \eqdef &
	\IFTb_{\GPT,\thetav\thetav}^{1/2} \, (\IFT_{\GPT}^{-1} \nabla \zeta)_{\thetav}
	\, .
\label{d7fje3gv7hre3je8gb53jEO}
\end{EQA}

\item Estimation bias: 
\begin{EQA}
	\bigl\| \QP \, \bigl\{ \tarpvs_{\GPT} - \tarpvs + (\IFT_{\GPT}^{-1} \, \GPT^{2} \prmtvs)_{\tarpv} \bigr\} \bigr\| 
	& \leq & 
	\| \QP \, \DPN^{-1} \| \, \frac{3 \dmax^{2} \dltwu_{3}}{4} \, \bias_{\DFM}^{2} \, .
	\qquad
\label{0mkvhgjnrw3dfwe3u8gEO}
\end{EQA}

\item PAC loss bound: on \(\Omega(\xx) \)
\begin{EQA}
	\!\!\!\!\!\!
	\bigl\| 
		\QP \bigl\{ \tilde{\tarpv}_{\GPT} - \tarpvs 
			- (\IFT_{\GPT}^{-1} \nabla \zeta)_{\tarpv} + (\IFT_{\GPT}^{-1} \, \GPT^{2} \prmtvs)_{\tarpv} 
		\bigr\} 
	\bigr\|
	& \leq &
	\| \QP \, \DPN^{-1} \| \, \frac{3 \dmax^{2} \dltwu_{3}}{4} \, 
	(\| \DFM \IFT_{\GPT}^{-1} \nabla \zeta \|^{2} + \bias_{\DFM}^{2}) \, .
	\qquad \qquad
\label{rG1r1d3GbvGbseEO}
\end{EQA}

\item \( L_{2} \)-norm risk bound:
\begin{EQA}
	&& \nquad
	\E \, \bigl\{ \| \QP (\tilde{\tarpv}_{\GPT} - \tarpvs) \| \Ind_{\Omega(\xx)} \bigr\} 
	\leq 
	\riskt_{\QP}^{1/2}
	+ \| \QP \, \DPN^{-1} \| \, \frac{3 \dmax^{2} \dltwu_{3}}{4} \, \bigl( \dimAfull_{\DFM} + \bias_{\DFM}^{2} \bigr) 
	\, ,
\label{EtuGus11md3GEO}
\end{EQA}
where \( \dimAfull_{\DFM} = \E \| \DFM \IFT_{\GPT}^{-1} \nabla \zeta \|^{2} \), and
\begin{EQA}[rcl]
	\riskt_{\QP} 
	& \eqdef &
	\E \bigl\{ 
		\bigl\| \QP (\IFT_{\GPT}^{-1} \nabla \zeta)_{\tarpv} - \QP (\IFT_{\GPT}^{-1} \, \GPT^{2} \prmtvs)_{\tarpv} \bigr\|^{2} \Ind_{\Omega(\xx)} 
	\bigr\} 
	\\
	& \leq &
	\tr \Var\bigl\{ \QP (\IFT_{\GPT}^{-1} \nabla \zeta)_{\tarpv} \bigr\}  
	+ \| \QP \, (\IFT_{\GPT}^{-1} \, \GPT^{2} \prmtvs)_{\tarpv} \|^{2} \, .
	\qquad
\label{t6skjid7ehfcr56tegfdEO}
\end{EQA}

\item Squared risk bound:
if 
\begin{EQA}[c]
	\E \bigl\{ \| \DFM \IFT_{\GP}^{-1} \nabla \zeta \|^{4} \Ind_{\Omega(\xx)} \bigr\} \leq \CONSTi_{4}^{2} \, \dimAfull_{\DFM}^{2} \, ,
\label{6hjdfv8e6hyefyeheew7skGPT}
\end{EQA}
and a constant \( \alp_{\QP} \) ensures 
\begin{EQ}[rcl]
	\alp_{\QP}
	& \eqdef & 
	\frac{\| \QP \, \DPN^{-1} \| \, (3/4) \dmax^{2} \dltwu_{3} \, (\CONSTi_{4} \, \dimAfull_{\DFM} + \bias_{\DFM}^{2})} {\sqrt{\riskt_{\QP}}} 
	< 1 \, ,
\label{6dhx6whcuydsds655srEO}
\end{EQ}
then 
\begin{EQA}
	(1 - \alp_{\QP})^{2} \riskt_{\QP}
	\leq 
	\E \bigl\{ \| \QP (\tilde{\tarpv}_{\GPT} - \tarpvs) \|^{2} \Ind_{\Omega(\xx)} \bigr\}
	& \leq &
	(1 + \alp_{\QP})^{2} \riskt_{\QP} \, .
\label{EQtuGmstrVEQtGEO}
\end{EQA}

\end{enumerate}
\end{theorem}

A great benefit of the suggested approach is that the scope of the results presented in Theorem~\ref{TsemieffEO}
is obtained as a straightforward corollary of the general theory developed in \cite{Sp2024}.
We only need to check four general conditions.
This check will be done later in Section~\ref{SEiOtools}. 
Fisher expansion \eqref{cdugec4e435dtfce525rs} plays the central role, it yields 
the loss and risk bounds provided that the remainder of this expansion is smaller in order than the leading term.
This condition will be discussed in the next section.
The finite sample expansion \eqref{cdugec4e435dtfce525rs} can be used for inference 
and for studying the limiting behavior of \( \tilde{\tarpv}_{\GPT} \).
It is very useful for proving bootstrap validity for many resampling procedures; see e.g. \cite{SpZh2014} and \cite{NSU2017} and references therein.
Similarly, Wilks expansion \eqref{vdud6ehry3hrf67ruy233EO} is important, e.g., for testing and model selection;
cf. to the classical Wilks theorem from \cite{FaZh2001,BoMa2011}.
A closed-form representation for the stochastic part \( (\IFT_{\GPT}^{-1} \nabla \zeta)_{\tarpv} \) and the bias
\( (\IFT_{\GPT}^{-1} \, \GPT^{2} \prmtvs)_{\tarpv} \) will be given in Section~\ref{SEiOtools}.
Section~\ref{ScritdimEO} discusses the variance of the semiparametric efficient score \( \xivr_{\GPT} \).
In particular, we derive finite-sample versions of semiparametric efficiency results for the considered setup.
The bias term \( (\IFT_{\GPT}^{-1} \, \GPT^{2} \prmtvs)_{\tarpv} \) is analysed in Section~\ref{SEiObias}
and detailed for the important cases of ridge and truncation penalties.

\Section{Stochastic term. Effective and critical dimension}
\label{ScritdimEO}
This section discusses the impact and applicability of the results of Theorem~\ref{TsemieffEO} in terms of the relation 
between the effective sample size \( \nEO \) from \eqref{tt7t7yrs3esfhiyfsftgh}
and the effective dimension \( \dimA_{\GPT} \) of the problem.
The main message of the analysis is that under the critical dimension condition
\( \dimA_{\GPT} \ll \nEO \), the target parameter \( \tarpv \) will be estimated 
with essentially the same accuracy as in the case of a precisely known operator \( \KSs \).

The full-dimensional parameter \( \prmtv \) is composed of
the \( \dimp \)-dimensional target parameter \( \tarpv \), 
\( \dimq \)-dimensional image \( \imav \),
and \( \dimp \, \dimq \)-dimensional operator \( \KS \) yielding the full parameter dimension
\( \dimfull = \dimp + \dimq + \dimp \, \dimq \).
Now we discuss the full effective dimension \( \dimAfull_{\GPT} = \tr \BBB_{\GPT} \) shown in \eqref{vhyvw323sdwe4w3rdesf23EO}.
The next result explains how the regularization term \( \| \KS \|_{\GPKS}^{2}/2 \) helps to reduce the value 
\( \dimAfull_{\GPT} \) compared with \( \dimfull \).

\begin{proposition}
\label{LdimfullEO}
For \( \nabla \zeta = (0,\errZv,\muA \, \errKS) \) with
\( \errZv = \Zv - \E \Zv \) and \( \errKS = \muA (\hKS - \KS) \), it holds
\begin{EQA}[c]
	\dimAfull_{\DFM}
	=
	\tr \Var(\DFM \, \IFT_{\GPT}^{-1} \nabla \zeta)
	\leq 
	\dimA_{\imav} + \dimA_{\KS} 
	\, ,
\label{ydfhuwir8ufug87g73jdeEO}
\end{EQA}
where
\begin{EQA}[c]
	\dimA_{\imav}
	\leq 
	\dmax^{4} \tr \Var\bigl( \errZv \bigr) ,
	\qquad
	\dimA_{\KS} 
	\leq 
	\dmax^{4} \tr \Var\bigl\{ (\Id_{\KS} + \dmax^{2} \muA^{-2} \GPKS^{2})^{-1/2} \errKS \bigr\} 
	\, .
\end{EQA}
For the homogeneous observation noise with \( \Var(\errZ_{\mm}) \leq \sigma_{\errZ}^{2} \),
\( \mm \leq \dimq \), it holds \( \dimA_{\KS} \leq \dmax^{4} \dimq \sigma_{\errZ}^{2} \).
Moreover, homogeneous operator noise \( \errKS = (\errKS_{\mm}) \) with 
\( \Var(\errKS_{\mm}) \leq \CONSTi_{\errKS} \Id_{\dimp} \) ensures
\begin{EQA}
	\dimA_{\KS}
	& \leq &
	\dmax^{4} \CONSTi_{\errKS} \sum_{\mm=1}^{\dimq} 
	\tr \Var\bigl\{ (\Id_{\dimp} + \dmax^{2} \muA^{-2} \GPKS_{\mm}^{2})^{-1} \bigr\}
	\, .
\label{u6r664dr535e5tu6868yug}
\end{EQA}
\end{proposition}

A practically relevant special case corresponds to a truncation penalty
\( \GPKS_{\mm}^{2} = \muA^{2} \gpks_{\mm}^{2} \Id_{\dimp} \) with \( \gpks_{\mm} = 0 \) 
for \( \mm \leq \dimm \) and \( \gpks_{\mm} = \infty \) for \( \mm > \dimm \).
Then
\begin{EQA}
	\dimA_{\KS}
	& \leq &
	\dmax^{4} \CONSTi_{\errKS} \, \dimp \, \dimm
	\, .
\label{u6r664dr535e5tu6868yugt}
\end{EQA}
Such a penalization replaces \( \dimq \) with \( \dimm \) and \( \dimp \, \dimq \) with 
\( \dimp \, \dimm \).
Theorem~\ref{TsemieffEO} requires \( \dimA_{\DFM} \ll \muA \sqrt{\nEO} \); see \eqref{y7jekc72hwjmgihy3h2ftt};
or, by \eqref{gdjdgredtgebegxgEO}
\begin{EQA}
	\nEO^{-1} \muA^{-2} \, \dimAfull_{\GPT} 
	\approx 
	\nEO^{-1} \muA^{-2} \, \dimp \, \mM
	& \ll & 
	1 .
\label{jef9fe8e47fu32oeoggitir4g}
\end{EQA}
Without penalization of the operator \( \KS \), it leads back to 
the full-dimensional condition \( \dimp \, \dimq \ll \muA^{2} \nEO \).

Now we discuss the leading term \( (\IFT_{\GPT}^{-1} \nabla \zeta)_{\tarpv} \) in Fisher expansion \eqref{cdugec4e435dtfce525rs}.
The next result shows that the variance of \( (\IFT_{\GPT}^{-1} \nabla \zeta)_{\tarpv} \)
is the same in order as for the ideal model with the precisely known operator \( \KSs \).

\begin{proposition}
\label{Leffscore}
Assume 
\( \| \Var(\errZv) \| \leq \CONSTi_{\errZ} \) and
\begin{EQA}[c]
	\bigl\| \Var\bigl\{ (\Id_{\KS} + \muA^{-2} \GPKS^{2})^{-1} \errKS \bigr\} \bigr\|
	\leq 
	\CONSTi_{\errKSe}
	\, .
\label{x6dc4c32c2ctv9necgsdr}
\end{EQA} 
Then 
\begin{EQA}[rcl]
	\Var \bigl( \IFT_{\GPT}^{-1} \nabla \zeta \bigr)_{\tarpv}
	& \leq &
	2 \dmax^{2} (\CONSTi_{\errZ} + 4 \CONSTi_{\errKSe} \eEO^{2}) 
	(\dmax^{-2} \DPN^{2} + \GP^{2})^{-1}
	\, ,
\label{du7dyfu3ijhrk398rffi33}
\end{EQA}
and for any \( \QP \) linear mapping \( \QP \) on \( \R^{\dimp} \)
\begin{EQA}[c]
	\E \bigl\| \QP \bigl( \IFT_{\GPT}^{-1} \nabla \zeta \bigr)_{\tarpv} \bigr\|^{2}
	\leq 
	2 \dmax^{2} (\CONSTi_{\errZ} + 4 \CONSTi_{\errKSe} \eEO^{2}) \tr \bigl\{ \QP (\dmax^{-2} \DPN^{2} + \GP^{2})^{-1} \QP^{\T} \bigr\}
	\, .
\label{judfvy6h34ibte35ugit}
\end{EQA}
Moreover, with \( \xivr_{\GPT} \) from \eqref{d7fje3gv7hre3je8gb53jEO}
\begin{EQA}[rcl]
	\E \| \xivr_{\GPT} \|^{2}
	&=&
	\tr \Var ( \xivr_{\GPT} )
	\leq 
	2 (\CONSTi_{\errZ} + 4 \CONSTi_{\errKSe} \eEO^{2}) 
	\tr \bigl\{ (\dmax^{-2} \DPN^{2} + \GP^{2})^{-1} \DPN^{2} \bigr\}
	\, .
\end{EQA}
\end{proposition}

\Section{A bound on the penalization bias and the estimation risk}
\label{SEiObias}
This section provides some details about the penalization bias 
\( \| \QP (\IFT_{\GPT}^{-1} \GPT^{2} \prmtvs)_{\tarpv} \| \).
Our intention is to show that the extension of the parameter space and penalization of the operator 
does not introduce any substantial bias for the target parameter.
We may use that the image vector \( \imav \) is not penalized 
and the full dimensional penalization matrix \( \GPT^{2} \) is of block-structure w.r.t. 
\( \tarpv \) and \( \KS \) parameters.
This enables us to separate the bias caused by \( \tarpv \)-penalty \( \| \GP \tarpv \|^{2} \)
and the bias due to operator penalty \( \| \KS \|_{\GPKS}^{2} \).

\begin{lemma}
\label{PsemiriskEO}
Let \( \GPKS^{2} = \blk\{ \GPKS_{1}^{2}, \ldots, \GPKS_{\dimq}^{2} \} \).
Then for any linear mapping \( \QP \) on \( \R^{\dimp} \)
\begin{EQA}
	\| \QP (\IFT_{\GPT}^{-1} \GPT^{2} \prmtvs)_{\tarpv} \|
	& = &
	\| \QP \, \IFTb_{\GPT,\tarpv\tarpv}^{-1} (\GP^{2} - S_{\GPKS}) \tarpvs \| 
	\leq 
	\| \QP \, \IFTb_{\GPT,\tarpv\tarpv}^{-1} \GP^{2} \tarpvs \|
	+ \| \QP \, \IFTb_{\GPT,\tarpv\tarpv}^{-1} \, S_{\GPKS} \, \tarpvs \| 
\label{vi458gteghe3hwjdxqqwwEO}
\end{EQA}
with \( \IFTb_{\GPT,\tarpv\tarpv} \) from \eqref{kf8r484r5iu8t7t7eurggEO} and
\begin{EQA}
	S_{\GPKS}
	& \eqdef &
	\frac{1}{2} \sum_{\mm=1}^{\dimq} \KSs_{\mm} {\KSs_{\mm}}^{\T} \GPKS_{\mm}^{2} 
	\bigl( \muA^{2} \Id_{\dimp} + \GPKS_{\mm}^{2} - \frac{1}{2} \tarpvs {\tarpvs}^{\T} \bigr)^{-1}
	\, .
\label{kchdhv63bvuedr43bvg0h3}
\end{EQA}
In particular, with \( \GPKS_{\mm}^{2} = \gpks_{\mm}^{2} \Id_{\dimp} \)
\begin{EQA}
	S_{\GPKS}
	& \eqdef &
	\frac{1}{2} 
	\sum_{\mm=1}^{\dimq} \frac{\gpks_{\mm}^{2}}{\muA^{2} + \gpks_{\mm}^{2} - \| \tarpvs \|^{2}/2} \, 
	\KSs_{\mm} {\KSs_{\mm}}^{\T} \, 
	\, .
\label{kchdhv63bfub83kfgibnurr}
\end{EQA}
\end{lemma}

Putting together bound \eqref{du7dyfu3ijhrk398rffi33} on the variance term and the obtained bias bound 
yields risk bound \eqref{EQtuGmstrVEQtGEO} with
\begin{EQA}[rcl]
	\riskt_{\QP}
	&=& \tr \Var\bigl\{ \QP (\IFT_{\GPT}^{-1} \nabla \zeta)_{\tarpv} \bigr\}  
	+ \| \QP (\IFT_{\GPT}^{-1} \, \GPT^{2} \prmtvs)_{\tarpv} \|^{2}
	\\
	& \lesssim &
	\tr\bigl\{ \QP (\dmax^{-2} \DPN^{2} + \GP^{2})^{-1} \QP^{\T} \bigr\} 
	+ \| \QP \, \IFTb_{\GPT,\tarpv\tarpv}^{-1} \GP^{2} \tarpvs \|^{2}
	+ \| \QP \, \IFTb_{\GPT,\tarpv\tarpv}^{-1} \, S_{\GPKS} \, \tarpvs \|^{2}
	\, .
\label{jcvyyvb7n99hj3jdfuycv}
\end{EQA}
Moreover, a proper choice of the operator penalty \( \GPKS \) enables us to ignore the last term reflecting the operator penalization. 

\Chapter{Risk bounds and choice of penalty}
\label{SpenaltyEO}
This section addresses an important question of choosing a proper penalty
to obtain optimal estimation accuracy.
It appears that the answer is different depending on the relation between 
operator regularity and smoothness of the signal.
For a very smooth operator, a simple ridge penalization does a good job.
In the other cases, one should apply one or another model reduction technique.
We discuss the spectral cut-off method in combination with the approximation spaces setup. 

\Section{Ridge penalty and a smooth operator}
An important example of penalty choice is a ridge penalty
\( \GP^{2} = \gp^{2} \Id_{\dimp} \).
It is basis and coordinate free and enforces the ``benign overfitting'' phenomenon
in high-dimensional regression;
see \cite{baLoLu2020,ChMo2022,NoPuSp2024} and references therein.
Later we consider the estimation problem with \( \QP = \Id_{\dimp} \).

Our study reveals an interesting phase transition effect. 
The use of ridge penalty leads to accurate results in the situation 
when the operator \( {\KSs}^{\T} \KSs \) is ``more regular'' than the signal \( \tarpvs \).
In this case, ridge penalization enforces nearly the same effect as 
a spectral cut-off method.
Moreover, with a properly chosen ridge parameter \( \gp^{2} \), one can achieve 
the bias-variance trade-off in estimation of the target parameter.
Error in the operator can be ignored and the final accuracy is as in the case of linear inverse problem with
the known operator \( \KSs \). 
The situation changes dramatically if the operator \( {\KSs}^{\T} \KSs \) is not smoother than the signal. 
This includes the case of a direct problem when the conditional number of the operator 
\( \KSs \) is bounded by a fixed constant. 
It is well known that the ridge penalty is not efficient in this case, 
and model reduction technique should be applied. 

In the rest of this section, we assume a smooth operator \( \KSs \)
and focus on ridge penalization for the target parameter \( \tarpv \) only, 
that is, \( \GPKS \equiv 0 \); cf. \cite{NoPuSp2024}.
Our study includes the case with \( \dimp = \dimq = \infty \).
Introduce the ordered eigenvalues \( \nEO_{1} \geq \nEO_{2} \geq \ldots \geq \nEO_{\dimp} \) of \( \DPN^{2} = {\KSs}^{\T} \KSs \).
By \eqref{EQtuGmstrVEQtGEO} and \eqref{t6skjid7ehfcr56tegfdEO} of Theorem~\ref{TsemieffEO},
the squared risk of \( \tilde{\tarpv}_{\GPT} \) can be approximated by \( \riskt \) with
\begin{EQA}[rcccl]
	\riskt
	&=& \tr \Var\bigl\{ (\IFT_{\GPT}^{-1} \nabla \zeta)_{\tarpv} \bigr\}  
	+ \| (\IFT_{\GPT}^{-1} \, \GPT^{2} \prmtvs)_{\tarpv} \|^{2} 
	& \lesssim &
	\tr\bigl\{ (\dmax^{-2} \DPN^{2} + \gp^{2} \Id_{\dimp})^{-1} \bigr\} 
	+ \gp^{4} \| \IFTb_{\GPT,\tarpv\tarpv}^{-1} \, \tarpvs \|^{2}
	\, .
\label{jcvyyvb7n99hj3jdfuycv}
\end{EQA}
Introduce an operator \( \BBH \eqdef \bigl( \Id_{\dimp} + \dmax^{-2} \gp^{-2} \DPN^{2} \bigr)^{-1} \).
It can be viewed as an approximation of the projector
in \( \R^{\dimp} \) on the subspace defined by the inequality
\( \DPN^{2} \geq \dmax^{2} \gp^{2} \Id_{\dimp} \).
This subspace is spanned by the eigenvectors \( \ev_{j} \)
corresponding to \( \nEO_{j} \geq \dmax^{2} \gp^{2} \).
The dimension of this space is given by 
\begin{EQA}[c]
	\jJ_{\gp} \eqdef \max\{ j \colon \nEO_{j} \geq \dmax^{2} \gp^{2} \} .
\label{dufhiw3jeio26vte3hj}
\end{EQA}
This relation can be inverted: 
given \( \jJ \), the corresponding ridge factor \( \gp^{2} \) can be given by 
\( \gp^{2} = \dmax^{-2} \nEO_{\jJ} \).

\begin{lemma}
\label{LvarEO}
Let \( \GP^{2} = \gp^{2} \Id_{\dimp} \) and \( \nEO_{1} \geq \nEO_{2} \geq \ldots \geq \nEO_{\dimp} \) 
be the ordered eigenvalues of \( \DPN^{2} \) while \( \ev_{j} \) be the corresponding eigenvectors.
With \( \jJ = \jJ_{\gp} \) from \eqref{dufhiw3jeio26vte3hj},
\begin{EQA}[rcl]
	\tr\bigl\{ (\dmax^{-2} \DPN^{2} + \GP^{2})^{-1} \bigr\}
	& \leq &
	\biggl( \dmax^{4} \sum_{j=1}^{\jJ} \frac{1}{\nEO_{j}} + \gp^{-4} \sum_{j=\jJ+1}^{\dimp} \nEO_{j} \biggr)
	\, ,
	\\
	\gp^{4} \| \IFTb_{\GPT,\tarpv\tarpv}^{-1} \, \tarpvs \|^{2}
	& \leq &
	\bigl\| \BBH \tarpvs \bigr\|^{2}
	=
	\sum_{j=1}^{\dimp} \frac{1}{(\dmax^{-2} \gp^{-2} \nEO_{j} + 1)^{2}} \langle \tarpvs, \ev_{j} \rangle^{2}
	\, .
\end{EQA}
\end{lemma}
\begin{proof}
Lemma~\ref{Leffscore} implies
\begin{EQA}[rcl]
	&& \hspace{-18pt}
	\tr \Var \bigl( \IFT_{\GPT}^{-1} \nabla \zeta \bigr)_{\tarpv}
	\leq 
	2 (\CONSTi_{\errZ} + 4 \CONSTi_{\errKSe} \eEO^{2}) 
	\tr \bigl\{ (\dmax^{-2} \DPN^{2} + \gp^{2} \Id_{\dimp})^{-1} \DPN^{2} (\dmax^{-2} \DPN^{2} + \gp^{2} \Id_{\dimp})^{-1} \bigr\}
	\\
	&=&
	2 (\CONSTi_{\errZ} + 4 \CONSTi_{\errKSe} \eEO^{2})  
	\sum_{j=1}^{\dimp} \frac{\nEO_{j}}{(\dmax^{-2} \nEO_{j} + \gp^{2})^{2}}
	\leq 
	2 (\CONSTi_{\errZ} + 4 \CONSTi_{\errKSe} \eEO^{2}) 
	\biggl( \dmax^{4} \sum_{j=1}^{\jJ} \nEO_{j}^{-1} + \gp^{-4} \sum_{j=\jJ+1}^{\dimp} \nEO_{j} \biggr)
	\, ,
\end{EQA}
and the first result follows.
The second one follows from Lemma~\ref{LIFTDFMEO}.
\end{proof}

The results can be further simplified if the values \( \nEO_{j} \) decay polynomially:
\( \nEO_{j} \asymp \nEO_{1} \, j^{-2s} \) for \( s > 1/2 \).
Then for any \( \jJ \), we may use
\begin{EQA}[c]
	\sum_{j=1}^{\jJ} \frac{1}{\nEO_{j}} \leq \CONSTi_{1} \, \frac{\jJ}{\nEO_{\jJ}}  \, ,
	\qquad
	\sum_{j=\jJ+1}^{\dimp} \nEO_{j} \leq \CONSTi_{2} \, \jJ \nEO_{\jJ}  \, .
\label{udvcghue3j3vc654rw2y}
\end{EQA}

\begin{proposition}
\label{PvarEO}
Let \( \GP^{2} = \gp^{2} \Id_{\dimp} \).
Assume \eqref{udvcghue3j3vc654rw2y}
and \( \jJ = \jJ_{\gp} \); see \eqref{dufhiw3jeio26vte3hj}.
Then
\begin{EQA}[rcl]
	\tr\bigl\{ (\dmax^{-2} \DPN^{2} + \gp^{2} \Id_{\dimp})^{-1} \bigr\}
	& \leq &
	\dmax^{4} 
	(\CONSTi_{1} + \CONSTi_{2}) \frac{\jJ}{\nEO_{\jJ}}
	\, .
	\qquad
\label{ckjivhiedi93wikvuy}
\end{EQA}
Further, let \( \tarpvs \) satisfy the smoothness condition
\begin{EQA}[c]
	\sum_{j=1}^{\dimp} \cgp_{j}^{2} \langle \tarpvs, \ev_{j} \rangle^{2}
	\leq 
	1,
\label{dcyhvc55t8tg764e33rf}
\end{EQA}
where \( \cgp_{j} \) grow while \( \cgp_{j} \nEO_{j} \) decrease with \( j \).
With \( \jJ = \jJ_{\gp} \) from \eqref{dufhiw3jeio26vte3hj}, it holds
\begin{EQA}[rcl]
	\bigl\| \bigl( \dmax^{-2} \gp^{-2} \DPN^{2} + \Id_{\dimp} \bigr)^{-1} \tarpvs \bigr\|^{2}
	& \leq &
	\max_{j \leq \dimp} \frac{\cgp_{j}^{-2}}{\dmax^{-2} \gp^{-2} \nEO_{j} + 1}
	\leq 
	\frac{1}{\cgp_{\jJ}^{2}}
	\, 
\label{ckjivhiedi93wikvuy2}
\end{EQA}
yielding
\begin{EQA}
	\riskt
	& \lesssim &
	\frac{\jJ}{\nEO_{\jJ}} + \frac{1}{\cgp_{\jJ}^{2}} 
	\, .
\label{c8f78vtv4c4tehfivef}
\end{EQA}
\end{proposition}

\begin{proof}
Statement \eqref{ckjivhiedi93wikvuy} follows Lemma~\ref{LvarEO} and \eqref{udvcghue3j3vc654rw2y}.
Further, 
\begin{EQA}[rcl]
	\| \BBH \tarpvs \|^{2}
	&=&
	\sum_{j=1}^{\dimp} \frac{1}{(\dmax^{-2} \gp^{-2} \nEO_{j} + 1)^{2}} \langle \tarpvs, \ev_{j} \rangle^{2}
	=
	\sum_{j=1}^{\dimp} \frac{\cgp_{j}^{-2}}{(\dmax^{-2} \gp^{-2} \nEO_{j} + 1)^{2}} \, \cgp_{j}^{2} \langle \tarpvs, \ev_{j} \rangle^{2}
	\\
	& \leq &
	\max_{j \leq \dimp} \frac{\cgp_{j}^{-2}}{(\dmax^{-2} \gp^{-2} \nEO_{j} + 1)^{2}}
	\sum_{j=1}^{\dimp} \cgp_{j}^{2} \langle \tarpvs, \ev_{j} \rangle^{2}
	\leq 
	\biggl\{ \min_{j \leq \dimp} (\dmax^{-2} \gp^{-2} \nEO_{j} \cgp_{j} + \cgp_{j}) \biggr\}^{-2} \, 
	\, .
\end{EQA}
As the values \( \nEO_{j} \cgp_{j} \) decrease and \( \cgp_{j} \) increase with \( j \),
it holds for any \( \jJ \)
\begin{EQA}[c]
	\min_{j \leq \jJ} (\dmax^{-2} \gp^{-2} \nEO_{j} \cgp_{j} + \cgp_{j})
	\geq 
	\dmax^{-2} \gp^{-2} \nEO_{\jJ} \cgp_{\jJ} \, ,
	\quad
	\min_{j > \jJ} (\dmax^{-2} \gp^{-2} \nEO_{j} \cgp_{j} + \cgp_{j})
	\geq 
	\cgp_{\jJ+1}
	\, .
	\qquad
\end{EQA}
Therefore, the choice of \( \jJ \) by \( \dmax^{-2} \gp^{-2} \nEO_{\jJ} \approx 1 \)
yields
\begin{EQA}[c]
	\min_{j \leq \dimp} (\dmax^{-2} \gp^{-2} \nEO_{j} \cgp_{j} + \cgp_{j})
	\geq 
	\cgp_{\jJ}
	\, ,
\end{EQA}
and \eqref{ckjivhiedi93wikvuy2} follows as well.
\end{proof}

Bound \eqref{ckjivhiedi93wikvuy2} requires that the values \( \cgp_{j} \) 
grow while \( \nEO_{j} \cgp_{j} \) decrease.
The latter means that the operator \( {\KSs}^{\T} \KSs \) is more regular than the signal \( \tarpvs \).
Now, consider the opposite case.
For instance, let the eigenvalues \( \nEO_{j} \) of \( \DPN^{2} \) decrease slowly with 
\( j \) or remain significantly positive. 
Then inversion of \eqref{dufhiw3jeio26vte3hj} for small \( \gp^{2} \) can be 
problematic. 
The corresponding cut-off index \( \jJ \) will be too large and cannot be used 
to mimic the bias-variance trade-off.
As a consequence, a ridge penalization is not efficient in such situations.
The following section explains how this situation can be handled using the spectral cut-off method.

\Section{Truncation penalty and approximation spaces}
This section introduces and studies the special class of truncation penalties.
We start with the spectral cut-off method, then we consider the approximation spaces 
setup.

\Subsection{Spectral cut-off}
Spectral cut-off is a standard tool in model reduction for inverse problems. 
It assumes that the eigenvector decomposition of \( \DPN^{2} = {\KSs}^{\T} \KSs \) 
with ordered eigenvalues \( \nEO_{1} \geq \nEO_{2} \geq \ldots \geq \nEO_{\jJ} \) and the corresponding eigenvectors
\( \ev_{j} \) are available.
By \( \Proj_{\jJ} \) we denote the canonical projector in the source space \( \R^{\dimp} \) onto the first \( \jJ \) coordinates \( \ev_{j} \).
Also, we assume that smoothness properties of the unknown source signal \( \tarpvs \) can be given in the canonical basis 
formed by these eigenvectors \( \ev_{j} \); see \eqref{dcyhvc55t8tg764e33rf}.
For a fixed cut-off parameter \( \jJ \), consider the truncation penalty \( \GP_{\jJ}^{2} = \diag\{ \GP_{1}^{2},\dots,\GP_{\dimp}^{2} \} \)
with \( \gp_{j}^{2} = 0 \) for \( j \leq \jJ \) and \( \gp_{j}^{2} = \infty \) for \( j > \jJ \).
Effectively, this penalty enforces \( \tilde{\tarp}_{\GPT,j} = \tarps_{\GPT,j} = 0 \) for \( j > \jJ \).

\begin{proposition}
\label{LvarEOspect}
Let \( \nEO_{1} \geq \nEO_{2} \geq \ldots \geq \nEO_{\dimp} \) 
be the ordered eigenvalues of \( \DPN^{2} \).
Let also \( \GP_{\jJ}^{2} \) be the spectral cut-off penalty at \( \jJ \). 
Then
\begin{EQA}[rcl]
	\tr\bigl\{ (\dmax^{-2} \DPN^{2} + \GP^{2})^{-1} \bigr\}
	& = &
	\dmax^{4} 
	\sum_{j=1}^{\jJ} \frac{1}{\nEO_{j}} 
	\, ,
	\\
	\| (\IFT_{\GPT}^{-1} \, \GPT^{2} \prmtvs)_{\tarpv} \|^{2}
	& = &
	\| (\Id_{\dimp} - \Proj_{\jJ}) \tarpvs \|^{2}
	\, 
\end{EQA}
yielding
\begin{EQA}[c]
	\riskt
	\lesssim 
	\dmax^{4} \sum_{j=1}^{\jJ} \frac{1}{\nEO_{j}} + \| (\Id_{\dimp} - \Proj_{\jJ}) \tarpvs \|^{2}
	\, .
\end{EQA}
Under \eqref{udvcghue3j3vc654rw2y} and \eqref{dcyhvc55t8tg764e33rf}, bound \eqref{c8f78vtv4c4tehfivef} applies. 
\end{proposition}

In the contrary to ridge penalization, 
the penalty coefficients \( \gp_{j}^{2} \) vanish for \( j \leq \jJ \).
This improves the bound on the penalization bias, the assumption of decay of \( \cgp_{j}^{2} \nEO_{j} \) is not required.

A proper choice of the cut-off parameter \( \jJ \) ensures a nearly optimal accuracy of estimation. 
This approach is widely used in linear inverse problems.
However, availability of the SVD for \( \DPN^{2} \) is a severe limitation, especially if the operator \( \KSs \) is unknown.
This issue can be resolved using the \emph{approximation space} approach.

\Subsection{Approximation spaces and truncation penalties}
This setup assumes specific basises in the parameter space \( \R^{\dimp} \) and the image space \( \R^{\dimq} \)
suitable for describing the smoothness properties of \( \tarpvs \) and regularity of \( \KSs \) simultaneously.
Without loss of generality, we apply the canonical basis in the source space \( \R^{\dimp} \) and the image space \( \R^{\dimq} \).
By \( \Proj_{\jJ} \) we denote the canonical projector in the source space \( \R^{\dimp} \) onto the first \( \jJ \) coordinates.
Similarly, \( \proj_{\mM} \) is the canonical projector in the image space \( \R^{\dimq} \).

Consider a truncation penalty \( \GP^{2} = \diag\{ \GP_{1}^{2},\dots,\GP_{\dimp}^{2} \} \) with 
 \( \GP_{j}^{2} = 0 \) for \( j \leq \jJ \) and
\( \GP_{j}^{2} = \infty \) for \( j > \jJ \).
Effectively, this penalty enforces \( \tilde{\tarp}_{\GPT,j} = \tarps_{\GPT,j} = 0 \) for \( j > \jJ \).
Similarly, the operator truncation penalty \( \GPKS^{2} = (\GPKS_{\mm}^{2}) \) 
corresponds to \( \GPKS_{\mm}^{2} = 0 \) for \( \mm \leq \mM \) and
\( \GPKS_{\mm}^{2} = \infty \) for \( \mm > \mM \).
Joint truncation of the target parameter \( \tarpv \) and of the operator \( \KS \) reduces 
the original \( \dimp \times \dimq \)-dimensional problem 
to the \( \jJ \times \mM \)-dimensional non-penalized problem with the operator 
\( \KSs_{\jJ,\mM} = \proj_{\mM} \KSs \Proj_{\jJ} \).
In the contrary to \cite{Ho2011} and \cite{Trabs_2018}, we accept \( \mM > \jJ \) and do not require \( \jJ = \mM \).

\Subsection{Estimation risk}
A \( (\jJ,\mM) \)-truncation penalty reduces the original problem to the non-penalized MLE 
\begin{EQA}[rcccl]
	\tilde{\prmtv}_{\jJ,\mM}
	&=&
	\argmax_{\prmtv = (\tarpv_{\jJ},\KS_{\jJ,\mM},\zv_{\mM})} \LL_{\jJ,\mM}(\prmtv)
	\, ,
	\qquad
	\prmtvs_{\jJ,\mM}
	&=&
	\argmax_{\prmtv = (\tarpv_{\jJ},\KS_{\jJ,\mM},\zv_{\mM})} \E \LL_{\jJ,\mM}(\prmtv)
	\, ,
\label{8sdi2ineuiujqdydtw3hychj}
\end{EQA}
with
\begin{EQA}[c]
	\LL_{\jJ,\mM}(\prmtv)
	\eqdef
	- \frac{1}{2} \bigl\{ \| \Zv_{\mM} - \zv_{\mM} \|^{2} + \muA^{2} \| \hat{\KS}_{\jJ,\mM} - \KS_{\jJ,\mM} \|_{\Fr}^{2}
	+ \| \zv_{\mM} - \KS_{\jJ,\mM} \tarpv_{\jJ} \|^{2} \bigr\}
	\, .
	\qquad
\label{ihfvic8kj3cdfcfi3b}
\end{EQA}
The results of Theorem~\ref{TsemieffEO} apply to this setup.
Moreover, many terms entering the formulation of this theorem can be specified in more detail.
First, we discuss the stochastic component. 
Denote 
\begin{EQA}[c]
	\DPN_{\jJ,\mM}^{2} \eqdef \KS_{\jJ,\mM}^{*\T} \KSs_{\jJ,\mM} 
	\, ,
	\quad
	\DFM_{\jJ,\mM}^{2} = \blk\{ \DPN_{\jJ,\mM}^{2},\Id_{\zv},\muA^{2} \Id_{\jJ,\mM} \} 
	\, ; 
\end{EQA}
cf. \eqref{iwdfhw787e3urne6fyh3o}.
Here \( \Id_{\jJ,\mM} \) stands for the identity operator on the space of \( \KS_{\jJ,\mM} \) matrices.
By \eqref{du7dyfu3ijhrk398rffi33} of Proposition~\ref{Leffscore}
\begin{EQA}[rcl]
	\Var \bigl( \IFT_{\GPT}^{-1} \nabla \zeta \bigr)_{\tarpv}
	& \leq &
	2 \dmax^{4} (\CONSTi_{\errZ} + 4 \CONSTi_{\errKSe} \eEO^{2}) 
	\DPN_{\jJ,\mM}^{-2} 
	\, .
\label{du7dyfu3ijhrktrunc}
\end{EQA}
Furthermore, restricting the source space \( \R^{\dimp} \) to \( \R^{\jJ} \) and the image space \( \R^{\dimq} \) to \( \R^{\mM} \) leads to 
reduction of the operator uncertainty measured by the value \( \dimA_{\KS} \) from Proposition~\ref{LdimfullEO}.
The \( \zv \)-component of the score vector reflecting the observation noise,
is not affected by the truncation of the source space.
However, the full dimensional operator noise \( \errKS \) is replaced with its restriction 
\( \errKS_{\jJ,\mM} \) to \( \R^{\jJ \times \mM} \) and \( \dimA_{\KS_{\jJ,\mM}} \leq \dmax^{4} \tr \Var(\errKS_{\jJ,\mM}) \).
This allows to reduce the critical dimension condition \eqref{jef9fe8e47fu32oeoggitir4g} to 
\( \muA^{-2} \, \jJ \, \mM \ll \nEO_{\jJ,\mM} = \lambda_{\min}(\DPN_{\jJ,\mM}^{2}) \).
The only negative effect of truncation is the bias \( \tarpvs_{\jJ,\mM} - \tarpvs \),
which can be easily analyzed.

\Subsection{Regularity of \( \KSs \) and smoothness of \( \tarpvs \)}
Regularity of \( \KSs \) will be described using the quantities
\begin{EQ}[rcccl]
	\nEO_{j}
	& \eqdef &
	\lambda_{j}(\Proj_{j} \DPN^{2} \Proj_{j})
	\, ,
	\qquad
	\nEOa_{\mm+1}
	& \eqdef &
	\| {\KSs}^{\T} (\Id_{\dimp} - \proj_{\mm}) \KSs  \|
	\, .
\label{vidf8vccil3l3wf6rbvy}
\end{EQ}
Here \( \lambda_{j}(\BBH) \) means the \( j \)th largest eigenvalue of the matrix \( \BBH \).
When considering \( \BBH = \Proj_{j} \DPN^{2} \Proj_{j} \) as a \( j \)-dimensional 
matrix, the value \( \nEO_{j} \) corresponds to its smallest eigenvalue. 
If the basis vectors \( \ev_{j} \) in \( \R^{\dimp} \) are defined as the ordered eigenvectors of \( \DPN^{2} \) then
\( \nEO_{j} = \nEOa_{j} \) coincide with the \( j \)th eigenvalue of \( \DPN^{2} \).
In general, \( \nEO_{j} \) and \( \nEOa_{j} \) might be different. 
For mildly/severely ill-posed problems, these values rapidly decrease with \( j \).
Given an operator \( \KS \), define for any \( j,\mm \)
\begin{EQA}[rclrcl]
	\KS_{j}
	& \eqdef &
	\KS \Proj_{j} \, ,
	\qquad
	&
	\DPN_{j}^{2}
	& \eqdef &
	\KS_{j}^{\T} \KS_{j} 
	\, ,
	\\
	\KS_{j,\mm}
	& \eqdef &
	\proj_{\mm} \KS \Proj_{j} \, ,
	\qquad
	&
	\DPN_{j}^{2}
	& \eqdef &
	\KS_{j,\mm}^{\T} \KS_{j,\mm} 
	\, .
\label{7d8cys8ike53efg83vcsh}
\end{EQA}

\begin{proposition}
\label{PappspaceEO}
Assume \eqref{vidf8vccil3l3wf6rbvy}.
Then for any \( \jJ \leq \dimp \) and \( \mM \leq \dimq \)
with \( \nEOa_{\mM+1} \leq \nEO_{\jJ}/2 \)
\begin{EQA}[rcl]
\label{fubynbunmn5323hgftwEOa}
	\| (\IFT_{\jJ,\mM}^{-1} \, \GPT_{\jJ,\mM}^{2} \prmtvs)_{\tarpv} \|
	& \leq &
	\| (\Id_{\dimp} - \Proj_{\jJ}) \tarpvs \| + \frac{\dmax^{2}\nEOa_{\mM+1}}{2\nEO_{\jJ}} \| \tarpvs \|
	\, ,
	\qquad
	\\
	\tr(\DPN_{\jJ}^{-2})
	& \leq &
	\sum_{j=1}^{\jJ} \frac{1}{\nEO_{j}} 
	\, , 
	\qquad
	\tr(\DPN_{\jJ,\mM}^{-2})
	\leq 
	\sum_{j=1}^{\jJ} \frac{1}{\nEO_{j} - \nEOa_{\mM+1}}
	\, . 
\end{EQA}
\end{proposition}

\begin{proof}
We split the proof into few steps.

\begin{lemma}
\label{LvarEOd}
For any \( \jJ \geq 1 \),
\begin{EQA}[c]
	\| ( \Proj_{\jJ} \IFTb_{\tarpv\tarpv} \Proj_{\jJ} )^{-} \|
	\leq 
	\frac{\dmax^{2}}{\nEO_{\jJ}}
	\, ,
	\qquad
	\tr(\DPN_{\jJ}^{-2})
	\leq 
	\sum_{j=1}^{\jJ} \frac{1}{\nEO_{j}} 
	\, .
	\qquad
\label{dfi786twe3t874r3ty}
\end{EQA}

\begin{proof}
As \( \DPN_{j}^{2} \) maps \( \R^{j} \) onto \( \R^{j} \), it follows from \( \| \DPN_{j}^{-2} \| \leq 1/\nEO_{j} \) that
\begin{EQA}[rcl]
	\tr(\DPN_{j}^{-2})
	& \leq &
	\tr(\DPN_{j-1}^{-2}) + 1/\nEO_{j}
	\, .
\end{EQA}
This and \eqref{vkr0roej83375e72hvirjrj} imply \eqref{dfi786twe3t874r3ty}
by induction.
\end{proof}

\end{lemma}
\begin{lemma}
\label{LDPNEiOa}
With \( \DPN_{j}^{2} \eqdef \Proj_{j} \DPN^{2} \Proj_{j} = 
\Proj_{j} \, \KS^{*\T} \KSs \Proj_{j} \),
\( \DPN_{j,\mm}^{2} = \Proj_{j} \KS^{*\T} \proj_{\mm} \KSs \Proj_{j} \), 
it holds 
\begin{EQA}[rcl]
	\| \DPN_{j,\mm}^{2} - \DPN_{j}^{2} \|
	& \leq &
	\nEOa_{\mm+1} 
	\, .
\end{EQA}
If \( \jJ \) and \( \mm \) are such that \( \nEOa_{\mm+1} \leq \nEO_{\jJ}/2 \) then
\begin{EQA}[c]
	\tr(\DPN_{\jJ,\mm}^{-2})
	\leq 
	\sum_{j=1}^{\jJ} \frac{1}{\nEO_{j} - \nEOa_{\mm+1}}
	\, .
\label{cyywegy3derd5erd53rt3}
\end{EQA}
\end{lemma}

\begin{proof}
It holds by definition
\begin{EQA}[rcl]
	\| \DPN_{j,\mm}^{2} - \DPN_{j}^{2} \|
	&=&
	\| \Proj_{j} \KS^{*\T} (\Id_{\dimq} - \proj_{\mm}) \KSs \Proj_{j} \|
	=
	\| (\Id_{\dimq} - \proj_{\mm}) \KSs \Proj_{j} \|^{2}
	\leq 
	\nEOa_{\mm+1}
	\, .
\end{EQA}
Hence, \( \lambda_{j}(\DPN_{j,\mm}^{2}) \geq \lambda_{j}(\DPN_{j}^{2}) - \nEOa_{\mm+1} 
= \nEO_{j} - \nEOa_{\mm+1} \)
and \eqref{cyywegy3derd5erd53rt3} follows as in Lemma~\ref{LDPNEiOa}.
\end{proof}

Further we evaluate the bias term due to  \( (\jJ,\mM) \)-truncation penalty.

\begin{lemma}
\label{PsemiriskEOt}
For the \( (\jJ,\mM) \)-truncation penalty operator \( \GPT_{\jJ,\mM}^{2} \),
it holds 
\begin{EQA}[c]
	\| (\IFT_{\jJ,\mM}^{-1} \, \GPT_{\jJ,\mM}^{2} \prmtvs)_{\tarpv} \|
	\leq 
	\| (\Id_{\dimp} - \Proj_{\jJ}) \tarpvs \| + \frac{\dmax^{2}\nEOa_{\mM+1}}{2\nEO_{\jJ}} \| \tarpvs \|
	\, .
	\qquad
\label{fubynbunmn5323hgftw}
\end{EQA}
\end{lemma}

\begin{proof}[Proof of Lemma~\ref{PsemiriskEOt}]
Consider the \( (\jJ,\mM) \)-truncation penalties \( \GPT_{\jJ,\mM}^{2} \).
For bias evaluation, we apply the general bound of Lemma~\ref{PsemiriskEO}.
The use of \( \GP^{2} = \gp^{2} (\Id_{\dimp} - \Proj_{\jJ}) \) with \( \gp^{2} = \infty \) yields
with \( \IFTb_{\tarpv\tarpv} = \IFT_{\tarpv\tarpv} - \IFT_{\tarpv\nupv} \, \IFT_{\GPT,\nupv\nupv}^{-1} \IFT_{\nupv\tarpv} \)
\begin{EQA}[rcl]
	\IFTb_{\GPT,\tarpv\tarpv}^{-1} \GP^{2}
	&=&
	\bigl( \IFTb_{\tarpv\tarpv} + \GP^{2} \bigr)^{-1} \GP^{2}
	=
	\Id_{\dimp} - \Proj_{\jJ}
	\, .
\end{EQA}
So, truncation in the source space yields the bias \( \IFTb_{\GPT,\tarpv\tarpv}^{-1} \GP^{2} \tarpvs = (\Id_{\dimp} - \Proj_{\jJ}) \tarpvs \).
Similarly, 
\begin{EQA}[rcl]
	\IFTb_{\GPT,\tarpv\tarpv}^{-1} 
	&=&
	\bigl( \IFTb_{\tarpv\tarpv} + \GP^{2} \bigr)^{-1} 
	=
	\bigl( \Proj_{\jJ} \IFTb_{\tarpv\tarpv} \Proj_{\jJ} \bigr)^{-}
	\, ,
\end{EQA}
and by Lemma~\ref{LIFTDFMEO}, 
\( \bigl( \Proj_{\jJ} \IFTb_{\tarpv\tarpv} \Proj_{\jJ} \bigr)^{-} \leq \dmax^{2} \bigl( \Proj_{\jJ} \DPN^{2} \Proj_{\jJ} \bigr)^{-} \).
Further,
\eqref{kchdhv63bfub83kfgibnurr} with \( \gpks_{\mm}^{2} = 0 \) for \( \mm \leq \mM \) and \( \gpks_{\mm}^{2} = \infty \) for \( \mm > \mM \)
yields
\begin{EQA}[rcl]
	S_{\mM}
	& = &
	\frac{1}{2} 
	\sum_{\mm=\mM+1}^{\dimq} \KSs_{\mm} {\KSs_{\mm}}^{\T}
	\, .
\end{EQA}
Therefore, truncation of the operator \( \KS \) results in the bias
\begin{EQA}[c]
	\IFTb_{\GPT,\tarpv\tarpv}^{-1} \, S_{\mM} \, \tarpvs
	=
	\frac{1}{2} 
	\bigl( \Proj_{\jJ} \IFTb_{\tarpv\tarpv} \Proj_{\jJ} \bigr)^{-} \sum_{\mm=\mM+1}^{\dimq} \KSs_{\mm} {\KSs_{\mm}}^{\T} \tarpvs
\end{EQA}
The norm of \( S_{\mM} \) can easily be bounded under \eqref{vidf8vccil3l3wf6rbvy}: for any unit vector \( \uv \in \R^{\dimp} \)
\begin{EQA}[c]
	\uv^{\T} S_{\mM} \uv
	=
	\frac{1}{2} 
	\sum_{\mm=\mM+1}^{\dimq} |{\KSs_{\mm}}^{\T} \uv|^{2}
	=
	\frac{1}{2} 
	\bigl\| (\Id_{\dimq} - \proj_{\mM}) \KSs \uv \bigr\|^{2}
	\leq 
	\frac{1}{2} \| \KS^{*\T} (\Id_{\dimq} - \proj_{\mM}) \KSs \|
	= 
	\frac{\nEOa_{\mM+1}}{2} 
	\, .
\end{EQA}
This and \eqref{dfi786twe3t874r3ty} imply
\begin{EQA}[rcl]
	\| ( \Proj_{\jJ} \IFTb_{\tarpv\tarpv} \Proj_{\jJ} )^{-} S_{\mM} \tarpvs \|
	& \leq &
	\| ( \Proj_{\jJ} \IFTb_{\tarpv\tarpv} \Proj_{\jJ} )^{-} \| \; \| S_{\mM} \| \; \| \tarpvs \|
	\\
	& \leq &
	\dmax^{2} \| ( \Proj_{\jJ} \DPN^{2} \Proj_{\jJ} )^{-} \| \; \| S_{\mM} \| \; \| \tarpvs \|
	\leq 
	\frac{\dmax^{2}\nEOa_{\mM+1}}{2\nEO_{\jJ}} \| \tarpvs \|
	\, ,
\end{EQA}
and the result follows.
\end{proof}

Putting together the obtained bounds completes the proof of the proposition.
\end{proof}

\Section{Rate of estimation in inverse problems}

This section illustrates the obtained results by showing that 
a ridge penalty for a smooth operator or a truncation penalty with properly
selected parameters \( \jJ,\mM \) lead to 
rate optimal accuracy of estimation.
As in the previous sections, we fix an approximation spaces setup.
Regularity of the operator \( \KSs \) is described by decreasing sequences 
\( \nEO_{1} \geq \nEO_{2} \geq \ldots \geq \nEO_{\dimp} \) and
\( \nEOa_{1} \geq \nEOa_{2} \geq \ldots \geq \nEOa_{\dimq} \) ensuring \eqref{vidf8vccil3l3wf6rbvy}.
Moreover, we only discuss the case of a mildly ill-posed problem when the values \( \nEO_{j} \) decrease
polynomially; see \eqref{udvcghue3j3vc654rw2y}.
We also assume a Sobolev smoothness of the signal \( \tarpv \) as in \eqref{dcyhvc55t8tg764e33rf}
with \( \cgp_{j} \) polynomially increasing.
A popular special case is given by \( \nEO_{j} \approx \nEO_{1} \, j^{-2s} \) 
and \( \cgp_{j}^{2} \approx \CGPz \, j^{2\beta} \).

\begin{proposition}
\label{PrateEiO}
Let the operator \( \KSs \) fulfill \eqref{vidf8vccil3l3wf6rbvy} with 
\( \nEO_{j} \geq \nEO_{1} \, j^{-2s} \) for all \( j > 1 \).
Further, let \( \tarpvs \) follow \eqref{dcyhvc55t8tg764e33rf} with \( \cgp_{j}^{2} = \CGPz \, j^{2\beta} \).
Define \( \jJ \) by \( \jJ \asymp (\nEO_{1}/\CGPz)^{1/(1 + 2\beta + 2s)} \) and \( \mM \) by 
\( \nEOa_{\mM+1} \leq \rho \nEO_{\jJ} \) for \( \rho \leq 1/2 \).
Then the risk \( \riskt_{\jJ,\mM} \) of the estimator \( \tilde{\tarpv}_{\jJ,\mM} \)
satisfies
\begin{EQA}[c]
	\riskt_{\jJ,\mM}
	\lesssim 
	\CGPz^{-\frac{2s+1}{1 + 2\beta + 2s}} \nEO_{1}^{-\frac{2\beta}{1 + 2\beta + 2s}}
	\, .
\end{EQA}
\end{proposition}
\begin{proof}
Define \( \jJ \) by \( \jJ \asymp (\nEO_{1}/\CGPz)^{1/(1 + 2\beta + 2s)} \).
\eqref{vidf8vccil3l3wf6rbvy} implies
\begin{EQA}[c]
	\tr(\DPN_{\jJ}^{-2})
	\leq 
	\sum_{j=1}^{\jJ} \frac{1}{\nEO_{j}} 
	\leq 
	\frac{1}{\nEO_{1}} \sum_{j=1}^{\jJ} j^{2s}
	\leq 
	\frac{1}{2s+1} \frac{\jJ^{2s+1}}{\nEO_{1}}
	\, .
\end{EQA}
If \( \mM \) is fixed to ensure \( \nEOa_{\mM+1} \leq \rho \nEO_{\jJ} \) for \( \rho \leq 1/2 \), 
then a similar bound applies to \( \tr(\DPN_{\jJ,\mM}^{-2}) \); see \eqref{cyywegy3derd5erd53rt3}.
Further, as \( \cgp_{j}^{2} = \CGPz \, j^{2\beta} \) increase with \( j \), it holds
\begin{EQA}[c]
	\| (\Id_{\dimp} - \Proj_{\jJ}) \tarpvs \|^{2}
	=
	\sum_{j=\jJ+1}^{\dimp} \langle \tarpvs, \ev_{j} \rangle^{2}
	\leq 
	\cgp_{\jJ}^{-2}
	\sum_{j=\jJ+1}^{\dimp} \cgp_{j}^{2} \langle \tarpvs, \ev_{j} \rangle^{2}
	\leq 
	\cgp_{\jJ}^{-2}
	= 
	\CGPz^{-1} \, \jJ^{-2\beta}
	\, .
\end{EQA}
This yields
\begin{EQA}[c]
	\riskt
	\lesssim  
	\tr(\DPN_{\jJ}^{-2}) + \cgp_{\jJ}^{-2}
	\leq 
	\nEO_{1}^{-1} \jJ^{2s+1} + \CGPz^{-1} \, \jJ^{-2\beta}
	\lesssim
	\CGPz^{-\frac{2s+1}{1 + 2\beta + 2s}} \nEO_{1}^{-\frac{2\beta}{1 + 2\beta + 2s}}
\end{EQA}
as stated.
\end{proof}

\iffourG{
\Subsection{4S bounds}
This section presents advanced bounds for the Fisher expansion and risk of estimation
based on the study from Section~\ref{Sfourthsemi}.
Consider the third-order tensor
\( \TensU(\uv) = \frac{1}{6} \langle \nabla^{3} \fs(\prmtvs_{\GP}), \uv^{\otimes 3} \rangle \)
and its gradient \( \frac{1}{2} \langle \nabla^{3} \fs(\prmtvs_{\GP}), \uv^{\otimes 2} \rangle \).
Here \( \uv = (\dtarpv,\dimav,\dKS) \) is a vector in the full parameter space \( \R^{\dimfull} \).
Definitions \eqref{dhfdfededrsdsrdsresswddd} and \eqref{iduejef634yhdjdekee} imply 
for \( \prmtv = (\tarpv,\imav,\KS) \)
\begin{EQA}
	\TensU(\prmtv,\uv) 
	& \eqdef & 
	\frac{1}{6} \langle \nabla^{3} \fs(\prmtv), \uv^{\otimes 3} \rangle
	=
	- \frac{1}{12} \frac{d^{3}}{dt^{3}} 
	\bigl\| \imav + t \dimav - (\KS + t \dKS) (\tarpv + t \dtarpv) \bigr\|^{2} 
	\bigg|_{t=0}
	\\
	&=&
	\frac{1}{2} (\dKS \tarpv + \KS \dtarpv - \dimav) \dKS \dtarpv
	\, .
\label{tuxuaiiiid5rww3wdfnh}
\end{EQA}
With \( \IFT_{\GPT} = \IFT_{\GPT}(\prmtvs_{\GPT}) \), define 
\begin{EQA}[rcccl]
	\svn_{\GPT} 
	&=&
	\IFT_{\GPT}^{-1} \bigl\{ \nabla \zeta + \nabla \TensU(\IFT_{\GPT}^{-1} \nabla \zeta) \bigr\} 
	&=&
	(\svn_{\GPT,\tarpv},\svn_{\GPT,\nupv}) \, ,
\label{8vfjvr43223efryfuweefgEO}
	\\
	\bvn_{\GPT}
	&=&
	\IFT_{\GPT}^{-1} \bigl\{ \GPT^{2} \prmtvs + \nabla \TensU(\IFT_{\GPT}^{-1} \GPT^{2} \prmtvs) \bigr\} 
	&=&
	(\bvn_{\GPT,\tarpv},\bvn_{\GPT,\nupv}) \, .
\label{8vfjvr43f8khg54ed5GEO}
\end{EQA}
We apply Theorem~\ref{Tsemieffp4}.

\begin{theorem}
\label{TseffEO}
Assume the conditions of Theorem~\ref{TsemieffEO} and let
\( \dmax^{2} \dltwu_{4} \, (\rr_{\DFM} \vee \bias_{\DFM})^{2} < \frac{1}{3} \).
For any linear mapping \( \QP \) of \( \tarpv \), it holds on \( \Omega(\xx) \)
\begin{EQ}[rcl]
	\| \QP \, ( \tilde{\tarpv}_{\GPT} - \tarpvs_{\GPT} - \svn_{\GPT,\tarpv} ) \|
	& \leq &
	\| \QP \, \DPN^{-1} \| \, \dmax^{2} 
	\Bigl( \frac{\dltwu_{4}}{2} + \dmax^{2} \dltwu_{3}^{2} \Bigr) \, \| \DFM \IFT_{\GPT}^{-1} \nabla \zeta \|^{3} \, ,
	\qquad
	\\
	\| \QP \, ( \tarpvs_{\GPT} - \tarpvs + \bvn_{\GPT,\tarpv}) \|
	& \leq &
	\| \QP \, \DPN^{-1} \| \, \dmax^{2} 
	\Bigl( \frac{\dltwu_{4}}{2} + \dmax^{2} \dltwu_{3}^{2} \Bigr) \, \bias_{\DFM}^{3} \, ,
	\qquad
\label{0mkvhgjnrwwe3u8gtygEO}
\end{EQ}
and
\begin{EQ}[rcl]
	\| \QP (\svn_{\GPT} - \IFT_{\GPT}^{-1} \, \nabla \zeta)_{\tarpv} \|
	& \leq &
	\| \QP \, \DPN^{-1} \| \, \frac{\dmax^{2} \dltwu_{3}}{2} \, 
	\| \DFM \IFT_{\GPT}^{-1} \nabla \zeta \|^{2} \, ,
	\qquad
	\\
	\| \QP (\bvn_{\GPT} - \IFT_{\GPT}^{-1} \, \GP^{2} \prmtvs)_{\tarpv} \|
	& \leq &
	\| \QP \, \DPN^{-1} \| \, \frac{\dmax^{2} \dltwu_{3}}{2} \, 
	\bias_{\DFM}^{2} \, .
	\qquad
\label{vuedy766t4e3bfvyt6eEO}
\end{EQ}
\end{theorem}

Now we study the risk of \( \tilde{\tarpv}_{\GP} \) using Theorem~\ref{Tseff42}.
Define 
\begin{EQA}
\label{hdvje39bug53ebfh8edxEO}
	\riskt_{\QP,2} 
	& \eqdef & 
	\E \bigl\{ \| \QP (\svn_{\GPT,\tarpv} - \bvn_{\GPT,\tarpv}) \|^{2} \Ind_{\Omega(\xx)} \bigr\} \, .
\label{yfhcvhched6chejdrteEO}
\end{EQA}

\begin{theorem}
\label{Tseff42EO}
Assume the conditions of Theorem~\ref{TseffEO} and let
\begin{EQA}[c]
	\E \bigl\{ \| \DFM \IFT_{\GP}^{-1} \nabla \zeta \|^{k} \Ind_{\Omega(\xx)} \bigr\} \leq \CONSTi_{k}^{2} \, \dimAfull_{\DFM}^{k/2} ,
	\qquad
	k=3,4,6 \, .
\label{6hjdfv8e6hyefyeheew7skGPT}
\end{EQA}
For a linear mapping \( \QP \) and \( \riskt_{\QP,2} \) from \eqref{hdvje39bug53ebfh8edxEO}, it holds
\begin{EQA}
	&& \nquad
	\E \bigl\{ \| \QP \, (\tilde{\tarpv}_{\GPT} - \tarpvs) \| \Ind_{\Omega(\xx)} \bigr\}
	\leq 	 
	\E \bigl\{ \| \QP (\svn_{\GPT,\tarpv} - \bvn_{\GPT,\tarpv}) \| \Ind_{\Omega(\xx)} \bigr\}
	\\
	&&
	+ \, \| \QP \DPN^{-1} \| \, \dmax^{2} 
	\Bigl( \frac{\dltwu_{4}}{2} + \dmax^{2} \dltwu_{3}^{2} \Bigr) \, \bigl( \CONSTi_{3}^{2} \, \dimAfull_{\DFM}^{3/2} + \bias_{\DFM}^{3} \bigr)
\label{0mkvhgjnrw3dfwe3u8gtygE1EO}
\end{EQA}
and 
\begin{EQA}
	&& \nquad
	\Bigl| \E \bigl\{ \| \QP (\svn_{\GPT,\tarpv} - \bvn_{\GPT,\tarpv}) \| \Ind_{\Omega(\xx)} \bigr\}
	- \E \bigl\{ \| \QP \, (\IFT_{\GPT}^{-1} \, \nabla \zeta - \IFT_{\GPT}^{-1} \, \GP^{2} \prmtvs)_{\tarpv} \| \Ind_{\Omega(\xx)} \bigr\} 
	\Bigr|
	\\
	& \leq &
	\| \QP \, \DPN^{-1} \| \, \frac{\dmax^{2} \dltwu_{3}}{2} \, \bigl( \dimAfull_{\DFM} + \bias_{\DFM}^{2} \bigr) \, .
\label{udtgecthwjdytdehduuc6EO}
\end{EQA}
Furthermore,
\begin{EQA}
	\bigl( 1 - \alp_{\QP,2} \bigr)^{2} \riskt_{\QP,2}
	\leq 
	\E \bigl\{ \| \QP \, (\tilde{\tarpv}_{\GPT} - \tarpvs) \|^{2} \Ind_{\Omega(\xx)} \bigr\}
	& \leq &
	\bigl( 1 + \alp_{\QP,2} \bigr)^{2} \riskt_{\QP,2} \, 
\label{6shx76whnjvyehfbvyfhEO}
\end{EQA}
provided that
\begin{EQA}
	\alp_{\QP,2}
	& \eqdef &
	\frac{\| \QP \, \DPN^{-1} \| \, \dmax^{2} ( \dltwu_{4}/2 + \dmax^{2} \dltwu_{3}^{2} ) \, ( \CONSTi_{6} \, \dimAfull_{\DFM}^{3/2} + \bias_{\DFM}^{3})}
	 	 {\sqrt{\riskt_{\QP,2}}} 
	< 1 \, .
\label{cuyhe4f8jfdebnuiEO}
\end{EQA}
If another value \( \alp_{\QP,1} < 1 \) is such that  
\begin{EQ}[rcl]
	&&
	\| \QP \, \DPN^{-1} \| \, \frac{\dmax^{2} \dltwu_{3}}{2} \, \bigl( \CONSTi_{4} \, \dimAfull_{\DFM} + \bias_{\DFM}^{2} \bigr)
	\leq 
	\alp_{\QP,1} \, \sqrt{\riskt_{\QP}} \, 
\label{6dhx6whcuydsds655srew4EO}
\end{EQ}
with \( \riskt_{\QP} \) from \eqref{t6skjid7ehfcr56tegfdEO} then
\begin{EQA}
	\riskt_{\QP} (1 - \alp_{\QP,1})^{2} 
	\leq 
	\riskt_{\QP,2}
	& \leq &
	\riskt_{\QP} (1 + \alp_{\QP,1})^{2} \, .
\label{EQtuGmstrVEQtGQEO}
\end{EQA}
\end{theorem}


}{}

\bibliography{exp_ts,listpubm-with-url}

\newpage
\appendix
\section{General conditions}
This section collects the conditions required for the general statements from \cite{Sp2024}.
\Section{Basic conditions}
\label{Scondgeneric}
Now we present our major conditions.
The most important one is about linearity of the stochastic component 
\( \zeta(\upsv) = L(\upsv) - \E L(\upsv) = L(\upsv) - \E L(\upsv) \).

\medskip
\begin{description}
    \item[\label{Eref} \( \bb{(\zeta)} \)]
      \textit{The stochastic component \( \zeta(\upsv) = L(\upsv) - \E L(\upsv) \) of the process \( L(\upsv) \) is linear in 
      \( \upsv \). 
      We denote by \( \nabla \zeta \equiv \nabla \zeta(\upsv) \in \R^{\dimp} \) its gradient
      }.
\end{description}

Below we assume some concentration properties of the stochastic vector \( \nabla \zeta \).
More precisely, we require that \( \nabla \zeta \) obeys the following condition.

\begin{description}
\item[\label{EU2ref}\( \bb{(\nabla \zeta)} \)]
	\textit{There exists \( \VP^{2} \geq \Var(\nabla \zeta) \) such that  
	for all considered  \( \BBH \in \Matr_{\dimp} \) and \( \xx > 0 \)
	}
\begin{EQA}
	\P\bigl( \| \BBH^{1/2} \VP^{-1} \nabla \zeta \| \geq \zq(\BBH,\xx) \bigr)
	& \leq &
	3 \ex^{-\xx} ,
\label{2emxGPm12nz122}
	\\
	\zq^{2}(\BBH,\xx)
	& \eqdef &
	\tr \BBH + 2 \sqrt{\xx \, \tr \BBH^{2}} + 2 \xx \| \BBH \| \,  .
\label{34rtyghuioiuyhgvftid}
\end{EQA}
\end{description}

This condition can be effectively checked if the errors in the data exhibit sub-gaussian or sub-exponential behavior; see 
\ifsqnorm{Section~\ref{SdevboundnonGauss}.}{\cite{Sp2023c}, \cite{Sp2023d}.}
The important special case corresponds to \( \BBH = \IF^{-1/2} \VP^{2} \IF^{-1/2} \) 
and \( \xx \approx \log n \) leading to the bound
\begin{EQA}
	\P\bigl( \| \IF^{-1/2} \nabla \zeta \| > \zq(\BBH,\xx) \bigr)
	& \leq &
	3/n .
\label{udyvfeyejff6777dj23}
\end{EQA}
The  value \( \dimG = \tr(\IF^{-1} \VP^{2}) \) can be called the \emph{effective dimension}; see \cite{SP2013_rough}.

We also assume that the log-likelihood \( L(\upsv) \) or, equivalently, its deterministic part 
\( \E L(\upsv) \) is a concave function.

\medskip
\begin{description}
    \item[\label{LLref} \( \bb{(\mathcal{C})} \)]
      \textit{The function \( \fs(\upsv) \eqdef \E L(\upsv) \) is concave on \( \Ups \) which is open and convex set in \( \R^{\dimp} \).      
      }
\end{description}

\Section{Gateaux smoothness and self-concordance}
Below we assume 
the function \( \fs(\upsv) \) to be strongly concave with the negative Hessian 
\( \IFN(\upsv) \eqdef - \nabla^{2} \fs(\upsv) \in \Matr_{\dimp} \) positive definite. 
Also, assume \( \fs(\upsv) \) three or sometimes even four times Gateaux differentiable in \( \upsv \in \Ups \).
For any particular direction \( \uv \in \R^{\dimp} \), we consider the univariate function 
\( \fs(\upsv + t \uv) \) and measure its smoothness in \( t \).
Local smoothness of \( \fs \) will be described by the relative error of the Taylor expansion 
of the third or fourth order.
Namely, define
\begin{EQ}[rcl]
	\dltw_{3}(\upsv,\uv) 
	&=& 
	\fs(\upsv + \uv) - \fs(\upsv) - \langle \nabla \fs(\upsv), \uv \rangle 
	- \frac{1}{2} \langle \nabla^{2} \fs(\upsv), \uv^{\otimes 2} \rangle , 
	\\
	\dltwd_{3}(\upsv,\uv) 
	&=&
	\langle \nabla \fs(\upsv + \uv), \uv \rangle - \langle \nabla \fs(\upsv), \uv \rangle 
	- \langle \nabla^{2} \fs(\upsv), \uv^{\otimes 2} \rangle \, ,
\label{dltw3vufuv12f2ga}
\end{EQ}
and
\begin{EQA}
	\dltw_{4}(\upsv,\uv)
	& \eqdef &
	\fs(\upsv + \uv) - \fs(\upsv) - \langle \nabla \fs(\upsv), \uv \rangle 
	- \frac{1}{2} \langle \nabla^{2} \fs(\upsv), \uv^{\otimes 2} \rangle
	- \frac{1}{6} \langle \nabla^{3} \fs(\upsv), \uv^{\otimes 3} \rangle \, .
\label{hvcduywgedfuyg2y1y35e3wweg}
\end{EQA}
Now, for each \( \upsv \), suppose to be given a positive symmetric operator 
\( \DFN(\upsv) \in \Matr_{\dimp} \) 
defining a local metric and a local vicinity around \( \upsv \):
\begin{EQA}
	\UVz_{\rr}(\upsv)
	&=&
	\bigl\{ \uv \in \R^{\dimp} \colon \| \DFN(\upsv) \uv \| \leq \rr \bigr\}
\label{ed7sycf7wedwgedq2ftwdfgtv}
\end{EQA}
for some radius \( \rr \).

Local smoothness properties of \( \fs \) at \( \upsv \) are given via the quantities
\begin{EQA}[rcccl]
    \dltwb(\upsv)
    & \eqdef &
    \sup_{\uv \colon \| \DFN(\upsv) \uv \| \leq \rr} \,
    \frac{2|\dltw_{3}(\upsv,\uv)|}{\| \DFN(\upsv) \uv \|^{2}} 
    \,\, ,
    \quad
    \dltwbd(\upsv)
    & \eqdef &
    \sup_{\uv \colon \| \DFN(\upsv) \uv \| \leq \rr} \, \frac{|\dltwd_{3}(\upsv,\uv)|}{\| \DFN(\upsv) \uv \|^{2}} \,\, . 
    \qquad
\label{dtb3u1DG2d3GPg}
\end{EQA}
The definition yields for any \( \uv \) with \( \| \DFN(\upsv) \uv \| \leq \rr \)
\begin{EQ}[rcccl]
	\bigl| \dltw_{3}(\upsv,\uv) \rangle \bigr|
	& \leq &
	\frac{\dltwb(\upsv)}{2} \| \DFN(\upsv) \uv \|^{2} 
	\, ,
	\qquad
	\bigl| \dltwd_{3}(\upsv,\uv) \bigr|
	& \leq &
	\dltwbd(\upsv) \| \DFN(\upsv) \uv \|^{2}
	\, .
	\qquad
\label{dta3u1DG2d3GPa1g}
\end{EQ}
%
The approximation results can be improved 
provided a third order upper bound on the error of Taylor expansion. 

\begin{description}
    \item[\label{LL3tref} \( \bb{(\mathcal{T}_{3})} \)]
      \textit{For some \( \dltwu_{3} \)}
\begin{EQA}
	\bigl| \dltw_{3}(\upsv,\uv) \bigr|
	& \leq &
	\frac{\dltwu_{3}}{6} \| \DFN(\upsv) \, \uv \|^{3} \, ,
	\quad
	\bigl| \dltwd_{3}(\upsv,\uv) \bigr|
	\leq 
	\frac{\dltwu_{3}}{2} \| \DFN(\upsv) \, \uv \|^{3} \, ,
	\quad
	\uv \in \UVz_{\rr}(\upsv).
\label{bd3xu16f3uo3st}
\end{EQA}
\end{description}
 
\begin{description}
    \item[\label{LL4tref} \( \bb{(\mathcal{T}_{4})} \)]
      \textit{For some \( \dltwu_{4} \)}
\begin{EQA}
	\bigl| \dltw_{4}(\upsv,\uv) \bigr|
	& \leq &
	\frac{\dltwu_{4}}{24} \| \DFN(\upsv) \, \uv \|^{4} \, ,
	\qquad
	\uv \in \UVz_{\rr}(\upsv).
\label{1mffmxum5st}
\end{EQA}
\end{description}

We also present a version of \nameref{LL3tref} resp. \nameref{LL4tref} in terms of the third (resp. fourth) derivative of \( \fs \).
\begin{description}
    \item[\label{LLsT3ref} \( \bb{(\mathcal{T}_{3}^{*})} \)]
    \emph{\( \fs(\upsv) \) is three times differentiable and 
	}
\begin{EQA}
    \sup_{\uv \colon \| \DFN(\upsv) \uv \| \leq \rr} \,\, \sup_{\zv \in \R^{\dimp}} \,\, 
    \frac{\bigl| \langle \nabla^{3} \fs(\upsv + \uv), \zv^{\otimes 3} \rangle \bigr|}
		 {\| \DFN(\upsv) \zv \|^{3}} 
	& \leq &
	\dltwu_{3} \, .
\label{jcxhydtferyu9j3d6vhew}
\end{EQA}

    \item[\label{LLsT4ref} \( \bb{(\mathcal{T}_{4}^{*})} \)]
    \emph{\( \fs(\upsv) \) is four times differentiable and 
	}
\begin{EQA}
    \sup_{\uv \colon \| \DFN(\upsv) \uv \| \leq \rr} \,\, \sup_{\zv \in \R^{\dimp}} \,\, 
    \frac{\bigl| \langle \nabla^{4} \fs(\upsv + \uv), \zv^{\otimes 4} \rangle \bigr|}
		 {\| \DFN(\upsv) \zv \|^{4}} 
	& \leq &
	\dltwu_{4} \, .
\label{jcxhydtferyu9j3d6vhew4}
\end{EQA}

\end{description}

%
\noindent
By Banach's characterization \cite{Banach1938}, \nameref{LLsT3ref} implies
\begin{EQA}
	\bigl| \langle \nabla^{3} \fs(\upsv + \uv), \zv_{1} \otimes \zv_{2} \otimes \zv_{3} \rangle \bigr|
	& \leq &	 
	\dltwu_{3} \| \DFN(\upsv) \zv_{1} \| \, \| \DFN(\upsv) \zv_{2} \| \, \| \DFN(\upsv) \zv_{3} \| \, 
\label{jbuyfg773jgion94euyyfg}
\end{EQA}
for any \( \uv \) with \( \| \DFN(\upsv) \uv \| \leq \rr \) and all \( \zv_{1} , \zv_{2}, \zv_{3} \in \R^{\dimp} \).
Similarly under \nameref{LLsT4ref}
\begin{EQA}
	\bigl| \langle \nabla^{4} \fs(\upsv + \uv), \zv_{1} \otimes \zv_{2} \otimes \zv_{3} \otimes \zv_{4} \rangle \bigr|
	& \leq &	 
	\dltwu_{4} \prod_{k=1}^{4} \| \DFN(\upsv) \zv_{k} \| \, ,
	\quad 
	\zv_{1} , \zv_{2}, \zv_{3}, \zv_{4} \in \R^{\dimp} \, .
	\qquad
\label{jbuyfg773jgion94euyyfg4}
\end{EQA}

\begin{lemma}
\label{LdltwLa3t}
Under \nameref{LL3tref} or \nameref{LLsT3ref},
the values \( \dltwb(\upsv) \) and \( \dltwbd(\upsv) \) from \eqref{dtb3u1DG2d3GPg} satisfy
\begin{EQA}[rcccl]
\label{gtcdsftdffrvsewsea}
	\dltwb(\upsv)
	& \leq &
	\frac{\dltwu_{3} \, \rr}{3 } \, ,
	\qquad
	\dltwbd(\upsv)
	& \leq &
	\frac{\dltwu_{3} \, \rr}{2} \, ,
	\qquad
	\upsv \in \Upsd .
\label{gtcdsftdfvtwdsefhfdvfrvsewseG}
\end{EQA}
\end{lemma}

\begin{proof}
For any \( \uv \in \UVz_{\rr}(\upsv) \) with \( \| \DFN(\upsv) \uv \| \leq \rr \)
\begin{EQA}
	\bigl| \dltw_{3}(\upsv,\uv) \bigr|
	& \leq &
	\frac{\dltwu_{3}}{6} \, \| \DFN(\upsv) \uv \|^{3} 
	\leq 
	\frac{\dltwu_{3} \, \rr}{6} \, \| \DFN(\upsv) \uv \|^{2},
\label{jrgeteteer2234587654}
\end{EQA}
and the bound for \( \dltwb(\upsv) \) follows.
The proof for \( \dltwbd(\upsv) \) is similar.
\end{proof}

The values \( \dltwu_{3} \) and \( \dltwu_{4} \) are usually very small.
Some quantitative bounds are given later in this section
under the assumption that the function \( \fs(\upsv) \) can be written in the form \( - \fs(\upsv) = n \hL(\upsv) \) 
for a fixed smooth function \( h(\upsv) \) with the Hessian \( \nabla^{2} \hL(\upsv) \). 
The factor \( n \) has meaning of the sample size%
\ifapp{; see \Chname \ref{ScritdimMLE} or \Chname \ref{SGBvM}.}{.}

\begin{description}
    \item[\label{LLtS3ref} \( \bb{(\mathcal{S}_{3}^{*})} \)]
      \emph{ \( - \fs(\upsv) = n \hL(\upsv) \) for \( \hL(\upsv) \) three times differentiable and
\begin{EQA}
	\sup_{\uv \colon \| \HL(\upsv) \uv \| \leq \rr/\sqrt{n}} 
	\frac{\bigl| \langle \nabla^{3} \hL(\upsv + \uv), \uv^{\otimes 3} \rangle \bigr|}{\| \HL(\upsv) \uv \|^{3}}
	& \leq &
	\hmax_{3} \, .
\end{EQA}
}
    \item[\label{LLtS4ref} \( \bb{(\mathcal{S}_{4}^{*})} \)]
      \emph{ the function \( \hL(\cdot) \) satisfies \nameref{LLtS3ref} and  
\begin{EQA}
	\sup_{\uv \colon \| \HL(\upsv) \uv \| \leq \rr/\sqrt{n}}
	\frac{\bigl| \langle \nabla^{4} \hL(\upsv + \uv), \uv^{\otimes 4} \rangle \bigr|}{\| \HL(\upsv) \uv \|^{4}}
	& \leq &
	\hmax_{4} \, .
\end{EQA}
}
\end{description}

\noindent
\nameref{LLtS3ref} and \nameref{LLtS4ref}
are local versions of the so-called self-concordance condition; see \cite{Ne1988} and \cite{OsBa2021}.
In fact, they require that each univariate function \( \hL(\upsv + t \uv) \) of \( t \in \R \)
is self-concordant with some universal constants \( \hmax_{3} \) and \( \hmax_{4} \).
Under \nameref{LLtS3ref} and \nameref{LLtS4ref}, with \( \DFN^{2}(\upsv) = n \, \HL^{2}(\upsv) \), the values 
\( \dltw_{3}(\upsv,\uv) \), \( \dltw_{4}(\upsv,\uv) \), and \( \dltwb(\upsv) \), \( \dltwbd(\upsv) \) can be
bounded.

\begin{lemma}
\label{LdltwLaGP}
Suppose \nameref{LLtS3ref}.
Then 
\nameref{LL3tref} follows with \( \dltwu_{3} = \hmax_{3} n^{-1/2} \).
Moreover, for \( \dltwb(\upsv) \) and \( \dltwbd(\upsv) \) from \eqref{dtb3u1DG2d3GPg}, it holds
\begin{EQA}[rcccl]
	\dltwb(\upsv)
	& \leq &
	\frac{\hmax_{3} \, \rr}{3 n^{1/2}} \, ,
	\qquad
	\dltwbd(\upsv)
	& \leq &
	\frac{\hmax_{3} \, \rr}{2 n^{1/2}} \, .
\label{gtcdsftdfvtwdsefhfdvfrvsewseGP}
\end{EQA}
Also \nameref{LL4tref} follows from \nameref{LLtS4ref} with \( \dltwu_{4} = \hmax_{4} n^{-1} \).
\end{lemma}

\begin{proof}
For any \( \uv \in \UVz_{\rr}(\upsv) \) and \( t \in [0,1] \), by the Taylor expansion of the third order
\begin{EQA}
	|\dltw(\upsv,\uv)|
	& \leq &
	\frac{1}{6} \bigl| \langle \nabla^{3} \fs(\upsv + t \uv), \uv^{\otimes 3} \rangle \bigr|
	=
	\frac{n}{6} \, \bigl| \langle \nabla^{3} \hL(\upsv + t \uv), \uv^{\otimes 3} \rangle \bigr|
	\leq 
	\frac{n \, \hmax_{3}}{6} \, \| \HL(\upsv) \uv \|^{3} 
	\\
	&=&
	\frac{n^{-1/2} \, \hmax_{3}}{6} \, \| \DFN(\upsv) \uv \|^{3}
	\leq 
	\frac{n^{-1/2} \, \hmax_{3} \, \rr}{6} \, \| \DFN(\upsv) \uv \|^{2} \, .
\label{jrgeteteer2234587654}
\end{EQA}
This implies \nameref{LL3tref} as well as \eqref{gtcdsftdfvtwdsefhfdvfrvsewseGP}; see \eqref{dta3u1DG2d3GPa1g}.
The statement about \nameref{LL4tref} is similar.
\end{proof}

%
%

\Chapter{Schur complement}
Consider a symmetric \( \dimttl \times \dimttl \) matrix \( \F \) with block representation 
\begin{EQA}[c]
	\F
	=
	\begin{pmatrix}
	\F_{\tarpv\tarpv} & \F_{\tarpv\nupv} \\
	\F_{\nupv\tarpv} & \F_{\nupv\nupv}
	\end{pmatrix} .
\label{HAuuuvvuvvTHS}
\end{EQA} 

\begin{lemma}
\label{LSchur}
Let the diagonal blocks \( \F_{\tarpv\tarpv} , \F_{\nupv\nupv} \) of \( \F \) be positive definite.
Define
\begin{EQA}[rcccl]
	\Fb_{\tarpv\tarpv} 
	& \eqdef & 
	\F_{\tarpv\tarpv} - \F_{\tarpv\nupv} \, \F_{\nupv\nupv}^{-1} \, \F_{\nupv\tarpv} ,
	\qquad
	\Fb_{\nupv\nupv}
	& \eqdef &
	\F_{\nupv\nupv} - \F_{\nupv\tarpv} \, \F_{\tarpv\tarpv}^{-1} \, \F_{\tarpv\nupv} 
	\, .
\label{y7djw38vjer4yw3jfclweus}
\end{EQA}
If \( \Fb_{\tarpv\tarpv} \) or \( \Fb_{\nupv\nupv} \) is also positive definite then \( \F \)
is positive definite as well.
It holds 
\begin{EQA}
	\begin{pmatrix}
		\F_{\tarpv\tarpv} & \F_{\tarpv\nupv} \\
		\F_{\nupv\tarpv} & \F_{\nupv\nupv}
	\end{pmatrix}^{-1}
	&=&
	\begin{pmatrix}
		\Id_{\dimp} & 0 \\
		- \F_{\nupv\nupv}^{-1} \, \F_{\nupv\tarpv} &\Id_{\dimq}
	\end{pmatrix}
	\,\,
	\begin{pmatrix}
		\Fb_{\tarpv\tarpv}^{-1} & 0 \\
		0 	& \F_{\nupv\nupv}^{-1}
	\end{pmatrix}
	\,\,
	\begin{pmatrix}
		\Id_{\dimp} &  - \F_{\tarpv\nupv} \, \F_{\nupv\nupv}^{-1} \\
		0 &	\Id_{\dimq}
	\end{pmatrix}
	\\
	&=&
	\begin{pmatrix}
		\Fb_{\tarpv\tarpv}^{-1} & - \Fb_{\tarpv\tarpv}^{-1} \, \F_{\tarpv\nupv} \, \F_{\nupv\nupv}^{-1} \\
		- \F_{\nupv\nupv}^{-1} \, \F_{\nupv\tarpv} \, \Fb_{\tarpv\tarpv}^{-1} & 
		\F_{\nupv\nupv}^{-1} + \F_{\nupv\nupv}^{-1} \, \F_{\nupv\tarpv} \, \Fb_{\tarpv\tarpv}^{-1} \, \F_{\tarpv\nupv} \, \F_{\nupv\nupv}^{-1} 
	\end{pmatrix}	
	\, 
\label{8jdkvtwfe6xhejdcedvscy}
\end{EQA}
and 
\begin{EQA}
	\begin{pmatrix}
		\F_{\tarpv\tarpv} & \F_{\tarpv\nupv} \\
		\F_{\nupv\tarpv} & \F_{\nupv\nupv}
	\end{pmatrix}^{-1}
	&=&
	\begin{pmatrix}
		\F_{\tarpv\tarpv}^{-1} + \F_{\tarpv\tarpv}^{-1} \, \F_{\tarpv\nupv} \, \Fb_{\nupv\nupv}^{-1} \, \F_{\nupv\tarpv} \, \F_{\tarpv\tarpv}^{-1} & - \F_{\tarpv\tarpv}^{-1} \, \F_{\tarpv\nupv} \, \Fb_{\nupv\nupv}^{-1} \\
		- \Fb_{\nupv\nupv}^{-1} \, \F_{\nupv\tarpv} \, \F_{\tarpv\tarpv}^{-1} & \Fb_{\nupv\nupv}^{-1} 
	\end{pmatrix}	
	\, .
\label{8jdkvtwfe6xhejdcedvscye}
\end{EQA}
In particular, this implies 
\( \Fb_{\tarpv\tarpv}^{-1} \, \F_{\tarpv\nupv} \, \F_{\nupv\nupv}^{-1} 
\equiv \F_{\tarpv\tarpv}^{-1} \, \F_{\tarpv\nupv} \, \Fb_{\nupv\nupv}^{-1} \),
\begin{EQA}
\label{t6rygfujio789uiew45ygi}
	\F_{\tarpv\tarpv}^{-1} + \F_{\tarpv\tarpv}^{-1} \, \F_{\tarpv\nupv} \, \Fb_{\nupv\nupv}^{-1} \, \F_{\nupv\tarpv} \, \F_{\tarpv\tarpv}^{-1}
	& \equiv &
	\Fb_{\tarpv\tarpv}^{-1}
	\, ,
	\\
	\F_{\nupv\nupv}^{-1} + \F_{\nupv\nupv}^{-1} \, \F_{\nupv\tarpv} \, \Fb_{\tarpv\tarpv}^{-1} \, \F_{\tarpv\nupv} \, \F_{\nupv\nupv}^{-1}
	& \equiv &
	\Fb_{\nupv\nupv}^{-1} 
	\, .
\end{EQA}
Moreover, for any \( \wv = (\tarpv,\nupv) \in \R^{\dimttl} \), it holds \( \| \F^{1/2} \, \wv \| \geq \| \Fb_{\nupv\nupv}^{1/2} \, \tarpv \| \) and
\begin{EQA}
	\| \F^{1/2} \, \wv \|^{2} 
	&=& 
	\| \Fb_{\nupv\nupv}^{1/2} \, \tarpv \|^{2} 
	+ \| \F_{\nupv\nupv}^{1/2} (\nupv - \F_{\nupv\nupv}^{-1} \, \F_{\nupv\tarpv} \, \tarpv) \|^{2} 
\label{jf76fugfyf5te345e4ghvhiu}
	\\
	\| \F^{-1/2} \, \wv \|^{2}
	&=&
	\| \Fb_{\tarpv\tarpv}^{-1/2} \, (\tarpv - \F_{\tarpv\nupv} \, \F_{\nupv\nupv}^{-1} \nupv)  \|^{2} 
	+ \| \F_{\nupv\nupv}^{-1/2} \nupv \|^{2} ;
\label{7ycdkjic8e38dd98le9dlo}
	\\
	\bigl( \F^{-1} \wv \bigr)_{\tarpv}  
	&=& 
	\Fb_{\tarpv\tarpv}^{-1} (\tarpv - \F_{\tarpv\nupv} \, \F_{\nupv\nupv}^{-1} \nupv) 
	=
	\Fb_{\tarpv\tarpv}^{-1} \tarpv - \F_{\tarpv\tarpv}^{-1} \, \F_{\tarpv\nupv} \, \Fb_{\nupv\nupv}^{-1} \nupv 
	\, .
\label{poryerjnjvyt65e6yjer}
\end{EQA}
Furthermore, suppose 
\begin{EQA}
	\| 
	\F_{\tarpv\tarpv}^{-1/2} \, \F_{\tarpv\nupv} \, \F_{\nupv\nupv}^{-1} \, \F_{\nupv\tarpv} \, \F_{\tarpv\tarpv}^{-1/2} 
	\|
	& \leq &
	\rhoIF^{2}
	< 
	1 .
\label{vcjcfvedtesqgghwqLco}
\end{EQA}
Then it holds for \( \F_{0} \eqdef \blk\{ \F_{\tarpv\tarpv},\F_{\nupv\nupv} \} \)
\begin{EQA}
	(1 - \rhoIF) \F_{0}
	\leq 
	\F
	& \leq &
	(1 + \rhoIF) \F_{0}
\label{f6eh3ewfd6ehev65r43tre}
\end{EQA}
and also
\begin{EQA}
	(1 - \rhoIF^{2}) \, \F_{\tarpv\tarpv}
	\leq 
	\Fb_{\tarpv\tarpv}
	& \leq &
	\F_{\tarpv\tarpv} \, ,
	\qquad
	(1 - \rhoIF^{2}) \, \F_{\nupv\nupv}
	\leq 
	\Fb_{\nupv\nupv}
	\leq 
	\F_{\nupv\nupv} \, .
\label{yg3w5dffctvyry4r7er7dgfe}
\end{EQA}
\end{lemma}

\begin{proof}
The block inversion follows by Schur's complement formula; see e.g. \cite{Boyd2004}[Appendix A.5.5].
Minimizing \( \| \F^{1/2} \, \wv \|^{2} = \tarpv^{\T} \, \F_{\tarpv\tarpv} \, \tarpv + 2 \tarpv^{\T} \, \F_{\tarpv\nupv} \, \nupv + \nupv^{\T} \, \F_{\nupv\nupv} \, \nupv \) 
w.r.t. \( \nupv \) leads to \( \nupv_{0} = - \F_{\nupv\nupv}^{-1} \, \F_{\nupv\tarpv} \tarpv \) 
and by quadraticity of \( \| \F^{1/2} \, \wv \|^{2} \) in \( \nupv \)
\begin{EQA}
	\| \F^{1/2} \, \wv \|^{2} 
	&=& 
	\tarpv^{\T} \, \F_{\tarpv\tarpv} \, \tarpv + 2 \tarpv^{\T} \, \F_{\tarpv\nupv} \, \nupv_{0} + \nupv_{0}^{\T} \, \F_{\nupv\nupv} \, \nupv_{0} 
	+ \| \F_{\nupv\nupv}^{1/2} (\nupv - \nupv_{0}) \|^{2} 
	\\
	&=&
	\| \Fb_{\nupv\nupv}^{1/2} \, \tarpv \|^{2} + \| \F_{\nupv\nupv}^{1/2} (\nupv - \nupv_{0}) \|^{2} .
\label{fji43324efgedygffhbu}
\end{EQA} 
This proves \eqref{jf76fugfyf5te345e4ghvhiu}.
Further, represent \( \F^{-1} \) using Gauss elimination \eqref{8jdkvtwfe6xhejdcedvscy}:
\begin{EQA}
	\F^{-1}
	&=&
	\begin{pmatrix}
		\Id_{\dimp} & 0 \\
		- \F_{\nupv\nupv}^{-1} \, \F_{\nupv\tarpv} &	\Id_{\dimq}
	\end{pmatrix}
	\,\,
	\begin{pmatrix}
		\Fb_{\tarpv\tarpv}^{-1} & 0 \\
		0 	& \F_{\nupv\nupv}^{-1}
	\end{pmatrix}
	\,\,
	\begin{pmatrix}
		\Id_{\dimp} &  - \F_{\tarpv\nupv} \, \F_{\nupv\nupv}^{-1} \\
		0 &	\Id_{\dimq}
	\end{pmatrix} \, .
\label{ucjhdcytew3hjf64ehfkesu}
\end{EQA}
Then 
\begin{EQA}
	\wv^{\T} \F^{-1} \wv 
	&=& 
	\begin{pmatrix}
		\tarpv - \F_{\tarpv\nupv} \, \F_{\nupv\nupv}^{-1} \nupv \\
		\nupv 
	\end{pmatrix}^{\T} \, 
	\begin{pmatrix}
		\Fb_{\tarpv\tarpv}^{-1} & 0 \\
		0 	& \F_{\nupv\nupv}^{-1}
	\end{pmatrix}
	\begin{pmatrix}
		\tarpv - \F_{\tarpv\nupv} \, \F_{\nupv\nupv}^{-1} \nupv \\
		\nupv
	\end{pmatrix} ,
\label{ucjd78renjg8we3vjfe}
\end{EQA}
and \eqref{7ycdkjic8e38dd98le9dlo} follows.
Also \eqref{8jdkvtwfe6xhejdcedvscy} implies \eqref{poryerjnjvyt65e6yjer}.

Next, define \( \F_{0} = \blk\{ \F_{\tarpv\tarpv},\F_{\nupv\nupv} \} \), 
\( U = \F_{\tarpv\tarpv}^{-1/2} \, \F_{\tarpv\nupv} \, \F_{\nupv\nupv}^{-1/2} \), and consider the matrix 
\begin{EQA}
	\F_{0}^{-1/2} \, \F \, \F_{0}^{-1/2} - \Id_{\dimttl}
	&=&
	\begin{pmatrix}
		0 & \F_{\tarpv\tarpv}^{-1/2} \, \F_{\tarpv\nupv} \, \F_{\nupv\nupv}^{-1/2}
		\\
		\F_{\nupv\nupv}^{-1/2} \, \F_{\nupv\tarpv} \, \F_{\tarpv\tarpv}^{-1/2} & 0
	\end{pmatrix}
	=
	\begin{pmatrix}
		0 & U
		\\
		U^{\T} & 0
	\end{pmatrix} .
\label{gswrew35e35e35e5txzgtqw}
\end{EQA}
Condition \eqref{vcjcfvedtesqgghwqLco} implies \( \| U U^{\T} \| \leq \rhoIF^{2} \) and hence,
\begin{EQA}
	- \rhoIF \, \Id_{\dimttl} 
	& \leq &
	\F_{0}^{-1/2} \, \F \, \F_{0}^{-1/2} - \Id_{\dimttl} 
	\leq 
	\rhoIF \, \Id_{\dimttl} \, .
\label{ye376evhgder52wedsytg}
\end{EQA}
Moreover,
\begin{EQA}
	\Fb_{\tarpv\tarpv} 
	& = &
	\F_{\tarpv\tarpv} - \F_{\tarpv\nupv} \, \F_{\nupv\nupv}^{-1} \, \F_{\nupv\tarpv} 
	=
	\F_{\tarpv\tarpv}^{1/2} (\Id_{\dimp} - U U^{\T}) \F_{\tarpv\tarpv}^{1/2} 
	\geq 
	(1 - \rhoIF^{2}) \F_{\tarpv\tarpv} \, ,
\label{ljhy6furf4jfu8rdfcweerdd}
\end{EQA}
and similarly for \( \Fb_{\nupv\nupv} \).
\end{proof}

For some situations, the nuisance variable \( \nuiv \) is by itself a composition of a few other subvectors.
We only consider the case of two variables \( \nuiv = (\zv,\nuov) \).
Denote by \( \F_{\targv\targv} \), \( \F_{\targv\zv} \), \( \F_{\targv\nuov} \), 
\( \F_{\zv\zv} \), \( \F_{\zv\nuov} \), \( \F_{\nuov\nuov} \) the corresponding blocks of \( \F \), that is,
\begin{EQA}
	\F
	&=&
	\begin{pmatrix}
		\F_{\targv\targv} & \F_{\targv\zv} & \F_{\targv\nuov} \\
		\F_{\zv\targv} & \F_{\zv\zv} & \F_{\zv\nuov} \\
		\F_{\nuov\targv} & \F_{\nuov\zv} & \F_{\nuov\nuov}
	\end{pmatrix} .
\label{576vner6734nf9h4ede}
\end{EQA}

\begin{lemma}
\label{LblocksIFT}
For the matrix \( \F \) from \eqref{576vner6734nf9h4ede}, suppose that 
\begin{EQ}[rcl]
	\| \F_{\targv\targv}^{-1/2} \F_{\targv\zv} \, \F_{\zv\zv}^{-1/2} \|
	& \leq &
	\rhoIF_{\targv\zv} \, ,
	\\
	\| \F_{\targv\targv}^{-1/2} \F_{\targv\nuov} \, \F_{\nuov\nuov}^{-1/2} \|
	& \leq &
	\rhoIF_{\targv\nuov} \, ,
	\\
	\| \F_{\nuov\nuov}^{-1/2} \F_{\nuov\zv} \, \F_{\zv\zv}^{-1/2} \|
	& \leq &
	\rhoIF_{\zv\nuov} \, ,
\label{vcjcfvedt7wdwhesqgLcompn}
\end{EQ}
and \( \max\{ \rhoIF_{\targv\zv} + \rhoIF_{\targv\nuov} \, , \rhoIF_{\targv\zv} + \rhoIF_{\zv\nuov} \, , 
\rhoIF_{\targv\nuov} + \rhoIF_{\zv\nuov} \} \leq 1 \).
Then 
\begin{EQA}
	\F
	& \geq &
	\begin{pmatrix}
		(1 - \rhoIF_{\targv\zv} - \rhoIF_{\targv\nuov}) \F_{\targv\targv} & 0 & 0 \\ 
		0 & (1 - \rhoIF_{\targv\zv} - \rhoIF_{\zv\nuov}) \F_{\zv\zv} & 0 \\
		0 & 0 & (1 - \rhoIF_{\targv\nuov} - \rhoIF_{\zv\nuov}) \F_{\nuov\nuov}
	\end{pmatrix}
\label{kjygy5t55gt2fc7w2qdxkwqiu}
\end{EQA}
\end{lemma}

\begin{proof}
By block-normalization we can reduce the proof to the case \( \F_{\targv\targv} = \Id \), \( \F_{\zv\zv} = \Id \),
\( \F_{\nuov\nuov} = \Id \).
Then it suffices to check positive semi-definiteness of the matrix
\begin{EQA}
	\BBH
	&=&
	\begin{pmatrix}
		(\rhoIF_{\targv\zv} + \rhoIF_{\targv\nuov}) \Id & \F_{\targv\zv} & \F_{\targv\nuov} \\
		\F_{\zv\targv} & (\rhoIF_{\targv\zv} + \rhoIF_{\zv\nuov}) \Id & \F_{\zv\nuov} \\
		\F_{\nuov\targv} & \F_{\nuov\zv} & (\rhoIF_{\targv\nuov} + \rhoIF_{\zv\nuov}) \Id
	\end{pmatrix}
\label{hythgwncuyuewuywhxjkqg}
\end{EQA}
with \( \| \F_{\targv\zv} \F_{\zv\targv} \| \leq \rhoIF_{\targv\zv}^{2} \),
\( \| \F_{\targv\nuov} \F_{\nuov\targv} \| \leq \rhoIF_{\targv\nuov}^{2} \),
\( \| \F_{\zv\nuov} \F_{\nuov\zv} \| \leq \rhoIF_{\zv\nuov}^{2} \).
Such a matrix is positive semi-definite by Gershgorin theorem or general results for diagonal dominant matrices.
We, however, present a simple proof.
For any \( \prmtv = (\targv,\zv,\nuov) \), we can use
\begin{EQA}
	2 |\targv^{\T} \F_{\targv\zv} \zv|
	& \leq &
	2 \| \targv \| \, \| \F_{\targv\zv} \zv \|
	\leq 
	2 \rhoIF_{\targv\zv} \| \targv \| \, \| \zv \|
	\leq 
	\rhoIF_{\targv\zv} (\| \targv \|^{2} + \| \zv \|^{2})
\label{gbtythhxct6wwnhxcd76}
\end{EQA}
and similarly for \( 2 \targv^{\T} \F_{\targv\nuov} \nuov \) and \( 2 \zv^{\T} \F_{\zv\nuov} \nuov \).
Therefore,
\begin{EQA}
	\prmtv^{\T} \BBH \prmtv
	&=&
	(\rhoIF_{\targv\zv} + \rhoIF_{\targv\nuov}) \| \targv \|^{2}
	+ (\rhoIF_{\targv\zv} + \rhoIF_{\zv\nuov}) \| \zv \|^{2}
	+ (\rhoIF_{\targv\nuov} + \rhoIF_{\zv\nuov}) \| \nuov \|^{2}
	\\
	&&
	+ \, 2 \targv^{\T} \F_{\targv\zv} \zv
	+ 2 \targv^{\T} \F_{\targv\nuov} \nuov
	+ 2 \zv^{\T} \F_{\zv\nuov} \nuov
	\\
	& \geq &
	(\rhoIF_{\targv\zv} + \rhoIF_{\targv\nuov}) \| \targv \|^{2}
	+ (\rhoIF_{\targv\zv} + \rhoIF_{\zv\nuov}) \| \zv \|^{2}
	+ (\rhoIF_{\targv\nuov} + \rhoIF_{\zv\nuov}) \| \nuov \|^{2}
	\\
	&&
	- \, \rhoIF_{\targv\zv} (\| \targv \|^{2} + \| \zv \|^{2})
	- \rhoIF_{\targv\nuov} (\| \targv \|^{2} + \| \nuov \|^{2})
	- \rhoIF_{\zv\nuov} (\| \zv \|^{2} + \| \nuov \|^{2})
	= 
	0
\label{bhxc4rt2yer6r67t478ds}
\end{EQA}
and the assertion follows.
\end{proof}

In some cases, it is more convenient to apply another block-diagonal matrix \( \DFN = (\DFN_{\targv\targv},\DFN_{\zv\zv},\DFN_{\nuov\nuov}) \).
\begin{lemma}
\label{LblocksDFN}
For the matrix \( \F \) from \eqref{576vner6734nf9h4ede}, suppose that 
\begin{EQ}[rcccccl]
	\| \DFN_{\targv\targv}^{-1} \, \F_{\targv\zv} \, \DFN_{\zv\zv}^{-1} \|
	& \leq &
	\aIF_{\targv\zv} \, ,
	\quad
	\| \DFN_{\targv\targv}^{-1} \, \F_{\targv\nuov} \, \DFN_{\nuov\nuov}^{-1} \|
	& \leq &
	\aIF_{\targv\nuov} \, ,
	\quad
	\| \DFN_{\nuov\nuov}^{-1} \, \F_{\nuov\zv} \, \DFN_{\zv\zv}^{-1} \|
	& \leq &
	\aIF_{\zv\nuov} \, .
	\qquad
\label{vcjcfvedt7wdwhesqgLcompnD}
\end{EQ}
Let also
\begin{EQ}[rcccccl]
	\| \DFN_{\targv\targv}^{-1} \, \F_{\targv\targv} \, \DFN_{\targv\targv}^{-1} \|
	& \geq &
	\bIF_{\targv\targv}^{2} \, ,
	\quad
	\| \DFN_{\zv\zv}^{-1} \, \F_{\zv\zv} \, \DFN_{\zv\zv}^{-1} \|
	& \geq &
	\bIF_{\zv\zv}^{2} \, ,
	\quad
	\| \DFN_{\nuov\nuov}^{-1} \, \F_{\nuov\nuov} \, \DFN_{\nuov\nuov}^{-1} \|
	& \geq &
	\bIF_{\nuov\nuov}^{2} \, ,
	\qquad
\label{vcjcfvedt7wdwhesqgLcompnd}
\end{EQ}
and 
\begin{EQ}[rcl]
	\bIF_{\targv\targv}^{2}  
	- \frac{\bIF_{\targv\targv} \, \aIF_{\targv\zv}}{\bIF_{\zv\zv}} - \frac{\bIF_{\targv\targv} \, \aIF_{\targv\nuov}}{\bIF_{\nuov\nuov}}
	& \geq &
	\dmax^{-2} \, ,
	\\
	\bIF_{\zv\zv}^{2}  
	- \frac{\bIF_{\zv\zv} \, \aIF_{\targv\zv}}{\bIF_{\targv\targv}} - \frac{\bIF_{\zv\zv} \, \aIF_{\zv\nuov}}{\bIF_{\nuov\nuov}}
	& \geq &
	\dmax^{-2} \, ,
	\\
	\bIF_{\nuov\nuov}^{2}  
	- \frac{\bIF_{\nuov\nuov} \, \aIF_{\targv\nuov}}{\bIF_{\zv\zv}} - \frac{\bIF_{\nuov\nuov} \, \aIF_{\zv\nuov}}{\bIF_{\zv\zv}}
	& \geq &
	\dmax^{-2} \, .
\label{dyc7f7yddfydrhye35wsss}
\end{EQ}
Then 
\begin{EQA}
	\F
	& \geq &
	\dmax^{-2} \DFM^{2} \, .
\label{kjygy5t55gt2fc7w2qDFM}
\end{EQA}
\end{lemma}

\begin{proof}
Conditions of the lemma imply \eqref{vcjcfvedt7wdwhesqgLcompn} with 
\begin{EQA}[c]
	\rhoIF_{\targv\zv}
	=
	\frac{\aIF_{\targv\zv}}{\bIF_{\targv\targv} \, \bIF_{\zv\zv}} \, , 
	\quad
	\rhoIF_{\targv\nuov}
	=
	\frac{\aIF_{\targv\nuov}}{\bIF_{\targv\targv} \, \bIF_{\nuov\nuov}} \, , 
	\quad
	\rhoIF_{\nuov\zv}
	=
	\frac{\aIF_{\nuov\zv}}{\bIF_{\nuov\nuov} \, \bIF_{\zv\zv}} \, , \quad
\label{d8dcjed7yw3md7ujdsdrdlk}
\end{EQA}
Now we apply Lemma~\ref{LblocksIFT} and note that by \eqref{dyc7f7yddfydrhye35wsss}
\begin{EQA}
	(1 - \rhoIF_{\targv\zv} - \rhoIF_{\targv\nuov}) \F_{\targv\targv} - \dmax^{-2} \DFM_{\targv\targv}^{2}
	& \geq &
	(1 - \rhoIF_{\targv\zv} - \rhoIF_{\targv\nuov}) \DFM_{\targv\targv}^{2} \bIF_{\targv\targv}^{2} - \dmax^{-2} \DFM_{\targv\targv}^{2}
	\geq 
	0 \, ,
	\\
	(1 - \rhoIF_{\targv\zv} - \rhoIF_{\zv\nuov}) \F_{\zv\zv} - \dmax^{-2} \DFM_{\zv\zv}^{2}
	& \geq &
	(1 - \rhoIF_{\targv\zv} - \rhoIF_{\zv\nuov}) \DFM_{\zv\zv}^{2} \bIF_{\zv\zv}^{2} - \dmax^{-2} \DFM_{\zv\zv}^{2}
	\geq 
	0 \, ,
	\\
	(1 - \rhoIF_{\targv\nuov} - \rhoIF_{\zv\nuov}) \F_{\nuov\nuov} - \dmax^{-2} \DFM_{\nuov\nuov}^{2}
	& \geq &
	(1 - \rhoIF_{\targv\nuov} - \rhoIF_{\zv\nuov}) \DFM_{\nuov\nuov}^{2} \bIF_{\nuov\nuov}^{2} - \dmax^{-2} \DFM_{\nuov\nuov}^{2}
	\geq 
	0 \, .
\label{ychywy7huyr78wjduy7djws}
\end{EQA}
This implies the assertion.
\end{proof}

{
\renewcommand{\Section}{\section}
\renewcommand{\Subsection}{\subsection}



\def\prmtvb{\accentset{\circ}{\prmtv}}
\def\tarpvb{\accentset{\circ}{\tarpv}}
\def\KSb{\accentset{\circ}{\KS}}
\def\zvb{\accentset{\circ}{\zv}}

\Section{Error-in-Operator model. Tools}
\label{SEiOtools}

This section collects some useful facts about the error-in-operator model.

\Subsection{Full dimensional information matrix and identifiability}
\label{SinfidentEO}
Given an image vector \( \Imav \) and a pilot \( \hKS \) of the operator \( \KS \),
consider the random function of the full parameter \( \prmtv = (\tarpv,\imav,\KS) \)
\begin{EQA}
	\fs(\prmtv)
	&=&
	- \frac{1}{2} \| \Imav - \imav \|^{2} 
	- \frac{1}{2} \| \imav - \KS \tarpv \|^{2} 
	- \frac{\muA^{2}}{2} \| \hKS - \KS \|_{\Fr}^{2} 
	\, .
	\qquad
\label{bh2gfefth3wqe3rEO}
\end{EQA} 

The next lemma describes the blocks of the matrix \( \IFT(\prmtv) = - \nabla^{2} \fs(\prmtv) \).

\begin{lemma}
\label{Lnabl2EO}
It holds for \( \prmtv = (\tarpv,\imav,\KS) \)
\begin{EQA}
	\IFT_{\tarpv\tarpv}(\prmtv)
	& = &
	- \nabla_{\tarpv\tarpv}^{2} \fs(\tarpv,\imav,\KS)
	=
	\KS^{\T} \KS 
	=
	\sum_{\mm=1}^{\dimq} \KS_{\mm} \, \KS_{\mm}^{\T} \, ,
\label{un4v7frue43kgv8uhkjEO}
	\\
	\IFT_{\imav\imav}(\prmtv)
	& = &
	- \nabla_{\imav\imav}^{2} \fs(\tarpv,\imav,\KS)
	\equiv 
	2 \Id_{\dimq} \, ,
\label{8dfnjey63f6ytjutz}
	\\
	\IFT_{\tarpv\imav}(\prmtv)
	& = &
	- \nabla_{\tarpv} \nabla_{\imav} \fs(\tarpv,\imav,\KS)
	=
	- \KS^{\T} \, .
	\qquad
\label{8dfnjey63f6ytjutzv}
\end{EQA}
The matrix \( \IFT_{\KS\KS}(\prmtv) \) is block-diagonal:
\begin{EQ}[rcl]
	\IFT_{\KS\KS}(\prmtv)
	&=& 
	\blk\bigl\{ \IFT_{\KS_{1} \KS_{1}}(\prmtv), \ldots,\IFT_{\KS_{\dimq} \KS_{\dimq}}(\prmtv) \bigr\} \, ,
	\\
	\IFT_{\KS_{\mm} \KS_{\mm}}(\prmtv) 
	& = &
	- \nabla_{\KS_{\mm} \KS_{\mm}}^{2} \fs(\tarpv,\imav,\KS)
	=
	\tarpv \tarpv^{\T} + \muA^{2} \Id_{\dimp} 
	\, .
\label{8dfnjey63swlpwwqEO}
\end{EQ}
If \( \ev_{\mm} \) being the \( \mm \)th basis vector in \( \R^{\dimq} \), then 
\begin{EQ}[rcl]
\label{7mdfgurt564634yufrt6re}
	\IFT_{\tarpv\KS_{\mm}}(\prmtv)
	= 
	- \nabla_{\tarpv} \nabla_{\KS_{\mm}} \fs(\tarpv,\imav,\KS)
	&=&
	(\KS_{\mm}^{\T} \tarpv - \imav_{\mm}) \Id_{\dimp} + \KS_{\mm} \, \tarpv^{\T} 
	\eqdef 
	\AFblk_{\mm}(\tarpv,\KS_{\mm}) \, ,
	\qquad
	\\
	\IFT_{\imav\KS_{\mm}}(\prmtv)
	=
	- \nabla_{\imav} \nabla_{\KS_{\mm}} \fs(\tarpv,\imav,\KS) 
	&=& 
	- \ev_{\mm} \, \tarpv^{\T} 
	\, .
\end{EQ}
and \( \IFT_{\tarpv\imav}(\prmtv) \IFT_{\imav\KS_{\mm}}(\prmtv) = \KS_{\mm} \tarpv^{\T} \).
\end{lemma}

\begin{lemma}
\label{Lpart2EiO}
Let \( \fs(\prmtv) \) be given by \eqref{bh2gfefth3wqe3rEO} for \( \prmtv = (\tarpv,\imav,\KS) \), 
and  \( \IFT(\prmtv) = \nabla^{2} \fs(\prmtv) = \nabla^{2} \fs(\tarpv,\imav,\KS) \).
For any \( \dprmtv = (\dtarpv,\dimav,\dKS) \) with \( \dtarpv \in \R^{\dimp} \), \( \dimav \in \R^{\dimq} \), and
\( \dKS \in \R^{\dimq \times \dimp} \), 
\begin{EQA}[rcl]
	\bigl\langle \IFT(\prmtv), \uv^{\otimes 2} \bigr\rangle
	&=&
	\frac{d^{2}}{dt^{2}} \fs(\prmtv + t \dprmtv) \bigg|_{t=0}
	\\
	&=&
	\| \dimav \|^{2} 
	+ \muA^{2} \| \dKS \|_{\Fr}^{2} 
	+ \| \dimav - (\KS \dtarpv + \dKS \tarpv) \|^{2} + 2 (\imav - \KS \tarpv)^{\T} \dKS \dtarpv
	\, .
	\qquad
\label{dyj3wvcrtwhvcwc33huh}
\end{EQA}
\end{lemma}

\begin{proof}
The proof is straightforward. 
The only non-quadratic term \( (\imav - \KS \tarpv)^{\T} \dKS \dtarpv \) corresponds to the second derivative of the product
\( (\imav + t \dimav) (\KS + t \dKS) \) multiplied with the discrepancy \( \imav - \KS \tarpv \).
\end{proof}

Let
\(
	\DPN^{2} 
	\geq 
	{\KSn}^{\T} \KSn \),
	\(
	\nEO
	= 
	\lambda_{\min}(\DPN^{2})
\).
With \( \eEO \leq 1/10 \), define the radius \( \REO = \muA \sqrt{\nEO} \eEO \) and
consider the local set \( \Upsd \)
\begin{EQA}
	\Upsd
	& = &
	\bigl\{ (\tarpv,\imav,\KS) \colon 
		\| \DPN \tarpv \| \leq \REO, \,
		\| \imav \| \leq \REO, \,
		\| (\KS - \KSn) \DPN^{-1} \|_{\Fr} \leq \eEO
	\bigr\} \, .
	\qquad
\label{rughur47478ytfg84g6ght}
\end{EQA}

Now we bound from below the Hessian of \( \fs(\prmtv) \) from \eqref{bh2gfefth3wqe3rEO} for \( \prmtv \in \Upsd \).

\begin{lemma}
\label{LIFTDFMEO}
With \( \DFM \) from \eqref{iwdfhw787e3urne6fyh3o} and \( \dmax = 2 \), it holds for any \( \prmtv \in \Upsd \)
\begin{EQA}[rcl]
\label{vkr0roej83375e72hvirjrj}
	\IFT(\prmtv)
	& \geq &
	\dmax^{-2} \DFM^{2} 
	=
	\dmax^{-2} \blk\{ \DPN^{2}, \Id_{\dimq}, \muA^{2} \Id_{\KS} \}
	\, ,
	\\
	\IFT_{\nupv\nupv}(\prmtv)
	& \geq &
	\blk\{ \Id_{\dimq}, \muA^{2} \Id_{\KS} \}
	\, .
\label{usxsxwgd6e6464td7dhc}
\end{EQA}
\end{lemma}

\begin{proof}
The first step of the proof is 
a bound on the deficiency \( \KS \tarpv - \imav \) over \( \prmtv \in \Upsd \).
Namely, we show that 
for any \( \prmtv \in \Upsd \), it holds 
\begin{EQA}[c]
	\| \KS \tarpv - \imav \|
	\leq 
	(2 + \eEO) \REO \, ,
\label{e6fc6e6e53er63yefgt7}
	\\
	- (2 \eEO - \eEO^{2}) \Id_{\dimp}
	\leq 
	\DPN^{-1} (\KS^{\T} \KS - {\KSs}^{\T} \KSs) \DPN^{-1}
	\leq 
	(2 \eEO + \eEO^{2}) \Id_{\dimp} \, .
\label{sd7w774bgruir3863e3ef}
\end{EQA}
Indded, by \( \| \DPN^{-2} \| = 1/\nEO \) and \( \DPN^{2} \geq {\KSn}^{\T} \KSn \)
\begin{EQA}
	\| \KS \tarpv - \imav \|
	& \leq &
	\| \KSn \tarpv \| + \| \imav \| + \| (\KS - \KSn) \tarpv \|
	\leq 
	2 \REO + \| (\KS - \KSn) \DPN^{-1} \| \; \| \DPN \tarpv \|
	\leq 
	(2 + \eEO) \REO 
	\, ,
\label{dstf6wd6y76wdet7wfvg}
\end{EQA}
and \eqref{e6fc6e6e53er63yefgt7} follows.
Further, it holds by \( \| (\KS - \KSn) \, \DPN^{-1} \| \leq \eEO \) and \( \| \KSn \DPN^{-1} \| \leq 1 \)
\begin{EQA}
	&& \nquad
	\DPN^{-1} (\KS^{\T} \KS - {\KSn}^{\T} \KSn) \DPN^{-1} 
	\\
	& = &
	\DPN^{-1}  (\KS - \KSn)^{\T} (\KS - \KSn) \DPN^{-1} + \DPN^{-1} {\KSn}^{\T} (\KS - \KSn) \DPN^{-1} 
	 + \DPN^{-1} (\KS - \KSn)^{\T} \KSn \DPN^{-1}
\label{yd6frty3eywuy773673647}
\end{EQA}
and \eqref{sd7w774bgruir3863e3ef} follows as well. 

To simplify derivations, later we assume that \( \DPN^{2} = {\KSs}^{\T} \KSs \).
Let us fix \( \prmtv \in \Upsd \) and \( \dprmtv = (\dtarpv,\dimav,\dKS) \) with \( \dtarpv \in \R^{\dimp} \), \( \dimav \in \R^{\dimq} \), and
\( \dKS \in \R^{\dimq \times \dimp} \).
As \( \| \DPN \tarpv \| \leq \REO \) and \( \| \DPN^{-2} \| \leq \nEO \), it holds with \( \eEO = \REO/(\muA \nEO^{1/2}) \)
\begin{EQA}[c]
	\| \dKS \tarpv \|
	\leq 
	\frac{1}{\sqrt{\nEO}} \| \dKS \|_{\Fr} \| \DPN \tarpv \| 
	\leq 
	\frac{\REO}{\muA \sqrt{\nEO}} \, \muA \| \dKS \|_{\Fr}
	\leq 
	\eEO \, \muA \| \dKS \|_{\Fr}
	\, .
\label{cuyvufjft4rtghewq2q2}
\end{EQA}
By \eqref{e6fc6e6e53er63yefgt7}, \( \| \imav - \KS \tarpv \| \leq (2 + \eEO) \REO \) yielding 
\begin{EQA}[rcl]
	2 \bigl| (\imav - \KS \tarpv)^{\T} \dKS \dtarpv \bigr|
	& \leq &
	2  \| \imav - \KS \tarpv \| \; \| \dKS \dtarpv \|  
	\leq 
	2(2 + \eEO) \REO \; \| \dKS \dtarpv \|
	\, .
	\qquad
\label{cudfjhed5dc3dfc3dskjy7}
\end{EQA}
Further, for any \( a,b,c \) and \( \rho,\rho_{1} \geq 0 \), the inequalities
\begin{EQA}[rcl]
	(a + b + c)^{2} + \rho a^{2} 
	& = &
	\Bigl\{ \sqrt{1 + \rho} \,\, a + \frac{b+c}{\sqrt{1 + \rho}} \Bigr\}^{2} + \frac{\rho}{1+\rho} (b+c)^{2}
	\geq 
	\frac{\rho}{1+\rho} (b+c)^{2} ,
	\\
	(b + c)^{2} + \rho_{1} b^{2}
	&=&
	\{ (1 + \rho_{1})^{1/2} b + (1 + \rho_{1})^{-1/2} c \}^{2} + \frac{\rho_{1}}{1+\rho_{1}} c^{2}
	\geq 
	\frac{\rho_{1}}{1+\rho_{1}} c^{2}
\end{EQA}
imply 
\begin{EQA}[c]
	(a + b + c)^{2} + a^{2} + \frac{\rho \rho_{1}}{1+\rho} b^{2}
	\geq 
	(1 - \rho) a^{2} + \frac{\rho \rho_{1}}{(1+\rho)(1+\rho_{1})} c^{2}
	\, .
\end{EQA}
In particular, with \( \rho = 3/4 \) and \( \rho_{1} = 7 \)
\begin{EQA}[c]
	(a + b + c)^{2} + a^{2} + 3 b^{2}
	\geq 
	\frac{1}{4} a^{2} + \frac{3}{8} c^{2} \, .
\label{dyx6x5543egdydbddd}
\end{EQA}
This and \eqref{sd7w774bgruir3863e3ef} yield
\begin{EQA}[rcl]
	\| \dimav - (\KS \dtarpv + \dKS \tarpv) \|^{2} + \| \dimav \|^{2} + 3 \| \dKS \tarpv \|^{2}
	& \geq &
	\frac{1}{4} \| \dimav \|^{2} + \frac{3}{8} \| \KS \dtarpv \|^{2}
	\geq 
	\frac{1}{4} \| \dimav \|^{2} + \frac{3}{8} (1 - \eEO)^{2} \| \DPN \dtarpv \|^{2} 
	\, .
\end{EQA}
By \eqref{dyj3wvcrtwhvcwc33huh}, \eqref{cuyvufjft4rtghewq2q2}, and \eqref{cudfjhed5dc3dfc3dskjy7} 
\begin{EQA}[rcl]
	\bigl\langle \IFT(\prmtv), \uv^{\otimes 2} \bigr\rangle
	& \geq &
	\frac{1}{4} \| \dimav \|^{2} + \frac{3}{8} \| \KS \dtarpv \|^{2} - 3 \| \dKS \tarpv \|^{2}
	+ \muA^{2} \| \dKS \|_{\Fr}^{2} - 2(2 + \eEO) \REO \; \| \dKS \dtarpv \|
	\\
	& \geq &
	\frac{1}{4} \| \dimav \|^{2} + \frac{3}{8} (1 - \eEO)^{2} \| \DPN \dtarpv \|^{2} 
	+ (1 - 3 \eEO^{2}) \muA^{2} \| \dKS \|_{\Fr}^{2}
	- 2(2 + \eEO) \REO \; \| \dKS \dtarpv \| 
	\, .
\end{EQA}
Also, for any \( s > 0 \)
\begin{EQA}[c]
	2 \| \dKS \dtarpv \|
	\leq 
	2 \| \dKS \DPN^{-1} \| \, \| \DPN \dtarpv \| 
	\leq 
	s^{-1} \| \dKS \DPN^{-1} \|_{\Fr}^{2} + s \| \DPN \dtarpv \|^{2} 
	\, .
\end{EQA}
We fix \( s \) by \( (2 + \eEO) \REO \, s = 6 \eEO^{2} \).
Then it holds by \( \REO = \muA \sqrt{\nEO} \eEO \)
\begin{EQA}[c]
	\frac{(2 + \eEO) \REO}{s} \| \dKS \DPN^{-1} \|_{\Fr}^{2}
	\leq 
	\frac{(2 + \eEO) \REO}{s \nEO} \| \dKS \|_{\Fr}^{2}
	\leq 
	\frac{(2 + \eEO)^{2} \REO^{2}}{6 \nEO \muA^{2} \eEO^{2}} \muA^{2}\| \dKS \|_{\Fr}^{2}
	\leq 
	\frac{(2 + \eEO)^{2}}{6} \,  \muA^{2}\| \dKS \|_{\Fr}^{2} 
	\, .
\end{EQA}
The use of \( \eEO \leq 1/10 \) yields
\begin{EQA}[rcl]
	\bigl\langle \IFT(\prmtv), \uv^{\otimes 2} \bigr\rangle
	& \geq &
	\frac{1}{4} \| \dimav \|^{2} 
	+ \Bigl( \frac{3}{8} (1 - \eEO)^{2} - 6 \eEO^{2} \Bigr) \| \DPN \dtarpv \|^{2} 
	+ \Bigl( 1 - \eEO^{2} - \frac{(2 + \eEO)^{2}}{6} \Bigr) \muA^{2} \| \dKS \|_{\Fr}^{2}
	\\
	& \geq &
	\frac{1}{4} \| \DPN \dtarpv \|^{2} + \frac{1}{4} \| \dimav \|^{2} + \frac{1}{4} \muA^{2} \| \dKS \|_{\Fr}^{2}
\end{EQA}
and assertion \eqref{vkr0roej83375e72hvirjrj} follows.
Similarly,
\begin{EQA}
	\bigl\langle \IFT_{\nupv\nupv}(\prmtv), (\dimav,\dKS)^{\otimes 2} \bigr\rangle
	&=&
	\frac{d^{2}}{dt^{2}} \fs(\tarpv,\imav + t \dimav, \KS + t \dKS) \bigg|_{t=0}
	=
	\| \dimav \|^{2} 
	+ \muA^{2} \| \dKS \|_{\Fr}^{2} 
	+ \| \dimav - \dKS \tarpv \|^{2} 
	\, .
\end{EQA}
This proves \eqref{usxsxwgd6e6464td7dhc}.
\end{proof}

\Subsection{Local smoothness}
\label{SsmoothEO}
The next step is to check \nameref{LLsT3ref} and \nameref{LLsT4ref} for all \( \prmtv \in \Upsd \);
see \eqref{rughur47478ytfg84g6ght}.
With \( \DPN^{2} \geq {\KSs}^{\T} \KSs \), the local geometry at 
\( \prmtv = (\tarpv,\imav,\KS) \in \Upsd \) with \( \KS = \KSs \)
is described by \( \DFM \) from \eqref{iwdfhw787e3urne6fyh3o}.

\begin{lemma}
\label{LS3checkEO}
For \( \DFM \) from \eqref{iwdfhw787e3urne6fyh3o} und \( \Upsd \) from \eqref{rughur47478ytfg84g6ght},
\( \fs(\prmtv) \) from \eqref{bh2gfefth3wqe3rEO} fulfills 
\nameref{LLsT3ref} and \nameref{LLsT4ref} at any \( \prmtv \in \Upsd \) with \( \rr = \REO \),
\( \KS = \KSs \), and
\begin{EQA}
	\dltwu_{3}
	&=&
	4.5 \, \nEO^{-1/2} \muA^{-1} ,
	\qquad
	\dltwu_{4}
	=
	3 \nEO^{-1} \muA^{-2} \, .
\label{gdjddheebfyygredtgebegxgEO}
\end{EQA}
\end{lemma}

\begin{proof}
Fix \( \dprmtv = (\dtarpv,\dimav,\dKS) \) with \( \dtarpv \in \R^{\dimp} \), \( \dimav \in \R^{\dimq} \), and
\( \dKS \in \R^{\dimq \times \dimp} \) s.t.
\begin{EQA}
	\dprmtv^{\T} \DFM^{2} \, \dprmtv
	=
	\| \DPN \dtarpv \|^{2} 
	+ \| \dimav \|^{2}
	+ \muA^{2} \| \dKS \|_{\Fr}^{2} 
	& = &
	1 .
	\qquad
\label{f564hjfgujg8747hfy76eheEO}
\end{EQA}
Consider \( \fs(\prmtv + t \dprmtv) = \fs(\tarpv + t \dtarpv,\imav + t \dimav,\KS + t \dKS) \).
The fourth derivative \( \frac{d^{4}}{dt^{4}} \fs(\prmtv + t \dprmtv) \) does not depend on \( \prmtv \) and
it holds for any \( t \)
\begin{EQA}
	- \frac{d^{4}}{dt^{4}} \fs(\prmtv + t \dprmtv) 
	& = &
	- \langle \nabla^{4} \fs(\prmtv), \dprmtv^{\otimes 4} \rangle
	=
	12 \| \dKS \,\dtarpv \|^{2} \, .
	\qquad
\label{ywkdofcfyweyenvcuduye}
\end{EQA}
By \eqref{f564hjfgujg8747hfy76eheEO} and \( \DPN^{-2} \leq \nEO^{-1} \Id_{\dimp} \),
it holds 
\begin{EQA}
	4 \| \dKS \,\dtarpv \|^{2}
	& \leq &
	\muA^{-2} \nEO^{-1} \biggl( \| \DPN \dtarpv \|^{2} + \muA^{2} \nEO \| \dKS \, \DPN^{-1} \|_{\Fr}^{2} \biggr)^{2}
	\\
	& \leq &
	\muA^{-2} \nEO^{-1} \bigl( \| \DPN \dtarpv \|^{2} + \muA^{2} \| \dKS \|_{\Fr}^{2} \bigr)^{2} 
	\leq 
	\muA^{-2} \nEO^{-1} \, .
\label{ugfr5544eehhjkudjdebhEO}
\end{EQA}
This and \eqref{ywkdofcfyweyenvcuduye} easily imply \nameref{LLsT4ref} with \( \dltwu_{4} = 3 \muA^{-2} \nEO^{-1} \).

Now we check \nameref{LLsT3ref}.
As in \eqref{dyj3wvcrtwhvcwc33huh}
\begin{EQA}
	&& \nquad
	- \frac{d^{3}}{dt^{3}} \fs(\prmtv + t \dprmtv) \bigg|_{t=0}
	=
	3 (\dKS \tarpv + \KS \dtarpv - \dimav)^{\T} \dKS \dtarpv
	=
	3 \bigl\{ \dKS \tarpv + \KSs \dtarpv - \dimav + (\KS - \KSs) \dtarpv \bigr\}^{\T} \dKS \dtarpv
	\, .
\label{ywkdofcfyweyenvcuduye3}
\end{EQA}
Further, by \eqref{cuyvufjft4rtghewq2q2}, it holds
\( \| \dKS \tarpv \| \leq \eEO \, \muA \| \dKS \|_{\Fr} \),
definition \eqref{rughur47478ytfg84g6ght} implies 
\( \| (\KS - \KSn) \DPN^{-1} \|_{\Fr} \leq \eEO \) yielding by \eqref{f564hjfgujg8747hfy76eheEO}
\begin{EQA}
	\| (\KS - \KSs) \dtarpv \|
	& \leq &
	\| (\KS - \KSs) \DPN^{-1} \DPN \dtarpv \|
	\leq 
	\eEO \, \| \DPN \dtarpv \| 
	\, ,
\label{ygysy6srx5x4exradtw}
\end{EQA}
while \eqref{ugfr5544eehhjkudjdebhEO} ensures 
\( 4 \| \dKS \,\dtarpv \|^{2} \leq \muA^{-2} \nEO^{-1} \).
Hence,
\begin{EQA}
	&& \nquad
	\bigl| (\dKS \tarpv + \KS \dtarpv - \dimav)^{\T} \dKS \dtarpv \bigr|
	\leq 
	\bigl( 
		\| \dKS \tarpv \| + \| \KSs \dtarpv \| + \| \dimav \| 
		+ \| (\KS - \KSs) \dtarpv \|
	\bigr) \frac{1}{\muA \sqrt{\nEO}}
	\\
	& \leq &
	\frac{1}{\muA \sqrt{\nEO}} \bigl( \eEO \muA \| \dKS \|_{\Fr}  
	+ \| \dimav \| + (1 + \eEO) \, \| \DPN \dtarpv \| \bigr) 
	\\
	& \leq &
	\frac{1}{\muA \sqrt{\nEO}} \sqrt{2 + 2 \eEO} \bigl( \muA \| \dKS \|_{\Fr}^{2} + \| \KSs \dtarpv \|^{2} + \| \dimav \|^{2} \bigr)^{1/2} 
	\leq 
	\frac{1}{\muA \sqrt{\nEO}} \sqrt{2 + 2 \eEO} 
\label{v8k3f7dnjhesrwehfg7ye}
\end{EQA}
and \nameref{LLsT3ref} follows with \( \dltwu_{3} = 4.5 /(\muA \sqrt{\nEO} ) \).
\end{proof}

\Subsection{Proof of Lemma~\ref{LdimfullEO}}
\label{SscoreEO}
The full dimensional score vector \( \nabla \zeta = (0,\errZv,\muA \, \errKS) \) involves
the observation noise \( \errZv = \Imav - \E \Imav \) and operator noise 
\( \errKS = \muA (\hKS - \E \hKS) \).
The full-dimensional noise energy is measured by
\( \dimAfull_{\DFM} \eqdef \tr \Var(\DFM \, \IFT_{\GPT}^{-1} \nabla \zeta) \); 
see \eqref{vhyvw323sdwe4w3rdesf23EO}.
The bound \( \IFT_{\GPT}^{-1} = (\IFT + \GPT^{2})^{-1} \leq (\dmax^{-2} \DFM^{2} + \GPT^{2})^{-1} \) for \( \dmax = 2 \) and the block-diagonal structure 
\( \DFM^{2} = \blk\{ \DPN^{2}, \Id_{\dimq}, \muA^{2} \Id_{\KS} \} \) and 
\( \GPT^{2} = \blk\{ \GP^{2}, 0, \GPKS^{2} \} \) enable us
to evaluate the full effective dimension via the corresponding values for each of the parameters \( \imav, \KS \).
Denote \( \DFM_{\dmax,\GPT}^{2} \eqdef \DFM^{2} + \dmax^{2} \GPT^{2} \), so that 
\( \DFM_{\dmax,\GPT} \, \IFT_{\GPT}^{-1} \DFM_{\dmax,\GPT} \leq \dmax^{2} \Id_{\dimp} \).
Also, \( \DFM^{2} \leq \DFM_{\dmax,\GPT}^{2} \) yielding \( \| \DFM \xiv \| \leq \| \DFM_{\dmax,\GPT} \, \xiv \| \) for any \( \xiv \in \R^{\dimp} \).
Hence,
\begin{EQA}[rcl]
	\tr \Var(\DFM \, \IFT_{\GPT}^{-1} \nabla \zeta)
	&=&
	\E \| \DFM \, \IFT_{\GPT}^{-1} \nabla \zeta \|^{2}
	\leq 
	\E \| \DFM_{\dmax,\GPT} \, \IFT_{\GPT}^{-1} \, \DFM_{\dmax,\GPT} \, \DFM_{\dmax,\GPT}^{-1} \nabla \zeta \|^{2}
	\\
	& \leq &
	\dmax^{4} \E \| \DFM_{\dmax,\GPT}^{-1} \nabla \zeta \|^{2}
	\leq 
	\dmax^{4} \tr \Var(\DFM_{\dmax,\GPT}^{-1} \nabla \zeta)	
	\\
	&=&
	\dmax^{4} \tr \Var(\errZv) + \dmax^{4} \tr \Var\bigl\{ (\muA^{2} \Id_{\KS} + \dmax^{2} \GPKS^{2})^{-1/2} \muA \errKS \bigr\}
\end{EQA} 
and \eqref{ydfhuwir8ufug87g73jdeEO} follows.
For homogeneous observation noise \( \Var(\errZ_{\mm}) \leq \sigma^{2} \), it holds
\begin{EQA}
	\dimA_{\imav}
	& = &
	\dmax^{4} \tr \Var(\errZv) 
	\leq 
	\dmax^{4} \sigma^{2} \, \dimq \, .
\label{uc7te737tvdsge5guru}
\end{EQA}
If the operator noise \( \errKS = \muA^{2} (\hKS - \E \hKS) \) is also homogeneous with 
\( \Var(\muA \hKS_{mj}) \leq \sigma_{\KS}^{2} \) then
\begin{EQA}[c]
	\dimA_{\KS}
	\leq 
	\tr \Var(\errKS)
	\leq 
	\sigma_{\KS}^{2} \, \dimp \, \dimq 
	\, .
\end{EQA}
The operator penalization can be used to reduce the value \( \dimA_{\KS} \):
\begin{EQA}
	\Var\bigl\{ (\dmax^{-2} \Id_{\KS} + \muA^{-2} \GPKS^{2})^{-1} \errKS \bigr\}
	&=&
	(\dmax^{-2} \Id_{\KS} + \muA^{-2} \GPKS^{2})^{-1} \Var( \errKS ) 
	(\dmax^{-2} \Id_{\KS} + \muA^{-2} \GPKS^{2})^{-1} 
	\, .
\end{EQA}
Let also \( \GPKS^{2} = \blk\{ \GPKS_{\mm}^{2} \} \).
Then the bound \( \| \Var( \errKS_{\mm}) \| \leq \CONSTi_{\errKSe} \) yields
\begin{EQA}
	\tr \Var\bigl\{ (\dmax^{-2} \Id_{\KS} + \muA^{-2} \GPKS^{2})^{-1} \errKS \bigr\}
	& \leq &
	\CONSTi_{\errKSe} \sum_{\mm=1}^{\dimq} 
	\tr (\dmax^{-2} \Id_{\dimp} + \muA^{-2} \GPKS_{\mm}^{2})^{-2} 
	\, .
\label{4bdce6w3nhdcf73jfyrhe}
\end{EQA}

\Subsection{Proof of Proposition~\ref{Leffscore}}
Now we study the stochastic term \( \bigl( \IFT_{\GPT}^{-1} \nabla \zeta \bigr)_{\tarpv} \)
for \( \nabla \zeta = (0,\nabla_{\ima},\nabla \KS) \).
Let \( \nupv = (\imav,\KS) \) be the aggregated nuisance parameter and 
\begin{EQA}
	\scorem{\nupv} \zeta 
	&=& 
	\binom{\errZv}{\errKS} \, ,
	\qquad
	\IFT_{\GPT,\nupv\nupv}
	=
	\begin{pmatrix}
		\IFT_{\imav\imav} & \IFT_{\imav\KS} \\
		\IFT_{\KS\imav} & \IFT_{\KS\KS} + \GPKS^{2}
	\end{pmatrix}
\label{sychhfejdbjwjttgvwcj}
\end{EQA}
Representation \eqref{8jdkvtwfe6xhejdcedvscy} yields by \( \nabla_{\tarpv} \zeta = 0 \)
\begin{EQA}
	\bigl( \IFT_{\GPT}^{-1} \nabla \zeta \bigr)_{\tarpv}
	&=&
	- \IFTb_{\GPT,\tarpv\tarpv}^{-1} \IFT_{\tarpv\nupv} \, \IFT_{\GPT,\nupv\nupv}^{-1} \scorem{\nupv} \zeta
	\, .
\label{yydult6guy76734esako}
\end{EQA}
Here \( \IFTb_{\GPT,\tarpv\tarpv}^{-1} \) is the \( \tarpv\tarpv \)-block of
\( \IFT_{\GPT}^{-1} \).
By Lemma~\ref{LIFTDFMEO}, \( \IFT_{\GPT} \geq \dmax^{-2} \DFM^{2} + \GPT^{2} \) 
with \( \DFM^{2} = \blk\{ \DPN^{2}, \Id_{\dimq}, \muA^{2} \Id_{\KS} \} \)
and \( \GPT^{2} = \blk\{ \GP^{2}, 0 , \GPKS^{2} \} \)
yielding
\begin{EQA}
	\IFTb_{\GPT,\tarpv\tarpv}^{-1}
	& \leq &
	(\dmax^{-2} \DPN^{2} + \GP^{2})^{-1} \, ,
	\qquad
	\IFT_{\GPT,\nupv\nupv}^{-1}
	\leq 
	\blk\bigl\{ \Id_{\dimq} \, , (\muA^{2} \Id_{\KS} + \GPKS^{2})^{-1} \bigr\} \, .
	\qquad
\label{90dr8f76yh3fyywigw76}
\end{EQA}
Further, \( \IFT_{\tarpv\nupv} \) consists of two blocks \( \IFT_{\tarpv\imav} \) and \( \IFT_{\tarpv\KS} \). 
It holds by \eqref{90dr8f76yh3fyywigw76}
\begin{EQA}
	&& \hspace{-20pt}
	\Var \bigl( \IFT_{\GPT}^{-1} \nabla \zeta \bigr)_{\tarpv}
	= 
	\IFTb_{\GPT,\tarpv\tarpv}^{-1} \, \IFT_{\tarpv\nupv} \, \IFT_{\GPT,\nupv\nupv}^{-1} \, \Var(\nabla_{\nupv} \zeta)
	\IFT_{\GPT,\nupv\nupv}^{-1} \, \IFT_{\nupv\tarpv} \, \IFTb_{\GPT,\tarpv\tarpv}^{-1} 
	\\
	& \leq &
	2 \IFTb_{\GPT,\tarpv\tarpv}^{-1} \, \IFT_{\tarpv\nupv} \, \IFT_{\GPT,\nupv\nupv}^{-1} \, 
	\blk\{ \Var(\errZv)\, , \Var(\muA \, \errKS) \} \, 
	\IFT_{\GPT,\nupv\nupv}^{-1} \, \IFT_{\nupv\tarpv} \, \IFTb_{\GPT,\tarpv\tarpv}^{-1}
	\\
	& = &
	2 \IFTb_{\GPT,\tarpv\tarpv}^{-1} \Bigl\{ 
		\IFT_{\tarpv\imav} \Var(\errZv) \IFT_{\imav\tarpv}
		+ 
		\IFT_{\tarpv\KS} \, 
		\Var\bigl\{ (\muA^{2} \Id_{\KS} + \GPKS^{2})^{-1} \muA \, \errKS \bigr\} \IFT_{\KS\tarpv} 
	\Bigr\} \IFTb_{\GPT,\tarpv\tarpv}^{-1}
	\, .
\label{ytwyuwswf5r4few3rfqikvwt}
\end{EQA}
By \eqref{8dfnjey63f6ytjutzv}, \( \IFT_{\tarpv\imav} =	- {\KSs}^{\T} \) and
\( \Var(\errZv) \leq \CONSTi_{\errZ} \Id_{\dimq} \) implies
\begin{EQA}[c]
	\IFT_{\tarpv\imav} \Var(\errZv) \IFT_{\imav\tarpv} 
	\leq 
	\CONSTi_{\errZ} {\KSs}^{\T} \KSs 
	\leq 
	\CONSTi_{\errZ} \, \DPN^{2} .
\end{EQA}
Further, for any \( \dtarpv \in \R^{\dimp} \) and \( \dKS \in \R^{\dimq \times \dimp} \), the use of \( {\KSs}^{\T} \tarpvs - \imavs = 0 \) yields
\begin{EQA}[c]
	\bigl\langle \IFT_{\tarpv\KS}, \dtarpv \otimes \dKS \bigr\rangle
	=
	- 2 \bigl\langle \KSs \dtarpv, \dKS \tarpvs \bigr\rangle 
	\, 
\end{EQA}
and by \( \| \DPN \thetavs \| \leq \REO \)
\begin{EQA}[c]
	\bigl| \bigl\langle \KSs \dtarpv, \dKS \tarpvs \bigr\rangle \bigr|
	\leq 
	\| \KSs \dtarpv \| \; \| \dKS \|_{\Fr} \, \nEO^{-1/2} \| \DPN \tarpvs \|
	\leq 
	\| \DPN \dtarpv \| \, \muA \, \| \dKS \|_{\Fr} \, \frac{\REO}{\muA \nEO^{1/2}} 
	=
	\muA \, \eEO \| \DPN \dtarpv \| \; \| \dKS \|
	\, .
\end{EQA}
Therefore, \( \muA^{-1} \bigl\| \DPN^{-1} \IFT_{\tarpv\KS} \bigr\| \leq 2 \eEO \) and
by \eqref{x6dc4c32c2ctv9necgsdr}
\begin{EQA}[rcl]
	\IFT_{\tarpv\KS} \, 
	\Var\bigl\{ (\muA^{2} \Id_{\KS} + \GPKS^{2})^{-1} \muA \, \errKS \bigr\} \IFT_{\KS\tarpv}
	&=&
	\muA^{-2} \IFT_{\tarpv\KS} \, 
	\Var\bigl\{ (\Id_{\KS} + \muA^{-2} \GPKS^{2})^{-1} \errKS \bigr\} \IFT_{\KS\tarpv}
	\\
	& \leq &
	\CONSTi_{\errKSe} \muA^{-2} \IFT_{\tarpv\KS} \, \IFT_{\KS\tarpv}
	\leq 
	\CONSTi_{\errKSe} \, \eEO^{2} \, \DPN^{2} 
	\, .
\end{EQA}
This implies
\begin{EQA}[rcl]
	\Var \bigl( \IFT_{\GPT}^{-1} \nabla \zeta \bigr)_{\tarpv}
	& \leq &
	2 (\CONSTi_{\errZ} + 4 \CONSTi_{\errKSe} \eEO^{2}) \, 
	\IFTb_{\GPT,\tarpv\tarpv}^{-1} \, \DPN^{2} \, \IFTb_{\GPT,\tarpv\tarpv}^{-1}
	\, .
\end{EQA}
By \eqref{90dr8f76yh3fyywigw76}, 
\( \IFTb_{\GPT,\tarpv\tarpv} \geq \DPN_{\dmax}^{2} \eqdef \dmax^{-2} \DPN^{2} + \GP^{2} \).
As \( \dmax^{-2} \DPN^{2} \leq \DPN_{\dmax}^{2} \), it holds
\begin{EQA}[rcl]
	\bigl\| (\DPN_{\dmax} \, \IFTb_{\GPT,\tarpv\tarpv}^{-1} \, \DPN_{\dmax}) \, (\DPN_{\dmax}^{-1} \, \DPN^{2} \, \DPN_{\dmax}^{-1}) \, 
	(\DPN_{\dmax} \, \IFTb_{\GPT,\tarpv\tarpv}^{-1} \, \DPN_{\dmax}) \bigr\|
	& \leq &
	\dmax^{2} 
\end{EQA}
and \( \IFTb_{\GPT,\tarpv\tarpv}^{-1} \, \DPN^{2} \, \IFTb_{\GPT,\tarpv\tarpv}^{-1} \leq \dmax^{2} \DPN_{\dmax}^{-2} \).
This yields \eqref{du7dyfu3ijhrk398rffi33} and \eqref{judfvy6h34ibte35ugit}.
Similarly for \( \xivr_{\GPT} \)
\begin{EQA}[rcl]
	\E \| \xivr_{\GPT} \|^{2}
	& \leq &
	2 (\CONSTi_{\errZ} + 4 \CONSTi_{\errKSe} \eEO^{2}) \, 
	\tr\bigl( \IFTb_{\GPT,\tarpv\tarpv}^{-1/2} \, \DPN^{2} \, \IFTb_{\GPT,\tarpv\tarpv}^{-1/2} \bigr)
	\leq 
	2 \dmax^{2} (\CONSTi_{\errZ} + 4 \CONSTi_{\errKSe} \eEO^{2}) \, 
	\tr (\DPN_{\dmax}^{-2} \, \DPN^{2})
	\, ,
\end{EQA}
and the second statement follows as well.

\Subsection{Proof of Lemma~\ref{PsemiriskEO}}

With \( \nupv = (\imav,\KS) \),
representation \eqref{8jdkvtwfe6xhejdcedvscy} leads to
\begin{EQA}
	(\IFT_{\GPT}^{-1} \GPT^{2} \prmtvs)_{\tarpv} 
	& = &
	\IFTb_{\GPT,\tarpv\tarpv}^{-1} 
	\bigl( \GP^{2} \tarpvs - \IFT_{\tarpv\nupv} \, \IFT_{\GPT,\nupv\nupv}^{-1} \, \GPY^{2} \nupvs \bigr) 
	\, ,
\label{vi458gghe3hwjdxqqwwEOr}
\end{EQA}
where \( \GPY^{2} \nupvs = (0,\GPKS^{2} \KSs) \) for \( \GPKS^{2} = \blk\{ \GPKS_{\mm}^{2} \} \).
As \( \IFT_{\imav\imav} = 2 \Id_{\dimq} \), it follows from
\eqref{8jdkvtwfe6xhejdcedvscye} of Lemma~\ref{LSchur} and \eqref{8dfnjey63f6ytjutz} of Lemma~\ref{Lnabl2EO} that
\begin{EQA}
	\IFT_{\GPT,\nupv\nupv}^{-1} \binom{0}{\GPKS^{2} \KSs}
	&=&
	\begin{pmatrix}
		- \IFT_{\imav\imav}^{-1} \, \IFT_{\imav\KS} \, \IFTb_{\GPKS,\KS\KS}^{-1} \, \GPKS^{2} \KSs  \\
		\IFTb_{\GPKS,\KS\KS}^{-1} \, \GPKS^{2} \KSs
	\end{pmatrix}
	=
	\begin{pmatrix}
		- \frac{1}{2} \, \IFT_{\imav\KS} \, \IFTb_{\GPKS,\KS\KS}^{-1} \, \GPKS^{2} \KSs  \\
		\IFTb_{\GPKS,\KS\KS}^{-1} \, \GPKS^{2} \KSs
	\end{pmatrix} \, ,
	\qquad
\label{ihcie83jec948f4yicje}
\end{EQA}
where \( \IFTb_{\GPKS,\KS\KS}^{-1} \) is the \( \KS\KS \)-block of \( \IFT_{\GPT,\nupv\nupv}^{-1} \):
\begin{EQA}[c]
	\IFTb_{\GPKS,\KS\KS} = \IFT_{\KS\KS} + \GPKS^{2} 
	- \frac{1}{2} \IFT_{\KS\imav} \IFT_{\imav\KS} \, .
\end{EQA}
By \eqref{7mdfgurt564634yufrt6re} and \( \imavs_{\mm} - {\KSs_{\mm}}^{\T} \tarpvs \equiv 0 \), 
the matrix \( \IFT_{\KS\imav} \IFT_{\imav\KS} \) is block-diagonal with 
the blocks \( \tarpvs {\tarpvs}^{\T} \).
This implies
\begin{EQA}[c]
	\IFTb_{\GPKS,\KS_{\mm}\KS_{\mm}} 
	=
	\muA^{2} \Id_{\dimp} + \GPKS_{\mm}^{2} - \frac{1}{2} \tarpvs {\tarpvs}^{\T} 
	\, , \qquad
	\mm=1,\ldots,\dimq .
\end{EQA}
Similarly, by \eqref{ihcie83jec948f4yicje} and \eqref{8dfnjey63f6ytjutzv}, \eqref{7mdfgurt564634yufrt6re} of Lemma~\ref{Lnabl2EO}
\begin{EQA}
	&& \nquad
	\IFT_{\tarpv\nupv} \, \IFT_{\GPT,\nupv\nupv}^{-1} \GPY^{2} \nupvs
	=
	\Bigl( - \frac{1}{2} \IFT_{\tarpv\imav} \IFT_{\imav\KS} + \IFT_{\tarpv\KS} \Bigr) 
	\IFTb_{\GPKS,\KS\KS}^{-1} \, \GPKS^{2} \KSs
	\\
	&=&
	\frac{1}{2} \sum_{\mm=1}^{\dimq} \KSs_{\mm} {\tarpvs}^{\T} 
	\bigl( \muA^{2} \Id_{\dimp} + \GPKS_{\mm}^{2} - \frac{1}{2} \tarpvs {\tarpvs}^{\T} \bigr)^{-1} \GPKS_{\mm}^{2} \KSs_{\mm}
	\\
	&=&
	\frac{1}{2} \sum_{\mm=1}^{\dimq} \KSs_{\mm} {\KSs_{\mm}}^{\T} \GPKS_{\mm}^{2} 
	\bigl( \muA^{2} \Id_{\dimp} + \GPKS_{\mm}^{2} - \frac{1}{2} \tarpvs {\tarpvs}^{\T} \bigr)^{-1} \tarpvs
	=
	S_{\GPKS} \tarpvs 
	\, .
\label{uedy8d38dhvkkrto5gj9e}
\end{EQA}
If \( \GPKS_{\mm}^{2} = \gpks_{\mm}^{2} \Id_{\dimp} \) then 
\begin{EQA}
	&& \nquad
	\IFT_{\tarpv\nupv} \, \IFT_{\GPT,\nupv\nupv}^{-1} \GPY^{2} \nupvs
	=
	\frac{1}{2} \sum_{\mm=1}^{\dimq} \KSs_{\mm} {\KSs_{\mm}}^{\T} \gpks_{\mm}^{2} 
	\Bigl\{ (\muA^{2} + \gpks_{\mm}^{2}) \Id_{\dimp} - \frac{1}{2} \tarpvs {\tarpvs}^{\T} 
	\Bigr\}^{-1} \tarpvs
	\\
	&=&
	\frac{1}{2} 
	\sum_{\mm=1}^{\dimq} \frac{\gpks_{\mm}^{2}}{\muA^{2} + \gpks_{\mm}^{2} - \| \tarpvs \|^{2}/2} \,
	\KSs_{\mm} {\KSs_{\mm}}^{\T} \, \tarpvs
	\, .
\label{fwydwyd3wt72te72yfywfcyef}
\end{EQA}
This completes the proof.


}

\end{document}